\documentclass[bj,authoryear,noshowframe]{imsart}

\usepackage{mathrsfs}
\usepackage{dsfont}
\usepackage{mathtools}
\usepackage{multicol}
\usepackage{subcaption}
\usepackage{graphicx}

\startlocaldefs

\theoremstyle{plain}
\newtheorem{theorem}{Theorem}[section]
\newtheorem{lemma}{Lemma}[section]
\newtheorem{remark}{Remark}[section]

\newcommand{\ZZ}{\mathcal{Z}}
\newcommand{\Z}{\mathbb{Z}}
\newcommand{\R}{\mathbb{R}}
\newcommand{\D}{\mathcal{D}}
\newcommand{\HH}{\mathcal{K}}
\newcommand{\Hc}{\mathcal{H}}
\newcommand{\B}{\mathcal{B}}

\newcommand{\U}{\mathcal{U}}
\newcommand{\V}{\mathcal{V}}
\newcommand{\N}{\mathbb{N}}

\newcommand{\M}{\mathcal{M}}
\newcommand{\C}{\mathcal{C}}
\newcommand{\Y}{\mathcal{Y}}
\newcommand{\X}{\mathcal{X}}
\newcommand{\PP}{\mathcal{P}}

\numberwithin{equation}{section}
\allowdisplaybreaks

\usepackage{xcite}
\usepackage{xr}
\makeatletter
\newcommand*{\addFileDependency}[1]{
  \typeout{(#1)}
  \@addtofilelist{#1}
  \IfFileExists{#1}{}{\typeout{No file #1.}}
}
\makeatother

\begin{document}

\begin{frontmatter}
\title{A nonparametric distribution-free test of independence among continuous random vectors based on \texorpdfstring{$L_1$}{}-norm}
\runtitle{A nonparametric distribution-free test of independence}

\begin{aug}
\author[A]{\inits{B.}\fnms{Nour-Eddine}~\snm{Berrahou}\ead[label=e1]{n.berrahou@uca.ac.ma}}
\author[B]{\inits{S.}\fnms{Salim}~\snm{Bouzebda}\ead[label=e2]{salim.bouzebda@utc.fr}\orcid{0000-0001-7801-4945}}
\author[A]{\inits{L.}\fnms{Lahcen}~\snm{Douge}\ead[label=e3]{lahcen.douge@uca.ac.ma}}
\address[A]{
Laboratoire de Modélisation des Systèmes Complexes, University Cadi Ayyad\\ Marrakesh, Morocco\printead[presep={,\ }]{e1,e3}}

\address[B]{Universit\'e de technologie de Compi\`egne, LMAC (Laboratory of Applied Mathematics of Compi\`egne)\printead[presep={,\ }]{e2}}
\end{aug}

\begin{abstract}
\noindent We propose a novel statistical test to assess the mutual independence of multidimensional random vectors. Our approach is based on the $L_1$-distance between the joint density function and the product of the marginal densities associated with the presumed independent vectors. Under the null hypothesis, we employ Poissonization techniques to establish the asymptotic normal approximation of the corresponding test statistic, without imposing any regularity assumptions on the underlying Lebesgue density function, denoted as $f(\cdot)$. Remarkably, we observe that the limiting distribution of the $L_1$-based statistics remains unaffected by the specific form of $f(\cdot)$. This unexpected outcome contributes to the robustness and versatility of our method. 
Moreover, our tests exhibit nontrivial local power against a subset of local alternatives, which converge to the null hypothesis at a rate of {${\tiny  n^{\tiny -1/2}h_n^{\tiny -{d/4}}}$}, $d\geq 2$, where $n$ represents the sample size and $h_n$ denotes the bandwidth. Finally, the theory is supported by a comprehensive simulation study to investigate the finite-sample performance of our proposed test. The results demonstrate that our testing procedure generally outperforms existing approaches across various examined scenarios.

\end{abstract} 
\begin{keyword} 
\kwd{Asymptotic normality}
\kwd{distribution-free tests}
\kwd{independence test}
\kwd{kernel density function estimator}
\kwd{$L_1$-distance}
\kwd{Poissonization}.\\
\noindent\href{https://mathscinet.ams.org/mathscinet/msc/msc2020.html}{MSC2020 subject classifications}: 60F05; (60F15, 60F17, 62G07).
\end{keyword}
\end{frontmatter}

\section{Introduction}

One of the classical and important problems in statistics is testing the independence between two or more components of a random vector.   Testing for mutual independence, which characterizes the structural relationships between random variables and is strictly stronger than pairwise independence, is a fundamental task in inference. Independent component analysis of a random vector, which consists of searching for a linear transformation that minimizes the statistical dependence between its components, is one of the fields in which mutual independence plays a central role; for instance, see
\cite{Bach}, \cite{chen}, \cite{Samworth} and \cite{matt}. Testing independence also has many applications, including causal inference (\cite{pearl}, \cite{peters},
\cite{pfister2018}, \cite{chakraborty2019}), graphical modeling (\cite{Lauritzen}, \cite{gan}), linguistics (\cite{nguyen}), clustering \citep{MR2231170}, dimension reduction \citep{MR2247974,MR3474038}. The traditional approach for testing independence is based on Pearson's correlation coefficient; for instance, refer to \cite{Binet}, \cite{PEARSONKARL}, \cite{spearman1904}, \cite{kendall1938}. However, its lack of robustness to outliers and departures from normality eventually led researchers to consider alternative nonparametric procedures.
To overcome such a problem, a natural approach is to consider the functional difference between
the empirical joint distribution and the product of the empirical marginal distributions, see \cite{Hoeffding}, \cite{Blum} and \cite{Bouzebda2011}.
This approach can also use characteristic empirical functions; see \cite{csorgHo1985}. Inspired by the work of \cite{Blum} and \cite{dugue1975}, \cite{Deheuvels1981} studied a test of multivariate independence based on the {M}\"{o}bius decomposition, generalized in \cite{Bouzebda2014}. 
 \cite{rosenblatt1975}, \cite{MR1190260} and \cite{Ahmad1997} presented the application of nonparametric kernel estimation method to evaluate the independence of two random vectors. This method involves comparing the $L_2$ norm of the difference between the joint density and the product of their marginals. 
The work conducted by \cite{han2017} delves into two families of distribution-free test statistics, notably featuring Kendall's tau (a rank-type $U$-statistic) and Spearman's rho (a simple linear rank statistic), both serving as significant examples. \cite{leung2018} addressed the problem of testing the independence of $m$ continuous variables when $m$ is greater than the sample size $n$. The authors considered three types of test statistics constructed as the sums or sums of squares of pairwise rank correlations.
\cite{MR4185806} investigated the testing of the mutual independence of a $p$-dimensional random vector, when $p=p_n$ increases to infinity. The authors derive maximum-type tests based on pairwise rank correlation measures. The rank correlations are presented in an elegant $U$-statistic framework and then aggregated to create maximum-type test statistics with a null distribution that converges to a nonstandard Gumbel distribution.
Recent studies have focused on testing random vectors' pairwise independence. The concept of distance covariance was introduced by \cite{szekely2007} to quantify the dependence between two random vectors, and \cite{szekely2009} investigated it further.
This distance covariance is defined as a weighted $L_2$-norm between the joint distribution's characteristic function and the product of the marginal distribution's characteristic functions. To test the independence of two random vectors, \cite{gretton2005} proposed a kernel-based independence criterion called the Hilbert-Schmidt independence criterion (HSIC). 
Due to the ``curse of dimensionality,'' these tests may fail to detect the non-linear dependence when the random vectors are of high dimension. To address this issue, \cite{MR4544605} proposed a general framework for testing the dependence of two random vectors by randomly selecting two subspaces composed of vector components.
\cite{MR4527082} proposed a new class of independence measures based on the maximum mean discrepancy in Reproducing Kernel Hilbert Space. In the literature, additional methods for testing the independence of two multidimensional random variables have emerged, including those based on the $L_1$-norm between the distribution of the vector and the product of the distributions of its components (\cite{gretton2010}), on ranks of distances \citep{heller2012}, on nearest neighbor methods~\citep{berrett} and on applying distance covariance to center-outward ranks and signs \citep{shi}. In \cite{MR3858367},  three distinct measures of mutual dependency are proposed.  One of the approaches extends the concept of distance covariance from pairwise dependency to mutual dependence, and the remaining two measures are derived by summing the squared distance covariances. Finally, 
\cite{auddy2023exact} considered the problem of testing for independence using the correlation from \cite{MR4353729}.  

A common idea to measure the difference between two functions is to consider the $L_q$-norm, $q\geq 1$. \cite{csorgHo1988} was the first to prove a central limit for the $L_q$-norm, $q\geq 1$, between a density function and its estimator, and \cite{horvath1991} introduced a Poissonization technique into the study of central limit theorems for the $L_q$-norm. It is noteworthy that the key Poissonization techniques date back to \cite{MR0410832} and \cite{MR576407}. 
\cite{mason1995} developed a general method for deriving the asymptotic normality of the $L_q$-norm of empirical functionals without assuming restrictive regularity conditions. 
\cite{mason2001} established the weak convergence of an $L_1$-norm density estimator process indexed by kernel functions.

This work presents a different approach to the problem of testing the mutual independence of random vectors in arbitrary dimensions based on the $L_1$-norm between the joint density and the product of $p$ the marginal densities, for $p\geq2$. The motivation for $L_1$-norm is connected with the fact that $\int_{\mathbb R^d}|f(\tilde {\mathbf x})-g(\tilde{\mathbf x})| d \tilde {\mathbf x}$ is the total variation distance between the measures with densities $f(\cdot)$ and $g(\cdot)$. As a result, this criterion does not share any good qualities associated with the $L_2$ (see Remark \ref{remark4.1} below) or $L_\infty$ norms, both of which are mathematically simpler to work with; for instance, see \cite{Devroye12} and \cite{Scot1992} for more details. 
\cite{MR3513521} pointed out that the choice of $L_1$-norm for the comparison of curves poses several mathematical challenges, which are caused by the fact that the mapping $f \rightarrow \int_{\mathbb R^d}|f(\tilde {\mathbf x})| d \tilde {\mathbf x}$ is in general not (Hadamard-) differentiable. 
Since the Lebesgue density functions, by definition, sit in $L_1(\mathbb{R}^d, \mathcal{B}, \lambda)$, where $\lambda$ denotes the Lebesgue measure and  $\mathcal{B}$ is the class of all Borel sets of $\mathbb R^d$,  \cite{devroye1} has repeatedly stated that the $L_1$-norm is the most appropriate metric for characterizing the error in estimation between a density function $f(\cdot)$ and its estimator.

To our best knowledge, this is the first time that the general context $L_1$-norm for testing the independence appeared in the literature and gives the main motivation of the present work by responding to the open problems mentioned in \cite{gretton2010}. The main contribution of this paper is to establish the asymptotic distribution of the proposed test statistic under the null hypothesis and under the local alternatives that converge to the null at the rate of $n^{-1/2}h_n^{-d/4}$. As an important feature, the $ L_1$-based 
tests are all model-free. Furthermore, no regularity requirements for the densities are necessary to demonstrate the asymptotic normality of our statistic, a desired attribute. We conduct simulations to determine the size and power of the test. We illustrate that the proposed test has superior power characteristics compared to existing tests for various analyzed situations.  The proposed test encompasses all dependency types, including complicated dependence structures, particularly sinusoidal dependence. 
As it will be seen later, the problem requires much more than ‘simply’ combining ideas from existing papers. Delicate mathematical derivations will be required to cope with $L_1$-based tests. We highlight that the use of the Poissonization techniques due \cite{mason2001} is substantially more difficult in our context presenting the technical merits of our contribution, see Remark \ref{comparison_Mason} below. More precisely, several problems arise when the $L_1$-norm is used in nonparametric tests. Even when the sample consists of independent random vectors, the $L_1$-norm introduces a dependence structure. This is the fundamental challenge. 
The solution to these problems is broken down into three distinct steps: Truncation, Poissonization, and dePoissonization devices are appropriate tools to reach the results. Particularly, Poissonization involves randomizing the sample size using a Poisson random variable, allowing the application of techniques requiring the independence of increments.

The layout of the article is as follows.
Section \ref{notation} gives the notation and the definitions we need.
Section \ref{Main_results} provides the asymptotic behavior of the $L_1$-based test under the null hypothesis and gives the consistency against the alternative at the nominal level $\alpha$. 
Section \ref{localalterantive}
provides the limiting law of the $L_1$-based test under the contiguous hypothesis. 
We discuss a bandwidth choice for practical use in Section \ref{Bandwidth_Selection}.
Section \ref{simulation} summarizes comprehensive simulation tests to evaluate the finite-sample properties of our approach across various dependency scenarios and compare it
to several state-of-the-art approaches. A summary of the findings highlighting remaining open issues are given in Section \ref{conclusion}. All proofs are deferred to Section \ref{proofs}. Due to the lengthiness of the proofs,  we restrict our analysis to the most important arguments.  Finally, technical proofs can be found in the supplementary material \cite{L1testBBD}.

\section{Notation and setting} \label{notation}
\noindent We want to test the mutual independence of $p$ continuous random vectors {$\mathbb{X}^1, \ldots, \mathbb{X}^p$} of dimensions $d_1, \ldots, d_p$ respectively, based on the $L_1$-distance of the discrepancy between the product of their densities and the joint density. Let $d=d_1+\cdots+d_p$ be the dimension of the random vector {$(\mathbb{X}^{1\top}, \ldots, \mathbb{X}^{p\top})^\top$}, $p\geq 1$, where $\mathbf x^{\top}$ denotes the transpose of $\mathbf x$. 
Accordingly,
$\tilde{\mathbf x}=(\mathbf x_1^{\top},\ldots,\mathbf x_p^{\top})^{\top}\in \R^{d_1}\times\cdots \times\R^{d_{p}}$ is a vector of length $d$.
Let $f(\cdot)$ denote the density function of $\tilde{\mathbf X}=(X_1,\ldots,X_d)^\top=(\mathbb{X}^{1\top}, \ldots, \mathbb{X}^{p\top})^\top$ with respect to the Lebesgue measure $\lambda$ on $\R^{d}$. For each $l=1,\ldots,p$, let $f_l(\cdot)$ be the density function of $\mathbb X^{l}$ with respect to the Lebesgue measure $\lambda_l$ on $\R^{d_l}$. Let us formulate the null hypothesis 
\begin{eqnarray*}
	\Hc_0 : \int_{\R^{d}}\bigg|f(\tilde{\mathbf x}) - \prod_{l=1}^{p}f_{l}(\mathbf x_{l})\bigg|\,d\tilde{\mathbf x} = 0,
~~
\mbox{against the alternative}
~~
 \Hc_1 : \int_{\R^{d}}\bigg|f(\tilde{\mathbf x}) - \prod_{l=1}^{p}f_{l}(\mathbf x_{l})\bigg|\,d\tilde{\mathbf x} > 0.
\end{eqnarray*}
We next consider sequences $\left\{\tilde{\mathbf X}_i: i \geq 1\right\}$ and $\left\{\mathbb{X}_i^l: i \geq 1\right\}$, $l=1,\ldots,p$, of independent and identically distributed [iid] random copies of the random vector [rv] $\tilde{\mathbf X}$ and $\mathbb{X}^l$, $l=1,\ldots,p$, respectively. 
We now introduce  the kernel functions $K_l(\cdot)$ defined on $\mathbb{R}^{d_l}$,  $l=1,\ldots,p$; fulfilling the condition ({\bf A$_1$}) below. For each $n \geq 1$, and for each choice of the bandwidth $h_n > 0$, we define 
the classical
Akaike–Parzen–Rosenblatt [refer to \cite{Akaike1954,Parzen1962,Rosenblatt}] kernel estimators of the density functions $f(\cdot)$ and $f_l(\cdot)$, $l=1,\ldots,p$, respectively, by
\begin{eqnarray*}
	f_{n}(\tilde{\mathbf x})&=&f_{n,h_n}(\tilde{\mathbf x})=
	\frac{1}{nh_{n}^{d}}
	\sum\limits_{i=1}^{n}\prod_{l=1}^{p}
	K_{l}\Bigg(\frac{\mathbf x_l-\mathbb X^{l}_{i}}{h_{n}}\Bigg),~~~
	f_{n,l}(\mathbf x_l)=f_{n,h_n,l}(\mathbf x_l)=\frac{1}{nh_{n}^{d_{l}}}\sum\limits_{i=1}^{n}
	K_{l}\Bigg(\frac{\mathbf x_l-\mathbb X^{l}_{i}}{h_{n}}\Bigg).
\end{eqnarray*} 
\cite{Parzen1962} has shown, under some assumptions on the kernel function, that $f_n(\cdot)$ is an asymptotically unbiased and consistent estimator for $f(\cdot)$ whenever $h_n \rightarrow 0, n h_n^d \rightarrow \infty$ and $\tilde{\mathbf{x}}$ is a continuity point of $f(\cdot)$. Under some additional assumptions on $f(\cdot)$ and $h_n$, he obtained an asymptotic normality result, too. 
The above Akaike–Parzen–Rosenblatt kernel estimators have been extensively studied in the literature, see, e.g., 
\cite{silverman1986}, 
\cite{nadaraya1989}, \cite{hardle1990},  \cite{Scot1992},
 \cite{wand1995}, \cite{Eggermont2001}, 
\cite{devroye2001}, \cite{MR4562252}, \cite{taachoucheB} and the references therein. To perform the test of the hypothesis $\Hc_0$,
we consider the statistic
\begin{eqnarray}\label{statVn}
V_{n}:=\int_{\R^{d}} \bigg|f_{n}(\tilde{\mathbf x})-\prod_{l=1}^{p}f_{n,l}(\mathbf x_l)\bigg|\,d\tilde{\mathbf x}.
\end{eqnarray}
\begin{remark} \label{remark201}
	Careful bookkeeping throughout the proofs in the paper could give an adapted set of conditions on possibly different bandwidths. We do not pursue this here to avoid an unnecessary increase in the technicality of the presentation.  If one wants to use the vector bandwidths (see, in particular, Chapter 12 of \cite{devroye2001}), with obvious changes in notation,
	our results and their proofs remain true when $h_{n}$ is replaced by a vector bandwidth $\mathbf{ h}_{n} = (h^{(1)}_{n}, \ldots, h^{(d)}_{n})$, where $\min_{1\leq i \leq d} h^{(i)}_{n} > 0$. In this situation we set 
	$h_{n}^d=\prod_{i=1}^{d} h_{n}^{(i)}$, and for any vector $\mathbf{ v} = (v_{1} ,\ldots,v_{d})$ we replace $\mathbf{ v}/h_{n}$ by 
	$(v_{1}/h_{n}^{(1)},\ldots,v_{d}/h_{n}^{(d)})$. 
\end{remark} 
\section{Main results} \label{Main_results}

\noindent Before presenting our results, we introduce some notation and the underlying assumptions. For $s\geq 1$, let $\B_{s}=\{\mathbf u\in {\R}^{s}:
\| \mathbf u\|\leq 1\}$, where $\| \cdot\|$ is some norm on ${\R}^{s}$.
Consider the following condition upon the kernels $K_l(\cdot), l=1,\ldots,p$, that will be
used below. 

\begin{enumerate}
	\item[{({\bf A$_1$})}]
	For each $l=1,\ldots,p$,
	$K_l(\cdot)$ is a bounded Lebesgue density function having support contained in the closed ball of radius $1/2$ centered at zero of $\R^{d_{l}}$.
\end{enumerate}
Let $Z_{1}$ and $Z_{2}$ be independent standard normal random variables. In the sequel, the $L_2$-norm of a kernel
$K(\cdot)$ is denoted by $\|K\|_2$. For each $l=1,\ldots,p$, define the function
\begin{eqnarray*}
	\rho_{l}(\mathbf t_l)=\frac{\displaystyle \int_{\R^{d_l}} K_{l}(\mathbf u_l) K_{l}(\mathbf u_l+\mathbf t_l)\,d\mathbf u_l}{\displaystyle \| K_{l}\|_{2}^2 }, \qquad \quad \mathbf t_l\in \R^{d_l}.
\end{eqnarray*}
Denote $\tilde{\bf{K}}=\prod_{l=1}^{p}\| K_l \|_2^2$ and, for each $\tilde{\mathbf{t}}=(\mathbf t_1^{\top},\ldots,\mathbf t_p^{\top})^{\top}$, $
\rho(\tilde{\mathbf t})=\prod_{l=1}^{p}\rho_l(\mathbf t_l)$.
Further, for $\tilde{\B}=\B_{d_1}\times \cdots\times \B_{d_p}$, set
\begin{eqnarray*}
	\sigma^2:=\tilde{\bf{K}} \int_{\tilde{\B}}\mathrm{Cov}\bigg(\bigg|
	\sqrt{1-\rho^2(\tilde{\mathbf t})}Z_1
	+\rho(\tilde{\mathbf t})Z_2\bigg|,|Z_2|\bigg)\,d\tilde{\mathbf t}. 
\end{eqnarray*}
Below, we write $Z \stackrel{\mathcal{D}}{=} \mathcal{N}(\mu, \sigma^{2} )$ whenever the random variable $Z$ follows a normal law with expectation $\mu$ and variance $\sigma^{2}$. Let ${\stackrel{\mathcal{D}}{\rightarrow}}$ denote the  convergence in distribution. The main result to be proved here may now be stated precisely as follows.
\begin{theorem}\label{thm1}
	Assume that \textnormal{({\bf A$_1$})} holds and let $h_n>0$  be such that, as $n \rightarrow \infty$, $h_{n}\rightarrow0$ and $nh_{n}^{3d}\rightarrow\infty$. 
	Then, under $\Hc_0$, as $n \rightarrow \infty$, we have
$
		\sqrt{n}\big(V_{n} -\mathbb E V_{n}\big)\stackrel{\D}{\rightarrow} \mathcal{N}(0,\sigma^2).$
\end{theorem}

\begin{remark}
	{{
				Note that the choice $h_n=\lambda n^{-\vartheta},~ \lambda>0, ~~0<\vartheta<1/(3d)$, fulfills the theorem's conditions, see Section \ref{Bandwidth_Selection} for more discussions.}
		\textnormal{ The variance $\sigma^2$ has an interesting alternate representation. Making use of the formulas for the absolute moments of a bivariate normal random variable of \cite{nabeya}, we have 
			\begin{eqnarray*}
				\sigma^2=\tilde{\bf{K}} \int_{\tilde{\B}}\varphi\big(\rho(\tilde{\mathbf t})\big)\,d\tilde{\mathbf t},
			~~
			\mbox{where}
			~~
				\varphi(\rho)=\frac{2}{\pi}\left(\rho \arcsin \rho+\sqrt{1-\rho^{2}}-1\right), ~\mbox{for}~ \rho\in \lbrack-1,1\rbrack.
		\end{eqnarray*}}}
\end{remark}

\begin{remark}
{
 It is noteworthy that Theorem \ref{thm1} imposes any assumptions on the Lebesgue density.  
		The bandwidth condition $nh_{n}^{3d}\rightarrow\infty$ is mainly due to the evaluation of the variance in Lemma \ref{lem5}. It allows us also to provide some asymptotic approximations in Lemma \ref{lem51} and to approximate 
		the centering term in Lemma \ref{Lemma12}. This condition is considered in nonparametric regression testing proposed by \cite{lee2013} to establish the asymptotic normality of the one-sided $L_q$-test statistic in testing functional inequalities, $1\leq q < \infty$.}
\end{remark}
\noindent Let us discuss how the process works so that our test can be carried out. Let, for each $l=1,\ldots,p$, $\C_l$ be a bounded Borel set such that $\mathbb P(\mathbb X^{l}\in \C_{l})>1-\varepsilon $ for some $\varepsilon $ small enough 
and set, for a Borel set $A\in\R^{d}$,  
\begin{eqnarray*}
	V_{n}(A):=\int_{A} \bigg|f_{n}(\tilde{\mathbf x})-\prod_{l=1}^{p}f_{n,l}(\mathbf x_l)\bigg|\, d\tilde{\mathbf x},
~~\mbox{ and }~~ V(A):=\int_{A}\bigg|f(\tilde{\mathbf x})-\prod_{l=1}^{p}f_{l}(\mathbf x_l)\bigg|\, d\tilde{\mathbf x}.
\end{eqnarray*}
We denote $\tilde{\C}=\C_1\times\cdots\times\C_p$.
Making use of the arguments developed in Lemma \ref{lem2} in \cite{L1testBBD}, and Tchebychev's inequality, it follows that 
\begin{eqnarray*}
	\sqrt{n}\left[V_n-V_n(\tilde{\C})-\mathbb E\big(V_n-V_n(\tilde{\C})\big)\right]
	\stackrel{\PP}{\rightarrow}0,
\end{eqnarray*}
where $\stackrel{\PP}{\rightarrow}$ stands for the convergence in probability.
Therefore, an application of Theorem \ref{thm1} gives
\begin{eqnarray*}\label{eqprop1}
	\sqrt{n}\left(V_{n}(\tilde{\C}) -\mathbb E V_{n}(\tilde{\C})\right)\stackrel{\D}{\rightarrow} \mathcal{N}(0,\sigma^2)\ \ \mbox{as}\ \ n\rightarrow\infty.
\end{eqnarray*}
The centering term $\sqrt{n}\mathbb E V_n(\tilde{\C})$ depends on an unknown distribution. To circumvent this problem, we estimate $\sqrt{n}\mathbb E V_n(\tilde{\C})$ by
\begin{eqnarray}\label{andefinition}
	a_n(\tilde{\C}):= \mathbb E|Z_1|\int_{\tilde{\C}}\sqrt{\hat{\mathcal{L}}_{n}(\tilde{\mathbf x})}\, d\tilde{\mathbf x},
\end{eqnarray}
where
\begin{eqnarray*}
	\hat{\mathcal{L}}_{n}(\tilde{\mathbf x}):=\prod_{l=1}^{p} \hat{v}_{n,l}(\mathbf x_l)-\sum_{l=1}^{p}\hat{v}_{n,l}(\mathbf x_l) \prod_{j\neq l}g_{n,j}(\mathbf x_j)+(p-1)\prod_{l=1}^{p}g_{n,l}(\mathbf x_l),
\end{eqnarray*}
where, for each $l=1,\ldots,p$,
\begin{eqnarray*}
	\hat{v}_{n,l}(\mathbf x_l):=\frac{1}{nh_{n}^{2d_l}}\sum_{i=1}^{n} K^{2}_{l}\left(\frac{\mathbf x_l-\mathbb X^{l}_{i}}{h_{n}}\right),
\mbox{~
and
}
	g_{n,l}(\mathbf x_l):=
	\frac{1}{n(n-1)h_n^{2d_l}}
	\sum_{i=1}^{n}\sum_{j\neq i}K_{l}\left(\frac{\mathbf x_l-\mathbb X^{l}_{i}}{h_{n}}\right)K_{l}\left(\frac{\mathbf x_l-\mathbb X^{l}_{j}}{h_{n}}\right).
\end{eqnarray*}
Note that $\mathbb{E}|{Z}_1|=\sqrt{2 / \pi} \approx 0.7978$.
Our suggested test statistic is 
\begin{eqnarray*}
	T_{n}(\tilde{\C}):=\frac{\sqrt{n}V_{n}(\tilde{\C})-a_{n}(\tilde{\C})}{\sigma}.
\end{eqnarray*}
Let $z_{1-\alpha}$ be the upper $\alpha$ quantile of the standard normal distribution. In this paper, we propose the following test:
$\mbox{ Reject } \Hc_{0}\mbox{ if and only if } T_{n}(\tilde{\C})>z_{1-\alpha}.$
Now, for each $ l=1, \ldots,p$, for each $\delta>0$, let $\C_{l}^{\delta}$ be the $\delta$-neighborhood of $\C_{l}$, i.e., $ \C_{l}^{\delta}=\{ \mathbf x_l+ \mathbf u_l: \mathbf x_l\in \C_l, \|\mathbf u_l\|\leq \delta\}.$ We need the following additional notation.  Set, for each $a\in\{0, 1\}$, 
\begin{eqnarray*}
	\Gamma_{n,a}(\tilde{\mathbf x}):= f_{n}(\tilde{\mathbf x})-
	a\sum_{l=1}^{p} f_{n,l}(\mathbf x_l)\prod_{j\neq l }\mathbb E f_{n,j}(\mathbf x_j)
	+(pa-1)\prod_{l=1}^{p}\mathbb E f_{n,l}(\mathbf x_l)
\end{eqnarray*}
and 
\begin{eqnarray*}
	\Gamma_{\eta,a}(\tilde{\mathbf x}):= f_{\eta}(\tilde{\mathbf x})-
	a\sum_{l=1}^{p} f_{\eta,l}(\mathbf x_l)\prod_{j\neq l }\mathbb E f_{n,j}(\mathbf x_j)
	+(pa-1)\prod_{l=1}^{p}\mathbb E f_{n,l}(\mathbf x_l),
\end{eqnarray*}
where
\begin{eqnarray*}
	f_{\eta}(\tilde{\mathbf x})=
	\frac{1}{nh_{n}^{d}}
	\sum\limits_{i=1}^{\eta}\prod_{l=1}^{p}
	K_{l}\Bigg(\frac{\mathbf x_l-\mathbf X^{l}_{i}}{h_{n}}\Bigg) \,\, \mbox{ and } \,\,
	f_{\eta,l}(\mathbf x_l)=\frac{1}{nh_{n}^{d_{l}}}\sum\limits_{i=1}^{\eta}
	K_{l}\Bigg(\frac{\mathbf x_l-\mathbf X^{l}_{i}}{h_{n}}\Bigg),
\end{eqnarray*}
and $\eta$ is a Poisson random variable with mean $n$, independent
of $\tilde{\mathbf X},\tilde{\mathbf X}_1,\tilde{\mathbf X}_2, \ldots$. Set, for a Borel set $A\in\R^{d}$ and for each $a\in\{0, 1\}$, 
\begin{eqnarray*}
	U_{n,a}(A):=\int_{A} |\Gamma_{n,a}(\tilde{\mathbf x})|\,d\tilde{\mathbf x} \,\, \mbox{ and } \,\, U_{\eta,a}(A):=\int_{A} |\Gamma_{\eta,a}(\tilde{\mathbf x})|\,d\tilde{\mathbf x}.
\end{eqnarray*}
The following theorem shows that the test has an asymptotically valid size.
\begin{theorem}\label{th2}
	Assume that \textnormal{({\bf A$_1$})} holds and, for each $ l=1, \ldots,p$, there exists $\delta>0$ such that $f_l(\cdot)$ is bounded on $ \C_{l}^{\delta}$. 
	Let $h_n>0$  be such that, as $n \rightarrow \infty$, $h_{n}\rightarrow0$ and $nh_{n}^{3d}\rightarrow\infty$. 
	Then, under $\Hc_0$, we have 
	$
		\lim_{n\rightarrow \infty} \mathbb P\{T_{n}(\tilde{\C})>z_{1-\alpha} \}= \alpha.$
\end{theorem}

\noindent The last result of this section shows the consistency of the test against the alternative at nominal level $\alpha$.
\begin{theorem}\label{th3}
	Assume that for each $l=1,\ldots,p$,
	$K_l(\cdot)$ is a density function.
	Let $h_n>0$  be such that, as $n \rightarrow \infty$, $h_{n}\rightarrow0$ and $nh_{n}^{3d}\rightarrow\infty$. 
	If $ \Hc_1$ is true, then, for any bounded subset $\tilde{D}=D_1\times\cdots\times D_p\subset \mathbb{R}^d$ such that $V(\tilde{D})>0$, we have
	$
		\lim_{n\rightarrow \infty} \mathbb P\{T_{n}(\tilde{D})>z_{1-\alpha} \}= 1.$
\end{theorem}

\begin{remark}
	Our results can be generalized by considering the following weighted family 
	\begin{eqnarray*}
		V_{n}^{(1)}:=\int_{\R^{d}} \bigg|f_{n}(\tilde{\mathbf x})-\prod_{l=1}^{p}f_{n,l}(\mathbf x_l)\bigg|w(\tilde{\mathbf x})\,d\tilde{\mathbf x},
	\end{eqnarray*}
	where $w: \mathbb{R}^d \rightarrow[0, \infty)$ is a nonnegative weight function, or, for $1\leq q<\infty$, 
	\begin{eqnarray}\label{statVnwei2}
		V_{n}^{(q)}:=\int_{\R^{d}} \bigg|f_{n}(\tilde{\mathbf x})-\prod_{l=1}^{p}f_{n,l}(\mathbf x_l)\bigg|^qw(\tilde{\mathbf x})\,d\tilde{\mathbf x}.
	\end{eqnarray}{
Let us recall a popular family of $\phi$-divergences
$$
D_\phi(\mathbb{P}, \mathbb{Q}):=\int_{\mathbb R^d} \phi\left(\frac{d \mathbb{P}}{d \mathbb{Q}}\right) d \mathbb{Q} \text { if } \mathbb{P} \ll \mathbb{Q} \text {, }
$$
where  $\phi:[0, \infty) \rightarrow(-\infty, \infty]$ is a convex function. $\mathbb{P} \ll \mathbb{Q}$ denotes that $\mathbb{P}$ is absolutely continuous with respect to $\mathbb{Q}$. Well-known distance/divergence measures obtained by appropriately choosing $\phi$ include total variation distance $(\phi(t)=$ $|t-1|)$. Let us define $\mathscr{P}_\lambda$ as the set of all probability measures, $\mathbb{P}$,  that are absolutely continuous with respect to some $\sigma$-finite measure, $\lambda$ on $\mathbb R^d$. For $\mathbb{P}, \mathbb{Q} \in \mathscr{P}_\lambda$, let $ \prod_{l=1}^{p}f_{l}=\frac{d \mathbb{P}}{d \lambda}$ and $f=\frac{d \mathbb{Q}}{d \lambda}$ be the Radon-Nikodym derivatives of $\mathbb{P}$ and $\mathbb{Q}$ with respect to $\lambda$, respectively. In the interesting paper \cite{MR2988458}, it is mentioned {that} $\phi$-divergences
and integral
probability metrics are fundamentally different and intersect only at the total variation
distance. In the last situation, we have 
$$
D^{TV}(\mathbb{P}, \mathbb{Q}):=\int_{\mathbb R^d} \left|\frac{d \mathbb{P}}{d \mathbb{Q}}-1\right| d \mathbb{Q}=\int_{\mathbb R^d} \left|\frac{\prod_{l=1}^{p}f_{l}}{f}-1\right| d \mathbb{Q}.
$$
Let $\mathbb{Q}_n$ denote the empirical measure based on $\{\tilde{\mathbf X}_i: i \geq 1\}$. One can estimate $D_\phi(\mathbb{P}, \mathbb{Q})$ by 
$$
D_{n}^{TV}(\mathbb{P}, \mathbb{Q}):=\int_{\mathbb R^d} \left|\frac{\prod_{l=1}^{p}f_{n,l,i}}{f_{n,i}}-1\right| d \mathbb{Q}_n=\frac{1}{n}\sum_{i=1}^{n}\left|\frac{\prod_{l=1}^p f_{n, l,i}\left(\mathbb{X}_i^l\right)}{f_{n,i}\left(\mathbb{X}_i^1, \ldots, \mathbb{X}_i^p\right)}-1\right|,
$$
where
\begin{eqnarray*}
	f_{n,i}(\tilde{\mathbf x})=
	\frac{1}{(n-1)h_{n}^{d}}
	\underset{j\neq i}{\sum\limits_{j=1}^{n}}\prod_{l=1}^{p}
	K_{l}\Bigg(\frac{\mathbf x_l-\mathbb X^{l}_{j}}{h_{n}}\Bigg),~~~
	f_{n,l,i}(\mathbf x_l)=\frac{1}{(n-1)h_{n}^{d_{l}}}\underset{j\neq i}{\sum\limits_{j=1}^{n}}
	K_{l}\Bigg(\frac{\mathbf x_l-\mathbb X^{l}_{j}}{h_{n}}\Bigg).
\end{eqnarray*}
 \cite{moon2016nonparametric}, Chapter 4, investigated  the  statistic 
$$
D_{n,\phi}(\mathbb{P}, \mathbb{Q}):=\int_{\mathbb R^d} \phi\left(\frac{\prod_{l=1}^{p}f_{n,l,i}}{f_{n,i}}\right) d \mathbb{Q}_n. 
$$
If the densities are bounded and continuously differentiable up to order $s$, and $\phi(\cdot)$ is continuously differentiable with bounded derivatives. Under some additional assumptions on the kernels, \cite[Chapter 4]{moon2016nonparametric}  provided the bias of $D_{n,\phi}(\mathbb{P}, \mathbb{Q})$ is of the order $O\left(h_n^s+1/(n h_n^{d})\right)$ and the variance of order $O(1/n)$.  Making use of the Efron-Stein inequality, he also proved the asymptotic normality for $D_{n,\phi}(\mathbb{P}, \mathbb{Q})$. The investigation of the case $\phi(t)$=$ |t-1|$, can not treated by the same techniques and needs more intricate developments.}

\end{remark} 
\begin{remark} \label{comparison_Mason}
	Regarding the difference between this work and of \cite{mason2001}, the first to be resolved is that our statistic is not an integral of a sum of centered independent variables. The initial task was to establish the asymptotic equivalence, in probability, of our truncated statistic $\sqrt{n}(V_n({\bar{C}}) - \mathbb E V_n(\bar{C}))$ and $\sqrt{n}(U_{n,0}(\bar{C})- \mathbb E U_{n,0}(\bar{C}))$, where $\bar{C}= C_1\times \cdots\times C_p$ and, for each $l=1,\ldots,p$, $C_l$ is a bounded Borel set in $\R^{d_{l}}$ satisfying (\ref{lem1eqn3}) and (\ref{lem1eqn4}) of Lemma \ref{lem1} with $g(\cdot)=f_l(\cdot)$ and $\HH=\HH_l:=\{ K_l, K_l^{2}, K_l^{3}\}$. 
This equivalence has been established through Lemmas \ref{lem3} to \ref{lem51}, in \cite{L1testBBD},  so that certain techniques from \cite{mason2001} are extended, in a nontrivial way, to the multivariate case to obtain the asymptotic normality of $\sqrt{n}(U_{n,0}(\bar{C})- \mathbb E U_{n,0}(\bar{C}))$. Finally, we mention that the delicate approximation of $\mathbb EV_n(\tilde{C})$ by $a_n (\tilde{C})$, defined in \eqref{andefinition}, is not discussed in \cite{mason2001}.
\end{remark} 
\begin{remark}
	Recall, in the classical kernel estimation, that the standardizing factor is $\left(n h_n^d\right)^{1 / 2}$ with $h_n \rightarrow 0$, indicating a lower rate of convergence. This is the cost incurred when drawing conclusions regarding the local quantities. The use of the $L_1$-norm        avoids such a lower rate of convergence.	
\end{remark}
\begin{remark}\label{remark4.1}
In this remark, we will point out some differences between the 
$L_1$ and $L_2$ norms. We first recall the following example from Section 6.5 in \cite{devroye2001}. {Let $f(\cdot)$ be the uniform density on $[0,1]$. Let $\epsilon$ be a small positive number, and assume that both density estimates ignore the data: $f_\epsilon(\cdot)$ denotes our catastrophic candidate, that is $1-\epsilon$ on $[0,1]$ and $\epsilon^3$ on $\left[1,1+1 / \epsilon^2\right]$. Now, $g_\epsilon=1$ on $\left[-\epsilon^2, 1-\epsilon^2\right]$. Verify that $\int\left|g_\epsilon-f\right|d x=2 \epsilon^2<2 \epsilon=\int\left|f_\epsilon-f\right|d x$, so that the choices are not even close. However, $\int\left(f_\epsilon-f\right)^2d x=\epsilon^2+\epsilon^4<2 \epsilon^2=\int\left(g_\epsilon-f\right)^2d x$,} so that minimizing the $L_2$-norm, even with $f(\cdot)$ given, picks the wrong density from the set $\left\{f_\epsilon, g_\epsilon\right\}$.  Hence, it can be concluded that criteria relying on the $L_2$-norm are unsuitable for universal properties in density estimation.  In Section 2.3.2.1 of \cite{Scot1992}, for pointwise estimation of $f(\cdot)$ by $\hat{f}(\cdot)$, it has been noted that one appeal of the $L_1$-norm $\int|\hat{f}-f| d x$ is that it pays more attention to the tails of a density than the $L_2$-norm, which de-emphasizes the relatively small density values thereby squaring.  By the fact that the density function has inverse length as its unit, $L_1$ is a dimensionless quantity after integration. The $L_2$-norm, on the other hand, retains the units of inverse length after integration of the squared density error.  We highlight the fact that $0 \leq \int|f-g|d x \leq 2$ while $0 \leq \int(f-g)^2d x \leq \infty$, where $f(\cdot)$ and $g(\cdot)$ are density functions; see Problem 7 in \cite{Scot1992}.  Some asymptotic results for $L_1$ estimates by \cite{Hall1992}  and \cite{MR1118245} support the notion that the practical differences between $L_1$ and $L_2$ criteria are reasonably small except in extreme situations. 
A classical approach for assessing independence involves the $L_2$-norm between the joint density function with the product of the marginal density functions; for instance, see
\cite{rosenblatt1975}, \cite{MR1190260} and \cite{Ahmad1997}, that is defined in \eqref{statVnwei2} for $p=2$ and $q=2$. It is important to mention that Rosenblatt's test is not distribution-free.  
According to Rosenblatt, it has been observed that $L_2$-tests relying on density estimates generally exhibit lower statistical power compared to tests based on sample distribution functions, \cite{Hoeffding} and \cite{Blum}. This finding raises important concerns for many applications, refer to \cite{MR4205649}.
The \cite{Blum} independence test is consistent when the appropriate assumptions are satisfied. However, two challenges arise when utilizing this statistic in a test. Estimating quantiles of the null distribution poses challenges. Furthermore, it is crucial to note that the accuracy of the empirical distribution function estimations deteriorates significantly as the dimensionality of the spaces $\mathbb{R}^{d_1}$ and $\mathbb{R}^{d_2}$ grows. This limitation restricts the effectiveness of the statistic in a multivariate context; for instance, see \cite{gretton2010} for more details and discussions. In Remark \ref{remark41}, we provide some additional comments on important advantages of the $L_1$-norm.  It is seen that the $L_1$ test demonstrates superior performance compared to the $L_2$ test in Section \ref{simulation}. 
\end{remark}
\begin{remark}
To the best of our knowledge, there exists a limited number of nonparametric methods for assessing independence that exhibit a nearly linear computational complexity. For further details, we refer to \cite{deb2020measuring}, \cite{auddy2023exact}, \cite{MR4472834} and \cite{MR4472835}. 
The \citeauthor{MR4472834}'s approach possesses several benefits, with the most prominent being its multivariate nature. Additionally, the computational cost of this approach is generally characterized by a time complexity of $O(n \log n)$ in relation to the sample size $n$. Moreover, it is worth noting that the test threshold remains exact at any sample size, rather than being limited to an asymptotic approximation. In most situations, the existing nonparametric independence tests typically require computation
that scales at least quadratically with the sample size; for instance, we may refer to \cite{pfister2018}.  
In the last reference, it has been observed that the computational complexity of their approach can be reduced by employing linear time approximation techniques (large-scale approximations), as elucidated by the work of \cite{MR3741641} in the context of pairwise HSIC.  
\cite{MR4338371} proposed statistical test based on a $U$-statistic in general setting where the computational complexity is $O(n^2)$. 
In \cite{MR4571116},  the time taken to compute optimal transport-based multivariate ranks is $O(n^3)$ in the worst case. The computational cost with that of permutation $k$-nearest-neighbour methods \cite{berrett}, which run in $O(Bkn \log n)$ time, where in permutation tests any test statistic must be calculated $B$ times. In \cite{MR3858367}, the proposed statistics have the time complexity $O(n^2)$, by considering the dimension is a constant. 
In our setting,  the computational complexity is $O(m_n d n)$, where $m_n$ is the number of points used in the integral approximation, that may be reduced using the idea in \cite{MR1930019}. Identifying certain transformations that can effectively mitigate the computational complexity associated with our test would be interesting. However, obtaining such transformations necessitates a distinct methodology from that employed in the existing literature, and we defer this issue for further investigation. 
\end{remark}

\section{Local alternatives}\label{localalterantive}
	\noindent We determine the power of the test in (\ref{statVn}) against some sequences of local alternatives. Consider the following sequences of local alternatives converging to the null hypothesis at the rate $n^{-1/2}h_n^{-d/4}$, 
	$$               
	H_{\delta,n}: f(\tilde{\mathbf x}) = \prod_{l=1}^{p}f_{l}(\mathbf x_{l})+n^{-1/2}h_n^{-d/4} \delta(\tilde{\mathbf x}),
	$$
	where $\delta(\cdot)$ is an integrable real function and its integral is zero. 
The next result
gives the asymptotic distribution under the local alternatives $H_{\delta,n}$.

\begin{theorem}\label{th4}
		Suppose that $(\mathbf A_1)$ holds and that $h_n \rightarrow 0$ and $n h^{3 d}_n \rightarrow \infty$, as $n \rightarrow \infty$. Then, under $H_{\delta,n}$, we have
		$$
		\lim _{n \rightarrow \infty} \mathbb P\left\{T_{n}(\tilde{\C})>z_{1-\alpha}\right\}=1-\Phi(z_{1-\alpha}-\eta(\delta)),
		$$
where $\Phi(\cdot)$ denotes the  cdf of the standard normal distribution
and  
		$$\eta(\delta)=\frac{1}{\sqrt{2\pi}\sigma}\int_{\tilde{\C}}\delta_1(\tilde{\mathbf x})
		d\tilde{\mathbf x}\quad \textrm{and}\quad\delta_1(\tilde{\mathbf x})=(\delta(\tilde{\mathbf x}))^2\left\{\displaystyle\tilde{\bf{K}} \prod_{l=1}^{p}f_{l}(\mathbf x_l)\right\}^{-1/2}.$$
	\end{theorem}
	
\noindent 
{
\begin{remark}
   Following \cite{auddy2023exact}, let us consider the joint density of $\Tilde{\mathbf X}$ defined by
\begin{eqnarray}\label{eq4.1}
f_{r}(\tilde{\mathbf x}) = (1-r)\prod_{l=1}^{p}f_{l}(\mathbf x_{l})+ r  g(\tilde{\mathbf x}),   
\end{eqnarray}
where $g(\cdot)$ is a density function with marginals $f_l(\cdot)$, $l=1,\ldots,p$, and $r\in [0,1]$. For a sequence $\{r_n\}$ with $r_n\in[0,1]$ for all $n\geq 1$, consider the family of joint density $f_{r_n}(\cdot)$ and the following testing problem :
$
H_0 : r_n=0\quad \textrm{versus} \quad H_{1,n} : r_n>0.  $
Following the same steps in the proof of Theorem \ref{th4}, as $r_n\rightarrow 0$,  we obtain the limit of the power function 
\begin{eqnarray*}
\lim _{n \rightarrow \infty} \mathbb P\left\{T_{n}(\tilde{\C})>z_{1-\alpha}\right\}  =\left\{\begin{array}{lcr}
 \alpha, & \,\,\mbox{ if }\,\,&  n h_n^{d/2} r_n^2 \rightarrow 0,\\
1-\Phi(z_{1-\alpha}- c\, \eta(g)), &\,\, \mbox{ if }\,\,&  n h_n^{d/2} r_n^2 \rightarrow c, \\
1, &\,\,\mbox{ if }\,\,&  n h_n^{d/2} r_n^2 \rightarrow \infty,
 \end{array}\right.
\end{eqnarray*}
where 
$$\eta(g)=\frac{1}{\sqrt{2\pi}\sigma}\int_{\tilde{\C}} g_1(\tilde{\mathbf x})
		d\tilde{\mathbf x}\quad \textrm{and}\quad g_1(\tilde{\mathbf x})=\left(g(\tilde{\mathbf x})- \prod_{l=1}^{p}f_{l}(\mathbf x_l)\right)^2\left\{\displaystyle\tilde{\bf{K}} \prod_{l=1}^{p}f_{l}(\mathbf x_l)\right\}^{-1/2}.$$
\cite{auddy2023exact} obtained a similar result in the bivariate case; for instance, see Corollary 3.1. In their result, the limiting behavior of the power depends on whether $n^{1/4} r_n^2$ converges to $0$,  $\infty$ or some number in $(0,\infty)$ while, for $L_1$-test and $L_2$-test, studied by \cite{rosenblatt1975}, the power depends on the limit of $n h_n^{d/2} r_n^2$ ($d=2$ for $L_2$-test). 
Since $h_n$ is chosen to satisfy $n h_n^d \rightarrow \infty$, then $n^{-1/2} h_n^{-d/4}=o(n^{-1/4})$, and we conclude that the Chatterjee's test, studied in  \cite{auddy2023exact}, is powerless along the local alternative (\ref{eq4.1}) compared with the $L_1$-test. \cite{shi2022power} showed that Hoeffding's test, Blum-Kiefer-Rosenblatt's test and Bergsma-Dassios-Yanagimoto's test are more powerful, with the optimal rate of order $O(n^{-1/2})$, than $L_1$-test, $L_2$-test and Chatterjee's test for the local alternative (\ref{eq4.1}), but presenting some difficulties in estimating quantiles of the null distribution in high dimensional setting. The main question is now to find other local alternatives for which the $L_1$-test is more powerful than the others.
\end{remark}
}

\begin{remark}\label{remark41}
	This paper develops a general framework for testing the independence among continuous random vectors based on $L_1$-norm. Our test achieves the following desirable properties:
\begin{enumerate}
\item[$-$] {\it Full distribution-freeness}. Numerous statistical tests utilize asymptotic distribution-freeness for computationally efficient distributional approximations that result in pointwise asymptotic control over their size. We note that some existing tests of multivariate independence are not distribution-free, which has both computational and theoretical repercussions.
\item[$-$] {\it Transformation invariance.} It is well known that $L_1$-distance is invariant under any smooth monotone transformation; for instance, see \cite{devroye2001}. In multivariate statistics, this kind of invariance is highly interesting.

\item[$-$] {\it Computational efficiency.}
Modern applications require evaluating a dependence measure and the accompanying test to be as computationally efficient as possible, irrespective of their statistical characteristics. Therefore, the effectiveness of the numerical experiments must come first, we may refer to \cite{heller2012}. Furthermore, distribution-free procedures
can alleviate the computational cost of statistical
problems.
\item[$-$] {Consistency under absolute continuity.} The unique condition for the consistency of our tests is that the underlying distributions are absolutely continuous, without the need for any moment requirements.  This makes possible the nonparametric inference under heavy-tailed data-generating distributions such as stable laws \cite{MR3093953} and Pareto distributions \cite{MR2751846}, and it distinguishes our tests from commonly used techniques like traditional distance covariance and energy statistic, for more discussion, refer to \cite{MR4571116}.
\end{enumerate}

\end{remark}

\section{Bandwidth selection}\label{Bandwidth_Selection}
\noindent There are basically no restrictions on
the choice of the kernels $K_l(\cdot)$, $l=1,\ldots,p$, in our setup, apart from satisfying conditions ({\bf A$_1$}). 
In \cite{devroye2001}, it is noticed that for large sample sizes, the shape of the optimal kernel is unique. For example, for $\mathbb{R}$, classical $L_2$ theory \cite{MR148149} shows that for $L_2$ errors, among all positive kernels, the \cite{MR0250422} kernel, $K(x)=\max (3/ 4(1-x^2), 0)$ is best possible. For $\mathbb{R}^d$, \cite{MR539513} showed the $L_2$ optimality of $\max ((1-\|\mathbf x\|^2)^d, 1)$, $\mathbf x\in \mathbb{R}^d$. For the $L_1$-norm, there is evidence that the Epsne{\v{c}}nikov kernel is also best among all positive kernels; refer to \cite{devroye2001}, Chapters 16 and 17. The
selection of the bandwidth, however, is more problematic.
Although any choice of bandwidth $h_n$ satisfying $h_{n}\rightarrow0$ and $nh_{n}^{3d}\rightarrow\infty$ will deliver the
asymptotic distribution in Theorem \ref{thm1}, in practice we need some guidance on how
to select $h_n$. Notice that the quadratic criteria such as the ASE, the ISE, or the MISE are applicable when estimating the density as an element of the $L_2$-space, and $h_n$ can be selected using the cross-validation or plug-in approach; for instance, see \cite{sain1994},  \cite{duong2005}. \cite{Hall1992} demonstrated that the under-smoothing approach yields confidence intervals with greater coverage accuracy than bias-reduced density estimators. \cite{Hall1992} suggested to use $h_n=c 1.05 \hat \sigma n^{-1/5}$, where $0<c<1$ and $\hat \sigma$ is the sample standard deviation. The difficulty of selecting a bandwidth suited to point and interval estimates was examined in \cite{chan2010}. They recommended selecting the local bandwidth with the highest optimal rate and then developing bias reduction estimators for both pointwise and interval estimation. \cite{matt} evaluated spacing bandwidth with low sensitivity to estimate independent components using a nonparametric probability integral transformation and distance covariance. In nonparametric testing, \cite{lee2013} proposed a one-sided $L_q$ approach in testing nonparametric functional inequalities. To evaluate the performance of their test, they used $h_n=c \hat \sigma n^{-1/5}$, where $0.75\leq c\leq 2.5$. 
For testing a parametric model for conditional mean function against a nonparametric alternative, \cite{MR1828537} proposed an adaptive-rate-optimal rule. \cite{MR2504206} proposed, utilizing the Edgeworth expansion of the asymptotic distribution of the test, to select the bandwidth such that the power function of the test is maximized while the size function is controlled. 
\cite{Hall1992} looked at the asymptotic expansion of $\mathbb{E}\left\{\int \mid f_{n, h}-f \mid dx\right\}$ and minimized the main asymptotic terms to obtain a recipe for $h$ as a function of $n, f(\cdot)$ and the kernel function $K(\cdot)$. They then estimate the unknown quantity involving $f(\cdot)$ from the data and propose this as a plug-in bandwidth estimate. However, all other plug-in methods we are aware of are not universal bandwidths. 
While the plug-in procedure displays an analytical solution, which depends on unknown quantities that need to be estimated, the double kernel is performed empirically. Notice also that one may use the maximum likelihood cross-validation method to determine the smoothing parameter; however, this procedure performs very poorly, as indicated in \cite{devroye1997}. The double kernel method uses a pair of kernels, $K(\cdot)$ and $L(\cdot)$, and picks $H=\arg \min _h \int\left|f_{n, h}-g_{n, h}\right|dx,$ where $f_{n, h}(\cdot)$ and $g_{n, h}(\cdot)$ are the kernel estimates with kernels $K(\cdot)$ and $L(\cdot)$, respectively. Assume that $d=1$. If the characteristic functions of $K(\cdot)$ and $L(\cdot)$ do not coincide on an open interval about the origin, then the choice $H$ is consistent, refer to \cite{MR1045250}. 
The ideal bandwidth selection for nonparametric testing differs from that for nonparametric estimation because we must balance the test's size and power rather than the estimator's bias and variance. There are no methods for calculating the appropriate bandwidth for our test, and it is difficult to formulate a theory that provides the solution. The choice of bandwidth determines the sensitivity with which specific types of dependence can be identified and, thus, affects the practical performance of the test. Idealistically, we should select a bandwidth $h $ that provides the best power for a given sample size, but deriving this process is intricate enough to need a separate study. 
Thus, for the present purposes, we provide the combinatorial procedure of \cite{devroye2001} that will be used. To this end, for fixed $\epsilon>0$, define the \cite{MR790571} class
$
\mathcal{A}_\epsilon=\left\{\{\tilde{\mathbf x}: f(\tilde{\mathbf x})>g(\tilde{\mathbf x})\} : f, g \in \mathcal{G}_\epsilon\right\},
$ 
where for every $\epsilon>0$, there exists a finite number $N_\epsilon$ of densities in $\mathcal{F}$, a prespecified class of densities, such that the $L_1$ balls of radius $\epsilon$ centered at these densities cover $\mathcal{F}$, that is, if these chosen densities are $\mathcal{G}_\epsilon=\left\{g_1, \ldots, g_{N_\epsilon}\right\}$, then
$
\mathcal{F} \subseteq \cup_{i=1}^{N_\epsilon} {B_{g_i, \epsilon}}.
$ 
Let $m<n$, and define $\mathcal{A}_{\boldsymbol{\Theta}}$ as the Yatracos class of subsets of $\mathbb{R}^d$ (corresponding to the family of density estimates $f_{n, \theta}, \theta \in \boldsymbol{\Theta}\subset \mathbb R$ ) as the class of all sets of the form
$
A_{\theta_1, \theta_2}=\{\tilde{\mathbf x}: f_{n-m, \theta_1}(\tilde{\mathbf x})>f_{n-m, \theta_2}(\tilde{\mathbf x})\}, \theta_1, \theta_2 \in \boldsymbol{\Theta} .
$ 
We select a parameter $\theta_n$ from $\boldsymbol{\Theta}$ by minimizing the distance
$
\Delta_\theta=\sup _{A \in \mathcal{A}_{\boldsymbol{\Theta}}}\left|\int_A f_{n-m, \theta}dx-\mu_m(A)\right|, 
$
over all $\theta \in \boldsymbol{\Theta}$, where $\mu_m(\cdot)$ denotes the empirical measure defined by the subsample $\tilde{\mathbf X}_{n-m+1}, \ldots, \tilde{\mathbf X}_n$. If the minimum does not exist, we select $\theta_n$ in such a way that $\Delta_{\theta_n}<\inf _{\theta^* \in \boldsymbol{\Theta}} \Delta_{\theta^*}+1 / n.$ Define $f_n(\cdot)=f_{n-m, \theta_n}(\cdot)$. If $\int f_{n-m, \theta}=1$ for all $\theta \in \boldsymbol{\Theta}$, then for the minimum distance estimate $f_n(\cdot)$ as defined above, from Theorem 10.1 of \cite{devroye2001}, we have
\begin{eqnarray}\label{equadevroye}
\int\left|f_n-f\right|dx \leq 3 \inf _{\theta \in \boldsymbol{\Theta}} \int\left|f_{n-m, \theta}-f\right|dx+4 \Delta+\frac{3}{n},~\mbox{where}~ \Delta=\sup _{A \in \mathcal{A}_{\boldsymbol{\Theta}}}\left|\int_A fdx-\mu_m(A)\right| .
\end{eqnarray}
Note that if the estimates do not satisfy the condition $\int f_{n-m, \theta}dx=1$, (\ref{equadevroye}) remains valid, but with the factors of ``3'' replaced by ``5''. \cite{MR3352508}
proposed an adaptive method to estimate density based on data with small measurement errors inspired by the combinatorial method described before. A simulation study in \cite{MR3352508} shows that the proposed estimator converges faster, in terms of the $L_1$-norm than a standard kernel density estimator (Akaike–Parzen–Rosenblatt) to the unknown density.

\section{Monte Carlo experiments}\label{simulation}
\noindent In this section, we examine the empirical performance of the $L_1$-based independence test compared to alternative approaches. 
 For this purpose,
we assess $L_1$-based test's size and power for various non-standard distributions. We assess the power and size of the $L_1$ test in comparison to other pre-existing tests enumerated below. 
The dHSIC test described by \cite{pfister2018} and 
implemented in the \texttt{R} package \texttt{dHSIC} [\cite{Pfister2017}]. The test based on ranks of distances introduced in \cite{heller2012}
and implemented in the \texttt{R} package \texttt{HHG} [\cite{heller2012}]. The test based on the distance covariance discussed in \cite{szekely2009} (will be called Dcov) and implemented in the \texttt{R} package \texttt{energy} [\cite{rizzo2016}]. The mutual information-based test developed in \cite{berrett} (will be called Mintav) and 
implemented in the \texttt{R} package \texttt{IndepTest} [\cite{MINT1}]. {Three  measures based on the generalized distance covariance, 
asymmetric measure $R_n$, 
 symmetric measure $S_n$ based on distance covariance and
 simplified complete measure $Q^{*}_{n}$ based on incomplete $V$-statistics
 and implemented in the \texttt{R} package 
\texttt{EDMeasure} [\cite{MR3858367}]. Finally, the $L_2$-test of independence proposed by \cite{rosenblatt1975}. \cite{MR1190260} showed that if the second order partial derivatives of the density function $f(\cdot)$ are bounded and uniformly continuous on $\mathbb R^2$, 
\begin{eqnarray*}
h_n^{-1}\left(n h_n^{2} \int_{\R^{2}} \left[f_{n}(\tilde{\mathbf x})-\prod_{l=1}^{2}f_{n,l}( x_l)\right]^2\,d\tilde{\mathbf x} -A(n)\right)\stackrel{\D}{\rightarrow} \mathcal{N}(0, 2\vartheta^2)\ \ \mbox{as}\ \ n\rightarrow\infty,
\end{eqnarray*}
where, for $K(\cdot)=K_1(\cdot)=K_2(\cdot)$ be a bounded kernel,
\begin{eqnarray*}
A(n)=\left(\int_{\R} K^2(u_1) \, d u_1\right)^2 -h_n \int_{\R} K^2(u_1) \, d u_1\left(\int_{\R} \left[f_1^2(x)+f_2^2(x)\right]^2\, dx\right),
\end{eqnarray*}
and
\begin{eqnarray*}
\vartheta^2=\int_{\R} f_1^2(x)\, dx \int_{\R} f_2^2(x)\, dx \int_{\R^2} \left(\int_{\R^2} K(t_1)K(t_2)K(t_1+s_1)K(t_2+s_2)\, dt_1\, dt_2\right)^2 ds_1\,ds_2.
\end{eqnarray*}
To compare our test with the $L_2$-test, we use the result above and approximate $A(n)$ and $\vartheta^2$ by using Lemma \ref{Lemma12}.}
We present the results of a Monte Carlo experiment for evaluating the size and power of each test in two-dimensional testing scenarios. We examine a collection of data-generating mechanisms designed to illustrate various types of possible dependence:
\begin{description}
	\item[(i)] We consider first the six simulated examples of unusual bivariate distributions presented in the Supplementary Material of \cite{heller2012} and labeled '4 Clouds', 'W', 'Diamond', 'Parabola', '2 Parabolas', and 'Circle'. Following \cite{heller2012} and \cite{fan2017}, 
	for these examples, we simulated data under both dependence and independence structures.
	The models coded as GEVmodel1 and GEVmodel2 were considered in \cite{fan2017}. 
 \begin{enumerate}
   \item[$-$] The dependence model for GEVmodel1 is as follows: $X_1 = W + N_1$, $X_2 = X_1 + T + N_2$, where $W$ is Weibull distributed, $T$ is Fr\'echet distributed and $N_1$, $N_2$ are independent $\mathcal{N}(0, 0.2^2)$. 
 \item[$-$] The dependence model for GEVmodel2 is as follows: $X_1 = W + N_1$ and $X_2 = 3W + N_2$.
 \end{enumerate}
 These models have non-symmetric distributions with heavy tails.
 \begin{enumerate}
   \item[$-$]We include an example taken from \cite{szekely2007} given by $X_2=X_1 N_3$, with $X_1$ and $N_3$ being independent standard normal variables. 
 \end{enumerate}
	\item[(ii)] We consider the class of sinusoidally dependent data, proposed by \cite{sejdinovic2013}, for which the density function is expressed, for $\ell\in \mathbb{N}$ and $(x_1,x_2)\in [ -\pi, \pi]^2$, by 
	$f_\ell(x_1,x_2)=\frac{1}{4\pi^2}\big( 1+ \sin(\ell x_1)\sin(\ell x_2)\big). $ 
This class of densities is particularly interesting and is considered by \cite{sejdinovic2013} as being challenging to detect dependence; intuitively, this is because as $\ell$ increases, the dependence becomes increasingly localized, while the marginal densities are uniform on $[-\pi, \pi ]$ for each $\ell$. This example was further discussed in \cite{berrett} and \cite{bottcher2019}.
	\item[(iii)] For $\theta\in\lbrack0,1\rbrack$ and $(x_1, x_2) \in \lbrack-1,1\rbrack^{2}$, define the density function
	$
		f_{\theta}(x_1,x_2)= 4^{-1}(1+\theta(1-2 \mathds{1}_{\left\lbrace x^{2}_1+x^{2}_2\leq 2/\pi\right\rbrace }).$
	\item[(iv)]
	Let $X_1$ and $Z$ be independent with $X_1 \sim U\lbrack-1,1\rbrack$, the uniform law on $\lbrack-1,1\rbrack$, and $Z \sim \mathcal{N}(0,1)$, and for a parameter $\rho \in\lbrack 0, \infty)$, 
	let $X_2=|X_1|^{\rho}Z$.
{\item[(v)]
	Let $X_1$ and $Z$ be independent with $X_1 \sim \mathcal{N}(0,1)$ and $Z \sim \mathcal{N}(0,1)$, and for a parameter $\rho \in\lbrack 0, 1)$.
	 Let $X_2$ be defined by  $X_2=\rho X_1+\sqrt{1-\rho^2}Z$.}
\end{description}
In implementing our test, the computing program codes are in \texttt{R}. The joint density estimate and its marginals are computed, without loss of generality,  by using the univariate uniform kernel function, i.e., $K_l(u) =\mathds{1}_{\{|u|\leq 1/2\}}, l=1,\ldots,p,$  and the combinatorial procedure described above for choosing the bandwidth vector $\mathbf{h}_n$ discussed in Remark \ref{remark201}.
This last method was applied only to find bandwidth vector $\mathbf{h}_n$ which allows us to estimate the joint density and then we use each component of $\mathbf{h}_n$ to give an estimate of their marginals and the function $\hat{\mathcal{L}}_{n}(\cdot)$. 
The support set $\tilde{\C}$ is estimated from the sample, specifically the common support set of all functions used in our test statistic, in this
case, is estimated by the set $\prod_{i=1}^{d}\big\lbrack \min_{1\leq j\leq n}({X}_{i,j})-h_n^{(i)}/2, \max_{1\leq j\leq n}({X}_{i,j})+h_n^{(i)}/2\big\rbrack$. The multidimensional integration is based on the SUbregion-Adaptive Vegas Algorithm developed in \cite{HAHN200578} (called suave) and implemented in the \texttt{R} package \texttt{cubature}. We chose the bandwidth $h$ for the $L_2$ test through the cross-validation process, employing the Gaussian kernel. 
We test the null hypothesis $\Hc_0 $ with a significance level of $\alpha = 0.05$. We run $1,000$ replications with the adaptive permutation size $B=200$ for all empirical measures that require a permutation procedure, i.e., dHSIC, HHG, Mintav, 
$R_n$, $S_n$, $Q^{*}_{n}$ and Dcov.
\\
{
\begin{table}[!ht]
	\begin{center}
		\caption{Size and power of $L_1 $, $L_2$,  dHSIC, HHG, Mintav, $R_n$, $S_n$, $Q^{*}_{n}$ and Dcov  for some unusual relations with  $n=50 $. Results based on $1,000$ replications.}
		\label{tab:table1}
		\vspace{0.2cm} 
  	\begin{tabular}{llcccccccccc}
			\hline\hline
			Distribution &Indep&$L_1$&$L_2$&dHSIC &HHG &Mintav &$R_n$&$S_n$&$Q^{*}_{n}$ &Dcov\\ 
			\hline \hline
			IndNorm &Yes 
         &$0.023$&$0.012$ &
         $0.033$& $0.043 $& $ 0.039$& $   0.045 $& $ 0.049 $& $    0.045$& $  0.043$\\
			\hline
			$4$ Clouds &Yes &$0.055$&$0.059$ &$  0.056$ &$ 0.057$ &$  0.042$ &$   0.052$ &$  0.055  $ &$   0.041 $ &$ 0.047$\\ 
			$4$ Clouds &No  &$1$ &$1$ &$
   0.092   $ &$  1 $ &$     1 $ &$  0.095 $ &$ 0.096  $ &$   0.138$ &$  0.086$\\ 
			\hline
			W  &Yes &$0.058$&$0.044$ &$
    0.057 $&$0.051 $&$ 0.048 $&$  0.057  $&$0.054    $&$ 0.057   $&$0.05
   $ \\ 
			W & No &$1$ &$1$ &$ 0.989$&$ 0.999     $&$ 1 $&$  0.933   $&$0.93   $&$  0.863 $&$ 0.915$ \\ 
			\hline
			Parabola &Yes&$0.041$ &$0.028$&$0.053 $&$ 0.05   $&$0.04  $&$ 0.052 $&$ 0.065  $&$   0.058 $&$ 0.051$\\ 
			Parabola&No&$0.972$&$ 0.895 $&$  0.988 $&$0.998     $&$ 1  $&$ 0.958 $&$ 0.966    $&$ 0.381 $&$ 0.938$\\ 
			\hline
			$2$ Parabolas &Yes&$0.032$ &$0.025$ &$ 0.051 $&$0.046 $&$ 0.032 $&$  0.055 $&$ 0.053   $&$  0.055 $&$  0.05$\\ 
			$2$ Parabolas &No&$0.996$ &$1$&$   1  $&$   1      $&$1 $&$  0.312 $&$ 0.301   $&$  0.858 $&$ 0.274$ \\
			\hline
			Circle &Yes&$0.055$&$ 0.058  $ &$ 0.053 $&$0.056  $&$0.046  $&$ 0.047 $&$ 0.051   $&$  0.051 $&$ 0.041$ \\
   Circle &No&$1$ &$ 0.994 $&$ 0.897$&$ 0.985 $&$ 0.993  $&$ 0.094  $&$0.088 $&$    0.265 $&$ 0.076$ \\ 
			\hline
			$X_2 =X_1Z$ &Yes&$0.025$&$0.022$&$0.051 $&$0.056  $&$0.033 $&$  0.046 $&$ 0.048 $&$     0.06 $&$ 0.037$ \\ 
			$X_2=X_1 Z$&No&$0.892$&$0.843$ &$ 0.95 $&$0.992$&$  0.987 $&$  0.665  $&$0.664    $&$ 0.669 $&$ 0.633
  $ \\ 
		\hline
			Diamond &Yes&$0.033$&$0.017$&$ 0.044$&$ 0.054   $&$0.03  $&$ 0.035 $&$ 0.035   $&$  0.046 $&$ 0.034$\\ 
			Diamond &No&$0.473$&$0.248$&$ 0.51 $&$0.643 $&$ 0.095  $&$ 0.038 $&$ 0.037   $&$  0.141 $&$ 0.031$\\ 
			\hline
			GEVmodel$1$ &Yes&$0.038 $&$0.028$ &$0.053$&$ 0.059  $&$0.036 $&$  0.052 $&$ 0.053   $&$  0.053 $&$ 0.047$\\
			GEVmodel$1$&No&$0.851$&$0.786$ &$ 0.863$&$ 0.937  $&$0.916 $&$  0.355 $&$ 0.356   $&$  0.731 $&$ 0.337$\\
			\hline
			GEVmodel$2$ &Yes&$0.023$&$0.014$&$0.053$&$ 0.052  $&$0.035  $&$ 0.052  $&$ 0.04   $&$  0.051 $&$ 0.042$\\
GEVmodel$2$ &No&$ 0.216$&$0.099$&$ 0.268 $&$0.401 $&$ 0.356   $&$0.645 $&$ 0.648  $&$   0.293 $&$ 0.624$\\ 
			\hline\hline
		\end{tabular} 
 \end{center}
\end{table}}

Table \ref{tab:table1} represents size and power comparisons between $L_1$ and all testing procedures described above with sample size $n=50$. In Setting {\bf(i)}, all the cited tests have levels close to the nominal level of $0.05$ for independent simulated data. For examples with dependent data, the $L_1$ test shows the same or better power compared to $L_2$,
dHSIC, HHG and Mintav for `W', `Parabola', `2 Parabolas' and `Circle'. We see also that for `4 Clouds' where dHSIC, $R_n$, $S_n$, $Q^*_n$, Dcov perform poorly, $L_1$, $L_2$,  HHG and Mintav have good power. $R_n$, $S_n$ and Dcov present good performance  only for `GEVmodel2', `Parabola', `W' and  
\cite{szekely2007}'s setting. While $Q_n^*$  presents good performance only for `W', `2 Parabola', `GEVmodel2' and \cite{szekely2007}'s setting.

\begin{figure}[!ht]

\begin{subfigure}{0.497\textwidth}
   \includegraphics[width=\textwidth]{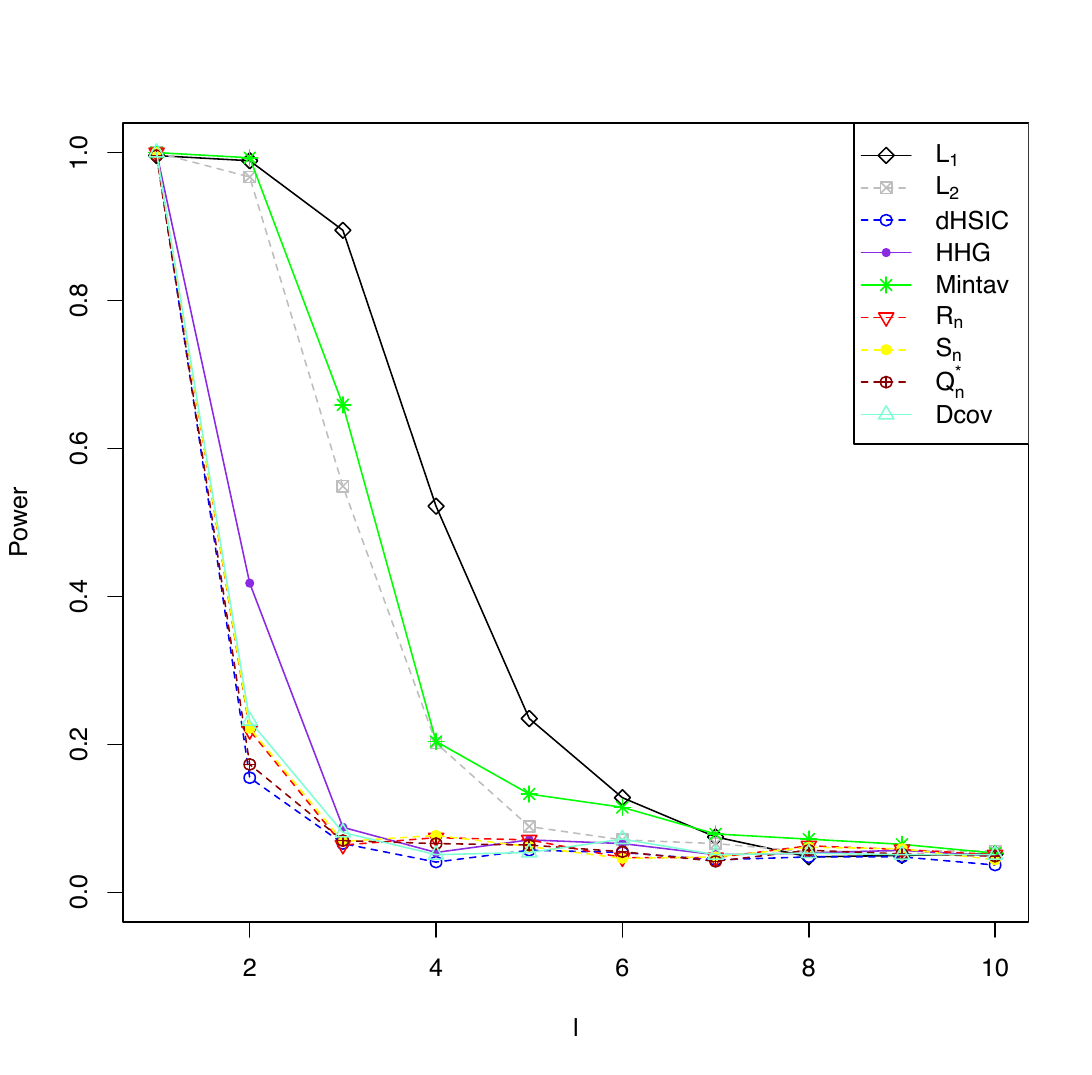}
  \caption{(ii), $n=200$}
  \label{fig:second}
\end{subfigure}
\hfill
\begin{subfigure}{0.497\textwidth}
  \includegraphics[width=\textwidth]{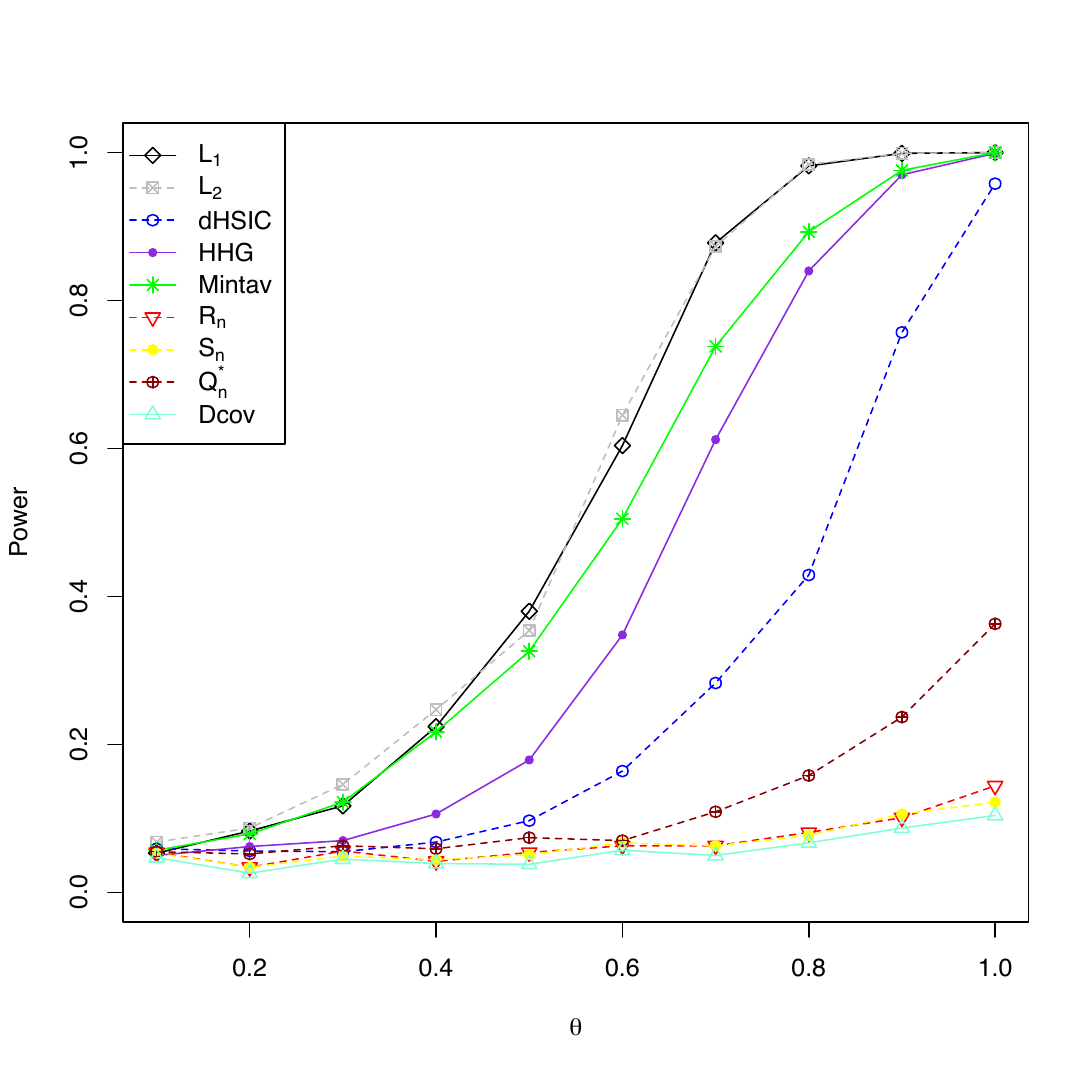}
  \caption{(iii), $n=200$}
  \label{fig:third}
\end{subfigure}
\hfill
\begin{subfigure}{0.497\textwidth}
   \includegraphics[width=\textwidth]{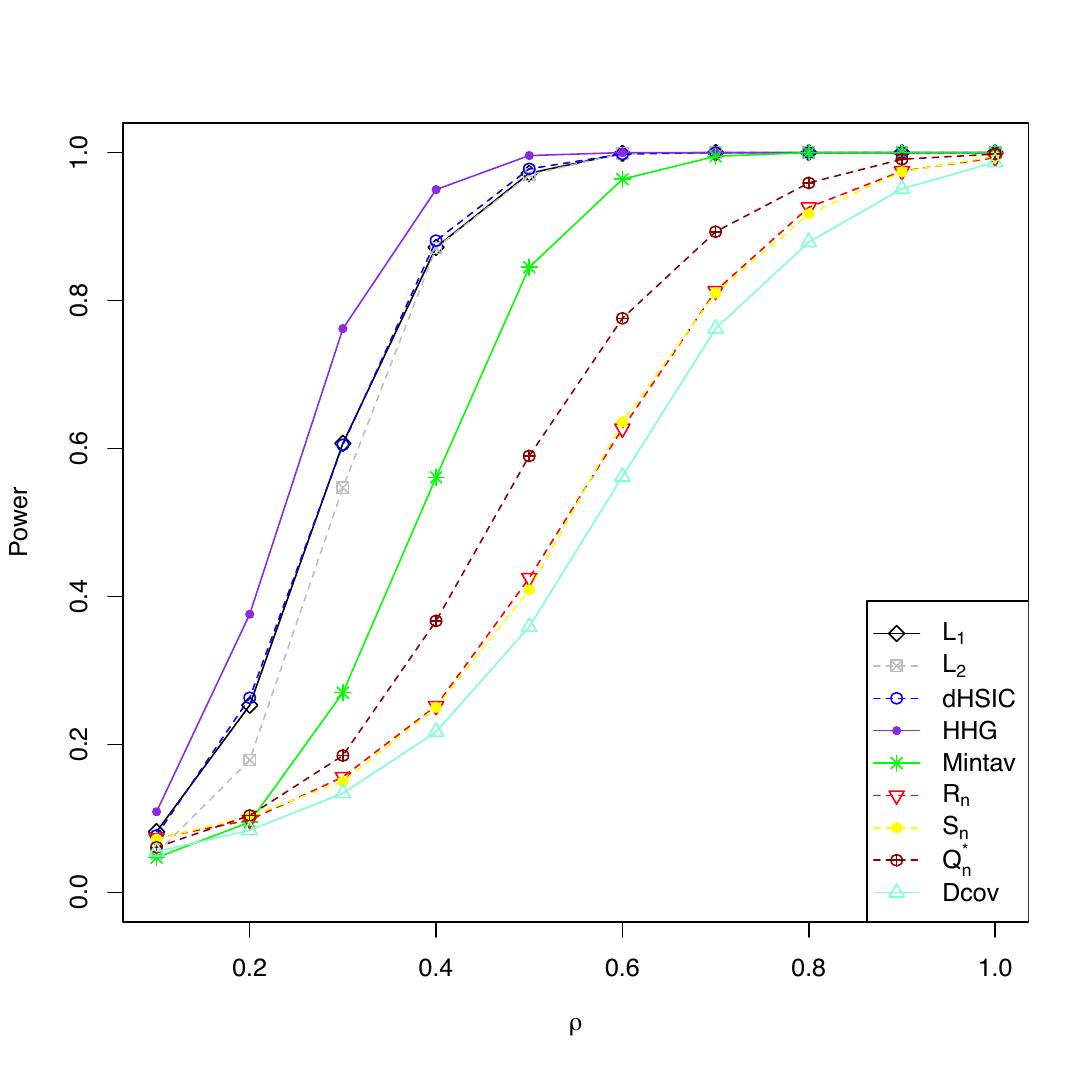}
  
  \caption{(iv),  $n=200$}
  \label{fig:four}
\end{subfigure} 
    \hfill
\begin{subfigure}{0.497\textwidth}

  \includegraphics[width=\textwidth]{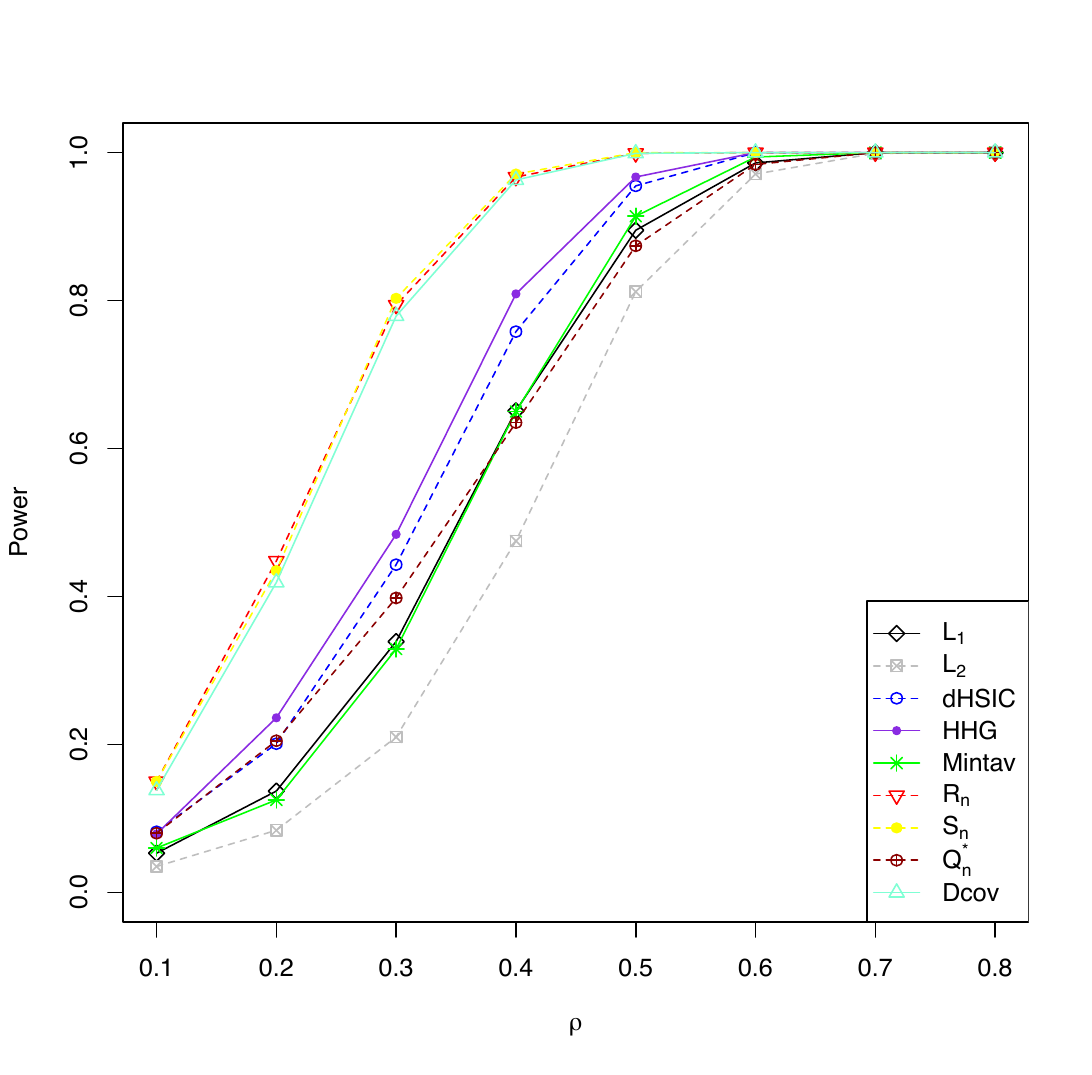}
  
  \caption{(v), $n=100$}
  \label{fig:five}
\end{subfigure} 
\caption{
Power curves for the different tests
}
\label{fig2}
\end{figure}
For the last three examples, we also compare our proposed test with all previously cited dependence measures, the Mintav test was proposed in \cite{berrett}, where it was shown that it outperforms many other tests for sinusoidal dependence (Setting (ii)) and underperforms in Setting (iv) when $\rho$ increases. Figure \ref{fig2} shows that, for $n=200$, the $L_1$ test has solid performances for Settings {\bf(ii)}, {\bf(iii)} and {(iv)}. In Setting {\bf(ii)}, we note that the $L_1$ test keeps a better power than Mintav, even if $\ell$ increases, while the other tests have almost no power when $\ell\geq 2$. In Setting {\bf(iii)}, $L_1$ and Mintav provide good power when $\theta$ increases. The $L_1$ test outperforms the other tests. In Setting {\bf(iv)}, $L_1$, HHG and dHSIC have similar power when $\rho$ increases and performs better than Mintav. 
In the last Setting {\bf(v)}, the $R_n$ and $S_n$ tests outperform the other tests. These measures share similar properties and asymptotic distributions with distance covariance
 quickly detecting the dependence in the Gaussian setting.  We highlight that our test necessitates the selection of bandwidth based on a priori information derived from the sample. Numerous tests of independence use a permutation or bootstrap strategy to determine the critical values, while other tests choose a posteriori the parameter with the greatest power. The aforementioned procedures render the decision uncertain, which is inconvenient for real-world data analysis applications. In conclusion,
due to the wide variety of potential forms of dependence, there is no universally most powerful test. However, if the specific nature of the dependence is known beforehand, it may be feasible to create a tailored test that yields high power.

\section{Concluding remarks}\label{conclusion}
\noindent The notion of independence is of the highest significance in the domain of probability and statistics. It sets probability apart from being only a part of measure theory and forms the foundation for statistical theory and the approach practitioners take to modelling. Statisticians often need to determine the realism of assumptions of independence. This is vital for assessing whether certain theoretical aspects of procedures would hold and for evaluating the fit quality of a statistical model, one can refer to \cite{MR4338371}. Although many existing tests for independence focus on two random vectors, they frequently do not have clear methods to evaluate mutual independence among more than two random vectors with different dimensions. This paper presents a new method to assess the independence of several random vectors, regardless of their dimensionality.  Our method uses the $L_1$-distance between the joint density and the product of the marginal densities. Under the null hypothesis, we use the Poissonization techniques to find, for the first time, the asymptotic normal approximation of the corresponding statistic. We do this without making assumptions about the regularity of the underlying Lebesgue density $f(\cdot)$. 
Also, and this was a surprise, the limiting distribution of the statistics based on $ L_1$-distance does not depend on $f(\cdot)$. We show that the tests have nontrivial local power against a subset of local alternatives that converge to the null at the rate of $n^{-1/2}h_n^{-d/4}$. Lastly, simulations are used to study how the tests behave for a moderate sample size.  {One aspect that remains unexplored in this article is the best selection of the smoothing parameters to maximize the power of the proposed tests. The subject at hand holds significant importance and warrants dedicated research effort.  We defer this matter to a forthcoming investigation.} 
Multiple avenues exist for developing our method further. 
When we look at modern machine learning algorithms, we may need to deal with kernel density estimation with complicated kernels that are not chosen by the user and may even be irregular and asymmetric on Riemannian manifolds with Riemann integrable kernels. A challenging task is to consider this setting. It would be interesting to extend the present work to the problem of testing conditional independence in an incomplete data setting that requires nontrivial mathematics that goes well beyond the scope of the present paper.
In the context of the serially dependent data, a future research direction would be to investigate the problem of testing independence as such that was examined in this work. 

\appendix
\section{Proof of  Theorems}

\label{proofs}
\noindent In this section, we give the proof of all theoretical results in this work. The proof of Theorem \ref{thm1} is quite involved.
For the sake of readability, the proof of Theorem \ref{thm1} is presented with a sequence of lemmas (given in \cite{L1testBBD}) providing the necessary techniques at every demonstration step. Lemma \ref{lem1} shows that Borel sets exist, allowing the use of the truncation procedure and providing some technical results on convolutions which are crucial for the proof. 
To prove the statements of this lemma, we suppose that the kernel functions have compact supports.
Lemma \ref{lem2} is useful to investigate the behavior of the difference between the normalized and truncated statistic $\sqrt{n} V_n(A)$ and the normalized statistic $\sqrt{n} V_n$, for any Borel set $A$ of $\R^{d}$, and shows that this difference is asymptotically negligible for large $n$ and large set $A$. 
Lemma \ref{lem3} and Lemma \ref{lem51} show that $\sqrt{n} V_n({\bar{C}})$ is asymptotically equivalent to $\sqrt{n}\int_{\bar{C}} |f_{n}(\tilde{\mathbf x})-\prod_{l=1}^{p}\mathbb E f_{n,l}(\mathbf x_l)|d\tilde{\mathbf x}.$ The problem in dealing with $V_n$ is its Poissonization\footnote{The interest of Poissonization relies on the nice properties of Poisson processes, namely the independence of their increments and the behavior of their moments. These properties considerably simplify calculations, to be more precise, if $\eta$ is a Poisson random variable independent of the i.i.d. sequence $\{X_i: i \in \mathbb{N}\}$, $X_0=0$, and if $A_k, k \in \mathbb{N}$, are disjoint measurable sets, then the processes $\sum_{i=0}^\eta \mathds 1\left(X_i \in A_k\right) \delta_{X_i}, k=1,2, \ldots$, are independent.} does not allow obtaining the independence of increments for disjoint sets, which are required in establishing Lemma \ref{lem5} by using a type of Berry-Esseen approximations for sums of independent random variables, due to \cite{sweeting}, for easy reference, see Theorem \ref{Thesweeting}. Intermediate results about the Poissonization step are given in five lemmas.
In this respect, Lemma \ref{lem4} establishes the asymptotic equivalence between the expectation of Poissonized and the non-Poissonized terms. Lemma \ref{lem5} allows us to derive the result of Lemma \ref{lem51} and provides the limit variance of the Poissonized version of $\sqrt{n} V_n(\bar{C})$. Lemma \ref{lem6} gives the asymptotic distribution of the Poissonized and normalized version of $\int_{\bar{C}} |f_{n}(\tilde{\mathbf x})-\prod_{l=1}^{p}\mathbb E f_{n,l}(\mathbf x_l)|d\tilde{\mathbf x}.$ This result follows from Lemma \ref{lem8} and Lemma \ref{lemC} by applying Theorem $1$ of \cite{shergin} to the sum of $m-$dependent random fields deduced by partitioning the integral into integrals on small disjoint domains. Lemma \ref{lem7} is devoted to the dePoissonization step.

\vspace{.3cm}

\noindent{\bf Proof of Theorem~\ref{thm1}.}
\label{sub:proof_of_theorem_thm1}
Notice that, by Lemma \ref{lem1}, for each $l=1,\ldots,p$, there exists a sequence of Borel sets $\{C_{l,k}\}_{k\geq 1}$ in $\R^{d_l}$ such that, for each $k\geq 1$, $C_{l,k}$ has finite Lebesgue measure and satisfies \textnormal{(\ref{lem1eqn3})} and \textnormal{(\ref{lem1eqn4})} of Lemma \ref{lem1} with $g(\cdot)=f_l(\cdot)$ and $\HH=\HH_l$, and
\begin{eqnarray}\label{eqn81}
\lim_{k\rightarrow \infty} \int_{C_{l,k}^c}f_{l}(\mathbf u_l)\, d\mathbf u_l=0,\quad \text{where}\, \, C_{l,k}^c \, \,\text{is the complement set of}\,\, C_{l,k}.
\end{eqnarray}
Let $\bar{C}_{k}=C_{1,k}\times\cdots\times C_{p,k}$ for $k\geq1$. By Lemma \ref{lem7}, for each $k\geq1$, as $n\rightarrow\infty$, we have 
\begin{eqnarray*}
  \frac{\sqrt{n}\big(V_{n}(\bar{C}_{k})-\mathbb E V_{n}(\bar{C}_{k})\big)}{\sqrt{n\mathrm{Var}\big({U_{\eta,0}(\bar{C}_{k})}\big)}}\stackrel{\D}{\rightarrow} { \mathcal{N}(0,1)},
\end{eqnarray*}
and, by (\ref{Varf}), we infer 
$
  \lim_{n\rightarrow\infty} n\mathrm{Var}\big(U_{\eta,0}(\bar{C}_{k}) \big)=\sigma^2\prod_{l=1}^{p}\int_{C_{l,k}}f_l(\mathbf u_l)\,d\mathbf u_l.
$ 
By Lemma \ref{lem2}, for each $k\geq1$, it readily follows that
\begin{eqnarray*}
  \limsup_{n\rightarrow \infty} n \mathrm{Var}\big(V_{n}(\bar{C}_{k}^{c})\big)
  \leq 64 \int_{\bar{C}_{k}^{c}} f(\tilde{\mathbf x})\,d\tilde{\mathbf x}\leq 64 \sum_{l=1}^{p}\int_{\bar{C}_{l,k}^{c}} f_l(\mathbf u_l)\,d\mathbf u_l.
\end{eqnarray*}
Now, by (\ref{eqn81}) and Theorem $4.2$ of \cite{billy}, we deduce that, as $n\rightarrow \infty$, $
  \sqrt{n}(V_{n} -\mathbb E V_{n})\stackrel{\D}{\rightarrow}  \mathcal{N}(0,\sigma^2).
$ 
We obtain the result of Theorem \ref{thm1}, as sought.
\hfill$\blacksquare$

\vspace{.3cm}
\noindent{\bf Proof of Theorem~\ref{th2}.} 
\label{sub:proof_of_theorem_ref_thm2}
Introduce $L_{n}(\cdot)$ by setting 
\begin{eqnarray*}
  L_{n}(\tilde{\mathbf x})&=& \prod_{l=1}^{p}K_{l}\Bigg(\frac{\mathbf x_l-\mathbb X^{l}}{h_n}\Bigg) 
  - \sum\limits_{l=1}^{p} K_{l}\Bigg(\frac{\mathbf x_l-\mathbb X^{l}}{h_n}\Bigg)\prod_{j\neq l} \mathbb EK_{j}\Bigg(\frac{\mathbf x_j-\mathbb X^{j}}{h_n} \Bigg)
  + (p-1)\prod_{l=1}^{p} \mathbb EK_{l}\Bigg(\frac{\mathbf x_l-\mathbb X^{l}}{h_n}\Bigg)  .
\end{eqnarray*}
By the statement (\ref{equa6.18Fact6.1}), we get, for some constant $\Upsilon_1>0$, 
\begin{eqnarray}\label{eqth22}
\left|\sqrt{n}\mathbb E {U_{n,1}(\tilde{\C})}
-\mathbb E|Z|\int_{\tilde{\C}}
\sqrt{{\mathcal{L}_{n,1}(\tilde{\mathbf x})}}\,d\tilde{\mathbf x} \right| 
\leq \frac{ \Upsilon_1 }{\sqrt{n}h_{n}^{3d}}\int_{\tilde{\C}}\frac{\displaystyle \mathbb E|L_{n}(\tilde{\mathbf x})|^{3} }{ \displaystyle\mathcal{L}_{n,1}(\tilde{\mathbf x}) }\,d\tilde{\mathbf x},
\end{eqnarray}
where $\mathcal{L}_{n,1}(\tilde{\mathbf x})= n\mathrm{Var}\big(\Gamma_{n,1}(\tilde{\mathbf x})\big)$. This last bound is, since, for each $l=1,\ldots,p$, $K_l(\cdot)$ is bounded, $\lambda(\bar{\C})<\infty$ and $h_{n}^{-2d}\mathrm{Var}\big(L_{n}(\tilde{\mathbf x})\big)=\mathcal{L}_{n,1}(\tilde{\mathbf x})$, less than or equals to $\frac{\upsilon_2 }{\sqrt{n}h_{n}^{d}}$ for some constant $\upsilon_2>0$.
Now, since, for each $l=1,\ldots,p$, $f_l(\cdot)$ is bounded on $\C_{l}^{\delta}$ for some $\delta>0$, thus, by arguing as in the proof of Lemma \ref{lem3}, we get 
\begin{eqnarray}\label{eqth23}
\sqrt{n}\mathbb E\big| V_{n}(\tilde{\C})- U_{n,1}(\tilde{\C})\big|=O\left
(\frac{1}{\sqrt{nh_{n}^{d}}}\right). 
\end{eqnarray}
Combining  (\ref{eqth22})-(\ref{eqth23}) and Lemma \ref{Lemma12}, we obtain $
  \sqrt{n}\big(V_{n}(\tilde{\C}) -a_n(\tilde{\C})\big)\stackrel{\D}{\rightarrow}  \mathcal{N}(0,\sigma^2)\ \ \mbox{as}\ \ n\rightarrow\infty.$ 
Therefore, for $n$ large enough, we 
obtain readily
\begin{eqnarray*}
  \mathbb P\Bigg\{\frac{\sqrt{n}}{\sigma}\big(V_{n}(\tilde{\C}) -a_n(\tilde{\C})\big)\leq z_{1-\alpha}\Bigg\}\rightarrow1-\alpha.
\end{eqnarray*}
Thus, the $\alpha$-level test rejects the null hypothesis if
$
  V_{n}(\tilde{\C}) > a_n(\tilde{\C})+ \frac{\sigma}{\sqrt{n}}z_{1-\alpha}.
$
This proves that the test has an asymptotic error probability equal to $\alpha$. \hfill$\blacksquare$

\vspace{.3cm}
\noindent{\bf Proof of Theorem \ref{th3}.} 
Making use of the triangular inequality, we readily obtain
\begin{eqnarray*} 
 |V_n(\tilde{D})-V(\tilde{D})|&\leq &  \int_{ \tilde{D}}\left|f_{n}(\tilde{\mathbf x})-\prod_{l=1}^{p}f_{n,l}(\mathbf x_l)-\Big(\mathbb Ef_{n}(\tilde{\mathbf x})-\prod_{l=1}^{p}\mathbb Ef_{n,l}(\mathbf x_l)\Big)\right|d\tilde{\mathbf x}\\&&+\int_{\tilde{D}} \left|f(\tilde{\mathbf x})-\prod_{l=1}^{p}f_{l}(\mathbf x_l)-\Big(\mathbb Ef_{n}(\tilde{\mathbf x})-\prod_{l=1}^{p}\mathbb Ef_{n,l}(\mathbf x_l)\Big)\right|d\tilde{\mathbf x}. 
\end{eqnarray*}
We now observe that 
\begin{eqnarray*} \lefteqn{\left|f_{n}(\tilde{\mathbf x})-\prod_{l=1}^{p}f_{n,l}(\mathbf x_l)-\Big(\mathbb Ef_{n}(\tilde{\mathbf x})-\prod_{l=1}^{p}\mathbb Ef_{n,l}(\mathbf x_l)\Big)\right|
  }\\&\leq& \left|f_{n}(\tilde{\mathbf x})-\mathbb Ef_{n}(\tilde{\mathbf x})\right|+
  \sum_{j=1}^{p}
  \left|f_{n,j}(\mathbf x_j)-\mathbb Ef_{n,j}(\mathbf x_j)\right|\prod_{l=1}^{j-1}\mathbb Ef_{n,l}(\mathbf x_l)\prod_{l=j+1}^{p}f_{n,l}(\mathbf x_l) ,
\end{eqnarray*}
and
\begin{eqnarray*} &&\left|f(\tilde{\mathbf x})-\prod_{l=1}^{p}f_{l}(\mathbf x_l)-\Big(\mathbb Ef_{n}(\tilde{\mathbf x})-\prod_{l=1}^{p}\mathbb Ef_{n,l}(\mathbf x_l)\Big)\right|
  \\&&\hspace{0.5cm}\leq \big|f(\tilde{\mathbf x})-\mathbb Ef_{n}(\tilde{\mathbf x})\big|+
  \sum_{j=1}^{p}\big|f_{j}(\mathbf x_j)-\mathbb Ef_{n,j}(\mathbf x_j)\big|\prod_{l=1}^{j-1}\mathbb Ef_{n,l}(\mathbf x_l)\prod_{l=j+1}^{p}f_{l}(\mathbf x_l) .
\end{eqnarray*}
Now, for each $l=1,\ldots, p$,  we have 
\begin{eqnarray*} 
  \int_{ {\R}^{d_l}} f_{n,l}(\mathbf x_l)\, d\mathbf x_l=1 \quad \textrm{and}\quad  \int_{ {\R}^{d_l}}\mathbb E f_{n,l}(\mathbf x_l)\, d\mathbf x_l \rightarrow   \int_{ {\R}^{d_l}} f_{l}(\mathbf x_l)\, d\mathbf x_l=1.
\end{eqnarray*}
Then, for $n$ large enough, there exists a constant $C>0$ such that
\begin{eqnarray*} 
  \mathbb E  |V_n(\tilde{D})-V(\tilde{D})|&\leq & \int_{ \tilde{D}}\mathbb E\big|f_{n}(\tilde{\mathbf x})-\mathbb Ef_{n}(\tilde{\mathbf x})\big|\,d\tilde{\mathbf x}+
  C \sum_{j=1}^{p}  \int_{D_j}\mathbb E\big|f_{n,j}(\mathbf x_j)-\mathbb Ef_{n,j}(\mathbf x_j)\big|\,d\mathbf x_j\\&&+
  \int_{ \mathbb{R}^{d}}\big|f(\tilde{\mathbf x})-\mathbb Ef_{n}(\tilde{\mathbf x})\big|\,d\tilde{\mathbf x}+
  C \sum_{j=1}^{p}  \int_{ \mathbb{R}^{d_j}} \big|f_{j}(\mathbf x_j)-\mathbb Ef_{n,j}(\mathbf x_j)\big|\,d\mathbf x_j.
\end{eqnarray*}
Making use of the Cauchy-Schwartz and Young's inequalities, we infer that 
\begin{eqnarray*} 
  \int_{\tilde{D}}\mathbb E\big|f_{n}(\tilde{\mathbf x})-\mathbb Ef_{n}(\tilde{\mathbf x})\big|\,d\tilde{\mathbf x}&\leq &  \int_{ \tilde{D}}\Big(\mathrm{Var}\big(f_{n}(\tilde{\mathbf x})\big)\Big)^{1/2}d\tilde{\mathbf x}\leq \frac{\big(\mu(\tilde{D})\big)^{1/2} }{\sqrt{n h^d_n}}\left(\int_{ \tilde{D}} \left\{\prod_{l=1}^{p}K_{l}
 \right\}^2\ast f(\tilde{\mathbf x})\, d\tilde{\mathbf x}\right)^{1/2}\\&= &O\bigg
  (\frac{1}{\sqrt{nh_{n}^{d}}}\bigg). 
\end{eqnarray*}
Likewise, for each $j=1,\ldots, p$, we have 
$\int_{D_j}\mathbb E\big|f_{n,j}(\mathbf x_j)-\mathbb Ef_{n,j}(\mathbf x_j)\big|\,d\mathbf x_j
  = O\bigg
  (\frac{1}{\sqrt{nh_{n}^{d_j}}}\bigg). 
$  
On the other hand, by Theorem $1$ in Chapter $2$ of \cite{devroye1}, we readily infer 
\begin{eqnarray*} 
  \int_{ \mathbb{R}^{d}}\big|f(\tilde{\mathbf x})-\mathbb Ef_{n}(\tilde{\mathbf x})\big|\,d\tilde{\mathbf x} \rightarrow 0 \quad \textrm{and}\quad  \int_{ \mathbb{R}^{d_j}} \big|f_{j}(\mathbf x_j)-\mathbb Ef_{n,j}(\mathbf x_j)\big|\,d\mathbf x_j \rightarrow 0,\, \, j=1,\ldots,p.
\end{eqnarray*}
Using the fact that  $n^{-1/2} h^{-d/2}\rightarrow 0$, we conclude that $V_n(\tilde{D})\stackrel{\PP}{\rightarrow} V(\tilde{D})$. To complete the proof of the theorem, it suffices to show that $
  \mathbb   P\Big(V_n(\tilde{D})<n^{-1/2}\big(\sigma z_{1-\alpha}+ a_n(\tilde{D})\big)\Big)\rightarrow 0.$ 
Now, we can check easily that $a_n(\tilde{D})=O(h^{-d/2})$, then for $n$ large enough, there exists some constant $\nu>0$ in such a way that $\frac{1}{\sqrt{n}}\big(\sigma z_{1-\alpha}+ a_n(\tilde{D})\big)- V(\tilde{D}) <-\nu.$
Therefore, we infer that
\begin{eqnarray*} 
  \mathbb   P\left(V_n(\tilde{D})<\frac{1}{\sqrt{n}}\big(\sigma z_{1-\alpha}+ a_n(\tilde{D})\big)\right)&\leq&\mathbb P\big(V_n(\tilde{D}) - V(\tilde{D})<-\nu\big)
  \leq  \mathbb P\big(\big|V_n(\tilde{D}) - V(\tilde{D})\big|>\nu\big)\rightarrow 0.
\end{eqnarray*}
Hence, we obtain the desired conclusion
  of Theorem \ref{th3}.
\hfill $\blacksquare$

\vspace{.3cm}
\noindent{\bf Proof of Theorem \ref{th4}.} 
\noindent Notice that 
\begin{eqnarray}\label{equacontigue}
 T_{n}(\tilde{\C})&=&\sqrt{n}{\sigma}^{-1}V_{n}(\tilde{\C})-{\sigma}^{-1}a_{n}(\tilde{\C})
\nonumber \\&=&\sqrt{n}{\sigma}^{-1} \left(V_{n}(\tilde{\C})-\mathbb EV_{n}(\tilde{\C})
)
\right)
+\sqrt{n}{\sigma}^{-1}\left(\mathbb EV_{n}(\tilde{\C})-n^{-1/2}a_{n}(\tilde{\C})\right).
\end{eqnarray}
To establish that the first term on the right side of (\ref{equacontigue}) converges to $\mathcal{N}(0,1)$ under $H_{\delta,n}$, we follow the proof of Theorem \ref{thm1} with some slight differences. Note that in this proof, we have not assumed independence. So for this reason we will use a type of Berry-Esseen bound for sums of independent random variables due to \cite{sweeting}, the theorem in \cite{pinelis} and Theorem 2 on page 63 of \cite{stein}. Consider the second term in (\ref{equacontigue}). Under $H_{\delta,n}$, since $nh_n^{3d}\rightarrow \infty$, by using the arguments of the proof of Lemma \ref{lem4}, we get {
$$
\lim_{n\rightarrow \infty}\sqrt{n}\mathbb EU_{n,1}(\tilde{\C})-\int_{\tilde{\C}} \mathbb E|\gamma_{n,1}(\tilde{\mathbf x})Z_1+\delta_n(\tilde{\mathbf x})|\, d\tilde{\mathbf x}=0,
$$
where $\gamma_{n,1}(\tilde{\mathbf x})=\sqrt{{\mathcal{L}_{n,1}(\tilde{\mathbf x})}}$ and $\delta_n(\tilde{\mathbf x})= h_n^{-d/4}\int_{\mathbb{R}^d} \delta(\tilde{\mathbf x}-\mathbf u h_n) 
\left\{\prod_{l=1}^{p}K_{l}(\mathbf u_l)
 \right\}\,d\mathbf u.$ Now we have 
\begin{eqnarray*}
  \lefteqn{\mathbb E|\gamma_{n,1}(\tilde{\mathbf x})Z_1+\delta_n(\tilde{\mathbf x})|- \mathbb E|\gamma_{n,1}(\tilde{\mathbf x})Z_1|}\\&=&2\delta_n(\tilde{\mathbf x}) \Phi\left(\frac{\delta_n(\tilde{\mathbf x})}{\gamma_{n,1}(\tilde{\mathbf x})}\right)+ 2\gamma_{n,1}(\tilde{\mathbf x})\phi\left(\frac{\delta_n(\tilde{\mathbf x})}{\gamma_{n,1}(\tilde{\mathbf x})}\right)
   -\delta_n(\tilde{\mathbf x})-2\gamma_{n,1}(\tilde{\mathbf x})\phi(0),
\end{eqnarray*}
where $\phi(\cdot)$   denotes the pdf of the standard normal distribution. By a Taylor expansion, we infer 
\begin{eqnarray*}
  \lefteqn{\mathbb E|\gamma_{n,1}(\tilde{\mathbf x})Z_1+\delta_n(\tilde{\mathbf x})|- \mathbb E|\gamma_{n,1}(\tilde{\mathbf x})Z_1|=2\delta_n(\tilde{\mathbf x})\left( \Phi(0)+\frac{\delta_n(\tilde{\mathbf x})}{\gamma_{n,1}(\tilde{\mathbf x})}\phi(0) +
  O\left(\frac{\delta_n(\tilde{\mathbf x})}{\gamma_{n,1}(\tilde{\mathbf x})}\right)^2\right)- \delta_n(\tilde{\mathbf x})}\\
  &&\qquad\qquad\qquad\qquad+ 2\gamma_{n,1}(\tilde{\mathbf x})\left( \sum_{k=0}^{3}\frac{1}{k!}\left(\frac{\delta_n(\tilde{\mathbf x})}{\gamma_{n,1}(\tilde{\mathbf x})}\right)^k \phi^{(k)}(0) +o\left\{\left(\frac{\delta_n(\tilde{\mathbf x})}{\gamma_{n,1}(\tilde{\mathbf x})}\right)^3\right\}\right) -2\gamma_{n,1}(\tilde{\mathbf x})\phi(0).
\end{eqnarray*}
By Lemma \ref{lem1}, Lemma \ref{lem4} and (\ref{var1}) we obtain, as $n\rightarrow \infty$, $$\lim_{n\rightarrow \infty} {\sigma}^{-1}\left(\sqrt{n}\mathbb E U_{n,1}(\tilde{\C})-\int_{\tilde{\C}} \mathbb E|\gamma_{n,1}(\tilde{\mathbf x})Z_1|\, d\tilde{\mathbf x}\right)=\eta(\delta), $$ this when combined with  Lemma \ref{lem3} and Lemma \ref{Lemma12} implies that $\sqrt{n}{\sigma}^{-1}\left(\mathbb EV_{n}(\tilde{\C})-n^{-1/2}a_{n}(\tilde{\C})\right)$ converges in probability to $\eta(\delta)$, as $n\rightarrow \infty$,  which completes the proof of Theorem \ref{th4}}.\hfill $\blacksquare$

\subsection*{Conflicts of interest} The authors declare that they have no known competing financial interests or personal relationships that could have appeared to influence the work reported in this paper.

\section*{Acknowledgments}
The authors would like to thank the Editor-in-Chief, an Associate-Editor, and two anonymous referees for their constructive remarks, which resulted in a substantial improvement of the work's original form and a more sharply focused presentation.

\section*{Supplementary material}
The supplementary material in \cite{L1testBBD}
contains the proofs of all technical lemmas.

\renewcommand{\thesection}{B}

\section{Technical lemmas}

\noindent The following lemma 
will be crucial for the proof of Theorem~\ref{thm1}.

\begin{lemma}\label{lem1}
	Let $g(\cdot)$ be a Lebesgue density function on $\R^{s}$, $s\geq 1,$ and $\HH$ be a finite class of bounded real-valued functions $K(\cdot)$ with compact support.
	Then, for each $K \in\HH$, we have 
	\begin{eqnarray}\label{lem1eqn1}
	|K_{h_n}\ast g(\mathbf z)-J(K)g(\mathbf z)|\rightarrow 0	\quad\mbox{as } h_n\rightarrow 0\mbox{ for almost all } \mathbf z\in\R^{s},
	\end{eqnarray}
	where
	\begin{eqnarray*}
		J(K)=\int_{\R^s}K(\mathbf u)\,d\mathbf u \, \mbox{ and } \,
		K_{h_n}\ast g(\mathbf z):= h_n^{-s}\int_{\R^s}K\Big(\frac{\mathbf z-\mathbf v}{h_n}\Big)g(\mathbf v)\,d\mathbf v.	
	\end{eqnarray*}
	Moreover, for all $0<\varepsilon<1$, there exist $M$, $\nu>0$ and a Borel set $C$ of finite Lebesgue measure such
	that
	\begin{eqnarray} \label{lem1eqn2}
	C\!\subset\!\lbrack -M+\nu, M-\nu\rbrack^{s}, \int_{ \R^{s}\setminus \lbrack -M, M\rbrack^{s}
	} g(\mathbf v)\,d\mathbf v=\mu>0, 	\int_{C} g(\mathbf v)\,d\mathbf v>1-\varepsilon, 
	\end{eqnarray}
	\begin{eqnarray}\label{lem1eqn3}
	g(\cdot) \mbox{ is bounded, continuous and bounded away from } 0
	\mbox{ on } C
	\end{eqnarray}
	and, for each $K \in\HH$,
	\begin{eqnarray}\label{lem1eqn4}
	\sup_{z\in C}| K_{h_n}\ast g(\mathbf z)-J(K)g(\mathbf z)|\rightarrow 0 \qquad \mbox{ as } h_n\rightarrow 0.
	\end{eqnarray}
\end{lemma}
\noindent{\bf{Proof of Lemma \ref{lem1}.}} 
The first statement is just a simple case of Theorem 3 in Chapter 2 of \cite{Devroye12}. The other statements can be proved in the same way as in Lemma 6.1 of \cite{mason2001}.
Since $g(\cdot)$ is Lebesgue integrable, the integral $\int_{ \R^{s}\smallsetminus \lbrack -M, M\rbrack^{s}
} g(\mathbf v)\,d\mathbf v$ is continuous in $M$ and converges to zero as $M\rightarrow\infty$. We may find 
$M>0$ and $\nu>0$ so that 
\begin{eqnarray*}
	\int_{ \R^{s}\setminus \lbrack -M, M\rbrack^{s}
	} g(\mathbf v)\,d\mathbf v=\frac{\varepsilon}{8} \quad \mbox{ and } \quad \int_{ \R^{s}\setminus \lbrack -M+\nu, M-\nu\rbrack^{s}
	} g(\mathbf v)\,d\mathbf v=\frac{\varepsilon}{4}.
\end{eqnarray*}
The existence of the desired set $C \subset \lbrack -M + \nu, M -\nu\rbrack^{s}$ can be inferred from Lusin's theorem followed by Egorov's theorem, for instance, see \cite{dudley}, Theorems $7.5.1$ and $7.5.2$. 
\hfill$\blacksquare$\\
~

\noindent The next lemma provides a truncation step and  will play an instrumental role in the sequel.

\begin{lemma}\label{lem2}
	Assume that assumption \textnormal{({\bf A$_1$})} holds. If $h_n\rightarrow 0$ and $nh_n^{d}\rightarrow\infty$, then, for all Borel set $A\in\R^{d}$, 
	\begin{eqnarray*}
		\limsup_{n\rightarrow \infty} n \mathrm{Var}(V_{n}(A))
		&\leq& 64 \int_{A} f(\tilde{\mathbf x})\,d\tilde{\mathbf x}.
	\end{eqnarray*}
\end{lemma}
\noindent{\bf{Proof of Lemma \ref{lem2}.}} 
Let $\mathbb E_0$ be the expectation $\mathbb E$ and let, for each $ i=1,\ldots, n,$  $\mathbb E_i$ be the conditional expectation given $(\mathbb X_{1}^{1},\ldots,\mathbb X_{i}^{1})$. Set, for each $i=1,\ldots,n,$
\begin{eqnarray*}
	\Upsilon_i:=\mathbb E_i V_{n}(A)-\mathbb E_{i-1} V_{n}(A)
	\, \mbox{ and } \, 
	\Phi_i:=\mathbb E_i\big( V_{n}(A)-V_{n,-i}(A)\big),
\end{eqnarray*}
where
\begin{eqnarray*}
	V_{n,-i}(A)&:=&\int_{A}\left|\frac{1}{nh_n^{d}} \sum_{k\neq i}^{n}\prod_{l=1}^{p}K_{l}\left(\frac{\mathbf x_l-\mathbb X_{k}^{l}}{h_n}\right) \right.\\
	&&\left.\hspace{2cm}-\frac{1}{n^{p}h_n^{d}}\sum_{k\neq i}^{n}K_{1}\left(\frac{\mathbf x_1-\mathbb X^{1}_k}{h_n}\right) \prod_{l=2}^{p}\sum_{k=1}^{n}K_{l}\left(\frac{\mathbf x_l-\mathbb X_{k}^{l}}{h_n}\right)\right| \,d\tilde{\mathbf x}.
\end{eqnarray*}
Clearly, we have 
\begin{eqnarray}\label{lem21}
\sum_{i=1}^n \Upsilon_i=V_{n}(A)-\mathbb E V_{n}(A)
\quad \mbox{ and }\quad \mathbb E_{i-1}\Upsilon_i=0, \quad i=1,\ldots,n.
\end{eqnarray} 
Also observe that, for each $i=1,\ldots,n,$ $\Upsilon_i=\Phi_i-\mathbb E_{i-1}\Phi_i$
and 
\begin{eqnarray}\label{PHI}
|\Phi_i|&\leq&\mathbb E_i \big|V_{n}(A)-V_{n,-i}(A)\big|\nonumber\\
&\leq& \frac{2}{nh_n^{d}}\int_{A} K_{1}\Bigg(\frac{\mathbf x_1-\mathbb X_{i}^{1}}{h_n}\Bigg)\prod_{l=2}^{p}\mathbb EK_{l}\Bigg(\frac{\mathbf x_l-\mathbb X^{l}}{h_n}\Bigg)\,d\tilde{\mathbf x}\nonumber\\
&=:&\Psi_i.
\end{eqnarray}
Hence, for each $i=1,\ldots,n,$
\begin{eqnarray}\label{lem212}
|\Upsilon_i	| \leq |\Phi_i|+\mathbb E_{i-1} |\Phi_i|
\leq \Psi_i+\mathbb E\Psi_i:=\xi_i.
\end{eqnarray}
For each $i=1,\ldots,n,$ let $\bar\xi_i$ be independent copy of $\xi_i$ and let
$(\varepsilon_i)$ be an i.i.d. sequence of Rademacher variables ($\mathbb P(\varepsilon_i=1)=\mathbb P(\varepsilon_i=-1)=1/2$), independent of the sequence $(\xi_i,\bar\xi_i)$. Using (\ref{lem21})-(\ref{lem212}), the convexity of $y=x^2$ and the fact that, for each $i=1,\ldots,n$,  $\xi_i$ is independent of $(\mathbb X^{1}_{1},\ldots,\mathbb X^{1}_{i-1})$, it follows from Theorem $3.1$ in \cite{Berg} that 
\begin{eqnarray}\label{lem22}
\mathbb E\left(\sum_{i=1}^n \Upsilon_i\right)^2
\leq \mathbb E \left(\sum_{i=1}^n\varepsilon_i(\xi_i+\bar\xi_i)\right)^{2}
\leq 4n \mathbb E (\Psi_1+E\Psi_1)^{2}
\leq 16n\mathbb E \Psi_1^2.
\end{eqnarray}
Combining (\ref{lem21})-(\ref{PHI}) with (\ref{lem22}), we obtain 
\begin{eqnarray*}
	n \mathbb E\left[V_{n}(A)-\mathbb E V_{n}(A)\right]^2&\leq& 64 \int_{A} \mathbb Ef_{n}(\tilde{\mathbf x})\, d\tilde{\mathbf x}. 
\end{eqnarray*}
By Theorem $1$ in Chapter $2$ of \cite{Devroye12}, we infer that, as $n\rightarrow\infty$, 
\begin{eqnarray*}
	\int_{A}\big|\mathbb E f_{n}(\tilde{\mathbf x})-f(\tilde{\mathbf x})\big| \,d\tilde{\mathbf x}\rightarrow0, 
\end{eqnarray*}
which completes the proof of the lemma. \hfill $\blacksquare$\\
~

\noindent Choose, for each $l=1,\ldots,p$, any bounded Borel set $C_l\subset\R^{d_{l}}$ satisfying (\ref{lem1eqn3}) and (\ref{lem1eqn4}) of Lemma \ref{lem1} with $g(\cdot)=f_l(\cdot)$ and $\HH=\HH_l:=\{ K_l, K_l^{2}, K_l^{3}\}$. Clearly, for each $l=1,\ldots,p$, for each $K\in\HH_l$, for all large enough $n$ uniformly in $\mathbf z_l \in C_l$ and for some constant
$D_0 > 0$,
\begin{eqnarray}\label{majora1}
\sup_{\mathbf z_l\in C_{l} }\int_{\R^{d_l}} K\bigg(\frac{\mathbf z_l-\mathbf u_l}{h_{n}}\bigg)f_{l}(\mathbf \mathbf u_l)\,d\mathbf u_l\leq D_0 h_{n}^{d_l}.
\end{eqnarray}
Let us recall, for each $a\in\{0, 1\}$, 
\begin{eqnarray*}
	\Gamma_{n,a}(\tilde{\mathbf x}):= f_{n}(\tilde{\mathbf x})-
	a\sum_{l=1}^{p} f_{n,l}(\mathbf x_l)\prod_{j\neq l }\mathbb E f_{n,j}(\mathbf x_j)
	+(pa-1)\prod_{l=1}^{p}\mathbb E f_{n,l}(\mathbf x_l)
\end{eqnarray*}
and set, for $\bar{C}= C_1\times \cdots\times C_p$,
\begin{eqnarray*}
	U_{n,a}(\bar{C}):=\int_{\bar{C}} |\Gamma_{n,a}(\tilde{\mathbf x})|\,d\tilde{\mathbf x}.
\end{eqnarray*}
The next lemma shows that $\sqrt{n}\big(V_{n}(\bar{C})-\mathbb E V_{n}(\bar{C})\big) $ is asymptotically equivalent in probability to $\sqrt{n}\big(U_{n,1}(\bar{C})-\mathbb E U_{n,1}(\bar{C})\big) $.

\begin{lemma}\label{lem3}
	Suppose that, for each $l=1,\ldots,p$, $C_l$ satisfies \textnormal{(\ref{lem1eqn3})} and \textnormal{(\ref{lem1eqn4})} of Lemma \ref{lem1} with $g(\cdot)=f_l(\cdot)$ and $\HH(\cdot)=\HH_l(\cdot)$. If
	$h_n\rightarrow 0$ and $nh_n^{d}\rightarrow\infty$, then 
	\begin{eqnarray*}
		\lim_{n\rightarrow\infty} \sqrt{n}\mathbb E\big| V_{n}(\bar{C})-U_{n,1}(\bar{C})\big|=0.
	\end{eqnarray*}
\end{lemma}
\noindent{\bf{Proof of Lemma \ref{lem3}.}} 
Notice that
\begin{eqnarray*} &&\Big|\big|f_{n}(\tilde{\mathbf x})-\prod_{l=1}^{p}f_{n,l}(\mathbf x_l)\big|
	-\big|
	\Gamma_{n,1}(\tilde{\mathbf x})\big|\Big|
	\\&&\hspace{2.5cm}\leq
	\sum_{\{I\subsetneq I_p,|I_p\setminus I|\geq2\}}\prod_{l\in I }\mathbb Ef_{n,l}(\mathbf x_l)\prod_{j\in I_p\setminus I}\big|f_{n,j}(\mathbf x_j)-\mathbb E f_{n,j}(\mathbf x_j)\big|,
\end{eqnarray*}
where $I_p =\{1,\ldots,p\}$, $|I|$ denotes the cardinality of the set $I$ and the empty product is defined to be $1$. Clearly, we have 
\begin{eqnarray}\label{maj1lem3}
\mathbb E\big| V_{n}(\bar{C})-U_{n,1}(\bar{C})\big|&\leq&
\sum_{\{I\subsetneq I_p,|I_p\setminus I|\geq2\}}\prod_{l\in I_p\setminus I}\int_{C_l}\mathbb E\big|f_{n,l}(\mathbf x_l)-\mathbb E f_{n,l}(\mathbf x_l)\big|\,d\mathbf x_l\nonumber\\&\leq &
\sum_{\{I\subsetneq I_p,|I_p\setminus I|\geq2\}}\prod_{l\in I_p\setminus I}\int_{C_l}\big(\mathbb E\lbrack f_{n,l}(\mathbf x_l)-\mathbb E f_{n,l}(\mathbf x_l)\rbrack^{2}\big)^{1/2}\,d\mathbf x_l.
\end{eqnarray}
Now, observe that, for each $l=1, \ldots, p$ and each $\mathbf x_l\in C_l$, 
\begin{eqnarray*}\label{maj2lem3}
	\mathbb E \left[f_{n,l}(\mathbf x_l)-\mathbb E f_{n,l}(\mathbf x_l)\right]^{2}&\leq& 
	\frac{1}{nh_{n}^{2d_{l}}}\int_{\R^{d_l}} K_{l}^{2}\bigg(\frac{\mathbf x_l-\mathbf u_l}{h_{n}}\bigg)f_{l}(\mathbf u_l)\,d\mathbf u_l.
\end{eqnarray*}
Therefore, by (\ref{majora1}), for each $l=1, \ldots, p$, we have 
\begin{eqnarray}\label{maj3lem3}
\sup_{\mathbf x_l\in C_l}\mathbb E \left[f_{n,l}(\mathbf x_l)-\mathbb E f_{n,l}(\mathbf x_l)\right]^{2}=O\Bigg(\frac{1}{nh_{n}^{d_{l}}}\Bigg).
\end{eqnarray}
Inequality (\ref{maj1lem3}) together with (\ref{maj3lem3}) give, for some constant $M_{0}>0$,
\begin{eqnarray*}
	\sqrt{n}\mathbb E\big| V_{n}(\bar{C})-U_{n,1}(\bar{C})\big|&\leq& M_{0} \sum_{\{I\subsetneq I_p,|I_p\setminus I|\geq2\}} 
	\frac{1}{\displaystyle 
		n^{(|I_p\setminus I|-1)/2}
		\prod_{l\in I_p\setminus I}h_{n}^{d_l/2}}\nonumber\\&\leq&
	\frac{M_{0}}{\sqrt{nh_{n}^{d}}}\sum_{\{I\subsetneq I_p,|I_p\setminus I|\geq2\}}\prod_{l\in I} h_{n}^{d_l/2}. 
\end{eqnarray*}
This completes the proof of Lemma \ref{lem3}. \hfill $\blacksquare$
\\
~

\noindent The next step in the proof of the asymptotic normality theorem is what is known as Poissonization. Let $\eta$ be a Poisson random variable with mean $n$, independent
of $\tilde{\mathbf X},\tilde{\mathbf X}_1,\tilde{\mathbf X}_2, \ldots$, and set
\begin{eqnarray*}
	f_{\eta}(\tilde{\mathbf x})=
	\frac{1}{nh_{n}^{d}}
	\sum\limits_{i=1}^{\eta}\prod_{l=1}^{p}
	K_{l}\Bigg(\frac{\mathbf x_l-\mathbf X^{l}_{i}}{h_{n}}\Bigg) \,\, \mbox{ and } \,\,
	f_{\eta,l}(\mathbf x_l)=\frac{1}{nh_{n}^{d_{l}}}\sum\limits_{i=1}^{\eta}
	K_{l}\Bigg(\frac{\mathbf x_l-\mathbf X^{l}_{i}}{h_{n}}\Bigg),
\end{eqnarray*}
where the empty sum is defined to be zero. For each $l=1,\ldots,p$, notice that 
\begin{eqnarray}\label{key1}
\mathbb E
f_{\eta,l}(\mathbf x_l)&=&\mathbb E
f_{n,l}(\mathbf x_l) =h_n^{-d_{l}}\int_{ {\R}^{d_l}} K_{l}\left(\frac{\mathbf x_l-\mathbf u_l}{h_n}\right)f_{l}( \mathbf u_l)\,d\mathbf u_l
\end{eqnarray}
and
\begin{eqnarray} \label{key4}
v_{n,l}(\mathbf x_l)&:=& n\mathrm{Var}\big(f_{\eta,l}(\mathbf x_l)\big)=h_n^{-2d_l}\int_{ {\R}^{d_l}} K_{l}^{2}\left(\frac{\mathbf x_l-\mathbf u_l}{h_n}\right)f_{l}( \mathbf u_l)\,d\mathbf u_l.
\end{eqnarray}
Set, for each $a\in\{0,1\}$,
\begin{eqnarray*}
	\Gamma_{\eta,a}(\tilde{\mathbf x})&:=& f_{\eta}(\tilde{\mathbf x})-
	a\sum_{l=1}^{p} f_{\eta,l}(\mathbf x_l)\prod_{j\neq l }\mathbb E f_{n,j}(\mathbf x_j)
	+(pa-1)\prod_{l=1}^{p}\mathbb E f_{n,l}(\mathbf x_l).
\end{eqnarray*}
Clearly, for each $a\in\{0,1\}$, we have 
\begin{eqnarray}\label{eqnvareta}
k_{n,a}(\tilde{\mathbf x})&:= & n\mathrm{Var}\big(\Gamma_{\eta,a}(\tilde{\mathbf x})\big)\nonumber\\
&=&
\prod_{l=1}^{p}v_{n,l}(\mathbf x_l)-a\sum_{l=1}^{p}v_{n,l}(\mathbf x_l)\prod_{j\neq l}^{p}\big(\mathbb E f_{n,j}(\mathbf x_j)\big)^{2}
+p(p-1)a\prod_{l=1 }^{p}\big(\mathbb E f_{n,l}(\mathbf x_l)\big)^{2}
\end{eqnarray}
and \begin{eqnarray*}
	\mathcal{L}_{n,a}(\tilde{\mathbf x})&:= & n\mathrm{Var}\big(\Gamma_{n,a}(\tilde{\mathbf x})\big)\\
	&=&
	\prod_{l=1}^{p}v_{n,l}(\mathbf x_l)-a\sum_{l=1}^{p}v_{n,l}(\mathbf x_l)\prod_{j\neq l}^{p}\big(\mathbb E f_{n,j}(\mathbf x_j)\big)^{2}+(pa-1)\prod_{l=1 }^{p}\big(\mathbb E f_{n,l}(\mathbf x_l)\big)^{2}.
\end{eqnarray*}
We now observe that
\begin{eqnarray*}\label{maj56}
	&&\left|h_{n}^{d}k_{n,a}(\tilde{\mathbf x})- \tilde{\bf{K}} \prod_{l=1}^{p}f_{l}(\mathbf x_l)\right|\nonumber\\
	&&\hspace{1cm}\leq 
	\sum_{I\subsetneq I_p}\prod_{l\in I }\|K_l\|_{2}^{2}f_{l}(\mathbf x_l
	)\prod_{j\in I_p\setminus I}
	\big|h_{n}^{d_j}v_{n,j}(\mathbf x_j)-\|K_{j}\|_{2}^{2}f_{j}(\mathbf x_j) \big|
	\nonumber\\&& 
	\hspace{1.2cm}+ ah_{n}^{d}\left\{\sum_{l=1}^{p}v_{n,l}(\mathbf x_l)\prod_{j\neq l}^{p}\big(\mathbb E f_{n,j}(\mathbf x_j)\big)^{2}+p(p-1)\prod_{l=1 }^{p}\big(\mathbb E f_{n,l}(\mathbf x_l)\big)^{2}\right\}.
\end{eqnarray*}
Combining this last observation with
(\ref{key1})-(\ref{key4}) and the fact that, for $l=1,\ldots,p$, $C_l$ satisfies (\ref{lem1eqn3}) and (\ref{lem1eqn4}) of Lemma \ref{lem1} with $g(\cdot)=f_l(\cdot)$ and $\HH(\cdot)=\HH_l(\cdot)$, we obtain
\begin{eqnarray}\label{var1}
\lim_{n\rightarrow \infty} \sup_{\tilde{\mathbf x}\in \bar{C} }\left|h_{n}^{d}k_{n,a}(\tilde{\mathbf x})-\tilde{\bf{K}} \prod_{l=1}^{p}f_{l}(\mathbf x_l)\right|=0.
\end{eqnarray}
Thus, for all large enough $n$ uniformly in $\tilde{\mathbf x}\in \bar{C}$ and for some constants $D_1>0$ and $D_2>0$, we have 
\begin{eqnarray}\label{l4eqn2}
D_1\leq h_{n}^{d}k_{n,a}(\tilde{\mathbf x})\leq D_2.
\end{eqnarray} 
Therefore, by (\ref{majora1}) and (\ref{l4eqn2}), we get 
\begin{eqnarray}\label{egali}
&&\sup_{\tilde{\mathbf x}\in \bar{C}}\left| \sqrt{\mathcal{L}_{n,a}(\tilde{\mathbf x})}- \sqrt{k_{n,a}(\tilde{\mathbf x})}\right|\nonumber \leq |2pa-p^{2}a-1| \sup_{\tilde{\mathbf x}\in\bar{C}}\frac{\displaystyle\prod_{l=1}^{p} \big(\mathbb Ef_{n,l}(\mathbf x_l)\big)^{2} }{\sqrt{ k_{n,a}(\tilde{\mathbf x})}}=O\left( \sqrt{h_{n}^{d}}\right).
\end{eqnarray}
Set, for each $a\in\{0,1\}$,
\begin{eqnarray*}
	U_{\eta,a}(\bar{C}):=\int_{\bar{C}} |\Gamma_{\eta,a}(\tilde{\mathbf x})|\,d\tilde{\mathbf x}.
\end{eqnarray*}
The next lemma shows that, for each $a\in\{0,1\}$, $\sqrt{n}\mathbb E U_{n,a}(\bar{C})$ is asymptotically equivalent to $\sqrt{n}\mathbb E U_{\eta,a}(\bar{C})$.
\begin{lemma}\label{lem4}
	Suppose that, for each $l=1,\ldots,p$, $C_l$ satisfies \textnormal{(\ref{lem1eqn3})} and \textnormal{(\ref{lem1eqn4})} of Lemma \ref{lem1} with $g(\cdot)=f_l(\cdot)$ and $\HH(\cdot)=\HH_l(\cdot)$. If
	$h_n\rightarrow 0$ and $\sqrt{n}h_n^{d}\rightarrow\infty$, as $n\rightarrow\infty$, then, for each $a\in\{0,1\}$, 
	\begin{eqnarray}\label{keyeqn1}
	\lim_{n\rightarrow \infty}
	\sqrt{n}\mathbb E U_{\eta,a}(\bar{C})
	-\mathbb E|Z_1|\int_{\bar{C}}
	\sqrt{	k_{n,a}(\tilde{\mathbf x})}\,d\tilde{\mathbf x}=0,
	\end{eqnarray}
	and
	\begin{eqnarray}\label{keyeqn2}
	\lim_{n\rightarrow \infty} \sqrt{n}\mathbb E U_{n,a}(\bar{C})
	-\mathbb E|Z_1|\int_{\bar{C}}
	\sqrt{k_{n,a}(\tilde{\mathbf x})}\,d\tilde{\mathbf x}=0.
	\end{eqnarray}
\end{lemma}
\noindent{\bf{Proof of Lemma \ref{lem4}.}} Let $\eta_1$ denote a Poisson random variable with mean $1$, independent of $\tilde{\mathbf X}_1,\tilde{\mathbf X}_2,\ldots$, and set
\begin{eqnarray*}
	\Y_{a,n}(\tilde{\mathbf x})&=&\frac{1}{h_{n}^{d}\sqrt{ k_{n,a}(\tilde{\mathbf x})}}\left[ \sum_{i\leq \eta_1}\bigg\{\prod_{l=1}^{p}K_{l}\bigg(\frac{\mathbf x_l-\mathbb X_{i}^{l}}{h_n}\bigg) 
	\right.\\&&\hspace{3.1cm}-a \sum\limits_{l=1}^{p} K_{l}\bigg(\frac{\mathbf x_l-\mathbb X_{i}^{l}}{h_n}\bigg)\prod_{j\neq l} \mathbb EK_{j}\bigg(\frac{\mathbf x_j-\mathbb X^{j}}{h_n} \bigg)\bigg\}
	\\&&\hspace{5.3cm}+ \left.(pa-1)\prod_{l=1}^{p} \mathbb EK_{l}\bigg(\frac{\mathbf x_l-\mathbb X^{l}}{h_n}\bigg) \right]	.
\end{eqnarray*}
Observe that $\mathrm{Var} [\Y_{a,n}(\tilde{\mathbf x})]=1$. Using the statement $(2.8)$ of Lemma $2.3$ of \cite{mason2001} and $c_r$-inequality, we get 
\begin{eqnarray*}\label{l4eqn1}
	\mathbb E|\Y_{a,n}(\tilde{\mathbf x})|^{3}\leq D_3\frac{\displaystyle h_{n}^{-\frac{3}{2}d} \prod_{l=1}^{p}\mathbb EK_{l}^{3}\left(\frac{\mathbf x_l-\mathbb X^{l}}{h_n}\right)  }{\displaystyle \big(h_{n}^{d} k_{n,a}(\tilde{\mathbf x})\big)^{3/2} },
\end{eqnarray*}
where $D_3>0$ is a constant. Using (\ref{majora1}) with (\ref{l4eqn2}), we obtain, for all large enough $n$ and for some constant $D_4>0$,
\begin{eqnarray}\label{l4eqn4}
\sup_{\tilde{\mathbf x}\in \bar{C}}\mathbb E|\Y_{a,n}(\tilde{\mathbf x})|^{3}\leq D_4h_{n}^{-\frac{d}{2}}.
\end{eqnarray}
Let $\Y_{a,n}^{(1)}(\tilde{\mathbf x}),\ldots,\Y_{a,n}^{(n)}(\tilde{\mathbf x})$ be independent copies of $\Y_{a,n}(\tilde{\mathbf x})$. Clearly,  we have
\begin{eqnarray*}\label{l4eqn5}
	T_{\eta,a}(\tilde{\mathbf x}):=	\frac{\sqrt{n}\Gamma_{\eta,a}(\tilde{\mathbf x})}{\sqrt{k_{n,a}(\tilde{\mathbf x})}}\stackrel{\D}{=} \frac{1}{\sqrt{n}} \sum_{i=1}^{n}\Y_{a,n}^{(i)}(\tilde{\mathbf x}),
\end{eqnarray*}
where $\stackrel{\D}{=}$ stands for equality in distribution. Therefore, by (\ref{equa6.18Fact6.1}) we have 
\begin{eqnarray}\label{l4eqn6}
\sup_{\tilde{\mathbf x}\in \bar{C}} \left| \sqrt{n}\frac{\mathbb E |\Gamma_{\eta,a}(\tilde{\mathbf x})| }{\sqrt{k_{n,a}(\tilde{\mathbf x})}} - \mathbb E\left| Z_1\right| \right| \leq \frac{D_5}{\sqrt{n}} \sup_{\tilde{\mathbf x}\in \bar{C}}\mathbb E|\Y_{a,n}(\tilde{\mathbf x})|^{3},
\end{eqnarray}
where $D_5$ is a universal positive constant. Now, by (\ref{l4eqn2}), in combination with (\ref{l4eqn4}) and (\ref{l4eqn6}), we get 
\begin{eqnarray*}
	\Bigg|\int_{ \bar{C}}\left\{\sqrt{n}\mathbb E|\Gamma_{\eta,a}(\tilde{\mathbf x})|-\mathbb E|Z_1|
	\sqrt{k_{n,a}(\tilde{\mathbf x})}\right\}\,d\tilde{\mathbf x}\Bigg|
	= O\Bigg( \frac{1}{\sqrt{n}h_{n}^{d}}\Bigg). 
\end{eqnarray*}
Similarly, we obtain, by the statement (\ref{equa6.18Fact6.1}),
\begin{eqnarray*}\label{l4eqn7}
	&&\Bigg|\int_{\bar{C}}\left\{\sqrt{n}\mathbb E|\Gamma_{n,a}(\tilde{\mathbf x})|-\mathbb E|Z_1|
	\sqrt{\mathcal{L}_{n,a}(\tilde{\mathbf x})}\right\}\,d\tilde{\mathbf x}	\Bigg| 
	=O\Bigg( \frac{1}{\sqrt{n}h_{n}^{d}}\Bigg),
\end{eqnarray*}
which by (\ref{egali}) implies
\begin{eqnarray*}
	\Bigg|\int_{\bar{C}}\left\{\sqrt{n}\mathbb E|\Gamma_{n,a}(\tilde{\mathbf x})|-
	\mathbb E|Z_1|\sqrt{k_{n,a}(\tilde{\mathbf x})}\right\}\,d\tilde{\mathbf x}\Bigg|=O\Bigg( \frac{1}{\sqrt{n}h_{n}^{d}}+\sqrt{h_{n}^{d}} \Bigg).
\end{eqnarray*}
This completes the proof.
\hfill $\blacksquare$
\\
~

\begin{lemma}\label{lem4b}
	Whenever $h_n\rightarrow 0$ and, for $l=1,\ldots,p$, $C_l$ satisfies \textnormal{(\ref{lem1eqn3})} of Lemma \ref{lem1} with $g(\cdot)=f_l(\cdot)$, we have
	\begin{eqnarray}\label{l4beqn1}
	&& \mathds{1}_{\bar{C}}(\tilde{\mathbf x}+h_n \tilde{\mathbf t} )\mbox{ converges in measure to }\mathds{1}_{\bar{C}}(\tilde{\mathbf x})=1 \mbox{ on } \bar{C}\times\tilde{\B}.
	\end{eqnarray}
\end{lemma}
\noindent{\bf{Proof of Lemma \ref{lem4b}.}} 
Notice that 
\begin{eqnarray*}
	&&\int_{\bar{C}}\int_{\tilde{\B}} \mathds{1}_{\bar{C}}(\tilde{\mathbf x}+h_n \tilde{\mathbf t})\,	d\tilde{\mathbf t} \,	 
	d\tilde{\mathbf x} =\frac{1}{h_n^{d}}	\int_{\bar{C}}\int_{\bar{C}}
	\mathds{1}_{ \tilde{\B}}\Bigg( \frac{\tilde{\mathbf x}-\tilde{\mathbf y}}{h_n}\Bigg) \,d\tilde{\mathbf y}\, d\tilde{\mathbf x}.
\end{eqnarray*} 
Now, by (\ref{lem1eqn1}), applied to
$K(\tilde{\mathbf t})=\dfrac{1}{ \lambda(\tilde{ \B})}\, \mathds{1}_{ \tilde{ \B}}(\tilde{\mathbf t})$ and $f(\tilde{\mathbf x})=\dfrac{1}{\lambda( \bar{C})} \mathds{1}_{\bar{C}}(\tilde{\mathbf x})$, for almost every $\tilde{\mathbf x}$, we have 
\begin{eqnarray*}
	&&
	\frac{1}{h_n^{d}}	\int_{\R^d}\mathds{1}_{\bar{C}}(\tilde{\mathbf y})
	\mathds{1}_{\tilde{ \B}}\Bigg( \frac{\tilde{\mathbf x}-\tilde{\mathbf y}}{h_n}\Bigg)\, d\tilde{\mathbf y} 
	\rightarrow \lambda( \tilde{ \B}) \mathds{1}_{\bar{C}}(\tilde{\mathbf x}).
\end{eqnarray*} 
Thus, by the dominated convergence theorem, we infer that 
\begin{eqnarray*}
	&&\frac{1}{h_n^{d}}	\int_{\bar{C}}\int_{\bar{C}}
	\mathds{1}_{ \tilde{\B}}\Bigg( \frac{\tilde{\mathbf x}-\tilde{\mathbf y}}{h_n}\Bigg) \,d\tilde{\mathbf y}\, d\tilde{\mathbf x}\rightarrow \lambda( \tilde{ \B})\,\lambda( \bar{C}),
\end{eqnarray*}
which, in other words, says
\begin{eqnarray*}
	&& (\lambda \times \lambda) \big \{ (\tilde{\mathbf x},\tilde{\mathbf t})\in \bar{C} \times \tilde{\B} : 1-\mathds{1}_{\bar{C}}(\tilde{\mathbf x}+h_n\tilde{\mathbf t})\neq 0 \big\}  \\
	&& \hspace{1cm}= \lambda( \tilde{ \B})\,\lambda( \bar{C})- \int_{C}\int_{\tilde{\B}} \mathds{1}_{\bar{C}}(\tilde{\mathbf x}+h_n\tilde{\mathbf t} ) \,	d\tilde{\mathbf t} \, d\tilde{\mathbf x}	 
	\rightarrow 0,
\end{eqnarray*}
yielding to (\ref{l4beqn1}). Hence the proof of Lemma \ref{lem4b} is complete.
\hfill$\blacksquare$
\\
~
\begin{lemma}\label{lem5}
	Suppose that, for each $l=1,\ldots,p$, $C_l$ satisfies \textnormal{(\ref{lem1eqn3})} and \textnormal{(\ref{lem1eqn4})} of Lemma \ref{lem1} with $g(\cdot)=f_l(\cdot)$ and $\HH(\cdot)=\HH_l(\cdot)$. If
	$h_n\rightarrow 0$ and $nh_n^{3d}\rightarrow\infty$, as $n\rightarrow\infty$, then, for all $a,b\in\{0,1\}$, we have 
	
	\begin{eqnarray*}
		\lim_{n\rightarrow\infty} n\mathrm{Cov}\big(U_{\eta,a}(\bar{C}), U_{\eta,b}(\bar{C})\big)=\sigma^2\prod_{l=1}^{p}\int_{C_l}f_l(\mathbf x_l)\,d\mathbf x_l.
	\end{eqnarray*}
\end{lemma}

\noindent{\bf{Proof of Lemma \ref{lem5}.}} For any Borel subset $A$ of $\R^{2d}$, set
\begin{eqnarray*}
	&&\sigma_{n,a,b}(\bar{C},A):=\int_{\bar{C}}\int_{\bar{C}}	 \mathrm{Cov}\big(|T_{\eta,a}(\tilde{\mathbf x})|, |T_{\eta,b}(\tilde{\mathbf y})|\big)\mathds{1}_{A}(\tilde{\mathbf x},
	\tilde{\mathbf y})
	\sqrt{k_{n,a}(\tilde{\mathbf x})k_{n,b}(\tilde{\mathbf y})}
	\,d\tilde{\mathbf x}\, d\tilde{\mathbf y}.
\end{eqnarray*} 
Notice that 
\begin{eqnarray*}
	n\mathrm{Cov}\big(U_{\eta,a}(\bar{C}), U_{\eta,b}(\bar{C})\big)&=&n \int_{\bar{C}}\int_{\bar{C}}	 \mathrm{Cov}\big(|\Gamma_{\eta,a}(\tilde{\mathbf x})|, |\Gamma_{\eta,b}(\tilde{\mathbf y})|\big)\,
	\,d\tilde{\mathbf x}\, d\tilde{\mathbf y}\\
	&=&\int_{\bar{C}}\int_{\bar{C}}	 \mathrm{Cov}\big(|T_{\eta,a}(\tilde{\mathbf x})|, |T_{\eta,b}(\tilde{\mathbf y})|\big)
	\sqrt{k_{n,a}(\tilde{\mathbf x})k_{n,b}(\tilde{\mathbf y})}
	\,d\tilde{\mathbf x}\, d\tilde{\mathbf y}\\	 &=&\sum_{I\subset I_p}\sigma_{n,a,b}(\bar{C},A_{n}(I)),
\end{eqnarray*}
where, for each $ I\subset I_p =\{1,\ldots,p\}$, 
\begin{eqnarray*}
	A_{n}(I):=\bigcap_{l\in I}
	\Big\{(\tilde{\mathbf x},\tilde{\mathbf y}): \|\mathbf x_l-\mathbf y_l\|\leq h_n\Big\}\bigcap\bigcap_{j\in I_p\setminus I}\Big\{(\tilde{\mathbf x},\tilde{\mathbf y}): \|\mathbf x_j-\mathbf y_j\|> h_n \Big\}.
\end{eqnarray*} 
Here, the empty intersection is defined to be $\R^{2d}$. We will show that, as $n\rightarrow\infty$,
\begin{eqnarray}\label{t1var}
\sigma_{n,a,b}(\bar{C},A_{n}(I_p))\rightarrow \sigma^{2}\prod_{l=1}^{p}\int_{C_l}f_l(\mathbf x_l)\,d\mathbf x_l
\end{eqnarray}
and 
\begin{eqnarray}\label{t2var}
\sum_{I\subsetneq I_p}\sigma_{n,a,b}(\bar{C},A_{n}(I))\rightarrow 0,
\end{eqnarray}
which will complete the proof of the lemma. First, consider (\ref{t1var}).
By (\ref{var1}), in combination with the fact that $\lambda(\bar{C})<\infty$ and
\begin{eqnarray*}
	\int_{\bar{C}}\int_{\bar{C}}	 
	\mathds{1}_{A_{n}(I_p)}(\tilde{\mathbf x},
	\tilde{\mathbf y})\,d\tilde{\mathbf x}\,d\tilde{\mathbf y}\leq \beta_0h_{n}^{d}\lambda(\bar{C}),
\end{eqnarray*}
for some constant $\beta_0>0$, we see that
\begin{eqnarray*}
	\sigma_{n,a,b}(\bar{C},A_{n}(I_p)
	)=\bar{{\sigma}}_{n,a,b}(\bar{C},A_{n}(I_p))+o(1),
\end{eqnarray*}
where
\begin{eqnarray*}
	\bar{{\sigma}}_{n,a,b}(\bar{C},A_{n}(I_p))&=&	\tilde{\bf{K}}\int_{\bar{C}}\int_{\bar{C}}\mathds{1}_{A_{n}(I_p)}(\tilde{\mathbf x},
	\tilde{\mathbf y})\mathrm{Cov}\big(|T_{\eta,a}(\tilde{\mathbf x})|,
	|T_{\eta,b}(\tilde{\mathbf y})|\big)
	\\&&\times h_{n}^{-d}	\prod_{l=1}^{p}\sqrt{f_{l}(\mathbf x_l)f_{l}(\mathbf y_l)} 		\,d\tilde{\mathbf x}\,d\tilde{\mathbf y}.
\end{eqnarray*}
Now, let
$(Z^{(a,b)}_{n,1}(\tilde{\mathbf x}),Z^{(a,b)}_{n,2}(\tilde{\mathbf y}))$, $\tilde{\mathbf x},
\tilde{\mathbf y}\in \R^{d}$, be a mean zero bivariate Gaussian process such that, for each $(\tilde{\mathbf x},
\tilde{\mathbf y})\in \R^{2d}$, $(Z^{(a,b)}_{n,1}(\tilde{\mathbf x}),Z^{(a,b)}_{n,2}(\tilde{\mathbf y}))$ and $(T_{\eta,a}(\tilde{\mathbf x}), T_{\eta,b}(\tilde{\mathbf y}))$
have the same covariance structure. In particular, we have
\begin{eqnarray*}
	\big(Z^{(a,b)}_{n,1}(\tilde{\mathbf x}),Z^{(a,b)}_{n,2}(\tilde{\mathbf y})\big)\stackrel{\D}{=}\bigg( \sqrt{1-\left(\rho_{n,a,b}(\tilde{\mathbf x},
		\tilde{\mathbf y})\right)^{2}}Z_{1}+\rho_{n,a,b}(\tilde{\mathbf x},
	\tilde{\mathbf y})Z_{2},Z_{2}\bigg),
\end{eqnarray*}
where 
\begin{eqnarray}\label{rho}
\rho_{n,a,b}(\tilde{\mathbf x},
\tilde{\mathbf y})&:=&\mathbb E [T_{\eta,a}(\tilde{\mathbf x}) T_{\eta,b}(\tilde{\mathbf y})]
\nonumber\\&=& \left[\prod_{l=1}^{p}w_{n,l}(\mathbf x_l,\mathbf y_l) +p(p-1)ab \prod_{l=1}^{p}\mathbb E f_{n,l}(\mathbf x_l) \mathbb E f_{n,l}(\mathbf y_l)
\nonumber\right.\\&&\hspace{0.2cm}	\left.+(ab-a-b)\sum_{l=1}^{p}w_{n,l}(\mathbf x_l,\mathbf y_l)\prod_{j\neq l}\mathbb Ef_{n,j}(\mathbf x_j)\mathbb Ef_{n,j}(y_j)\right]\nonumber \\&&\times\frac{1}{\sqrt{
		k_{n,a}(\tilde{\mathbf x})k_{n,b}(\tilde{\mathbf y})}},
\end{eqnarray}
where, for each $l=1,\ldots,p$,
\begin{eqnarray*}
	w_{n,l}(\mathbf x_l,\mathbf y_l):=h_{n}^{-2d_l} \mathbb E\Bigg[K_{l}\left(\frac{\mathbf x_l-\mathbb X^{l}}{h_n}\right)K_{l}\Bigg(\frac{\mathbf y_l-\mathbb X^{l}}{h_n}\Bigg)\Bigg].
\end{eqnarray*}
Set
\begin{eqnarray*}
	\bar{{\tau}}_{n,a,b}(\bar{C}, A_{n}(I_p))
	&=&\tilde{\bf{K}}\int_{\bar{C}}\int_{\bar{C}}\mathds{1}_{A_{n}(I_p)}(\tilde{\mathbf x},
	\tilde{\mathbf y})	 \mathrm{Cov}\big(|Z^{(a,b)}_{n,1}(\tilde{\mathbf x})|, |Z^{(a,b)}_{n,2}(\tilde{\mathbf y})|\big)
	\\&&\times h_{n}^{-d}	\prod_{l=1}^{p}\sqrt{f_{l}(\mathbf x_l)f_{l}(\mathbf y_l)} 		\,d\tilde{\mathbf x}\,d\tilde{\mathbf y},
\end{eqnarray*}
which by the change of variables $
\tilde{\mathbf y} = \tilde{\mathbf x}+
h_n\tilde{\mathbf t}$ equals
\begin{eqnarray*}
	\int_{\bar{C}}\int_{\tilde{\B}}
	g_{n,a,b}(\tilde{\mathbf x},\tilde{\mathbf t})
	\,d\tilde{\mathbf x}\,d\tilde{\mathbf t},
\end{eqnarray*}
where
\begin{eqnarray*} g_{n,a,b}(\tilde{\mathbf x},\tilde{\mathbf t})&:=&\tilde{\bf{K}}\mathds{1}_{\bar{C} }(\tilde{\mathbf x})\mathds{1}_{\bar{C}}(\tilde{\mathbf x}+
	h_n\tilde{\mathbf t}) \mathrm{Cov}\big(|Z^{(a,b)}_{n,1}(\tilde{\mathbf x})|, |Z^{(a,b)}_{n,2}( \tilde{\mathbf x}+
	h_n\tilde{\mathbf t} )|\big)
	\\&&\times	\prod_{l=1}^{p}\sqrt{f_{l}(\mathbf x_l)f_{1}(\mathbf x_l+h_n\mathbf t_l)}.
\end{eqnarray*}
We will show that, as $n\rightarrow \infty$,
\begin{eqnarray}\label{5eqn1}
\bar{{\tau}}_{n,a,b}(\bar{C}, A_{n}(I_p))\rightarrow \sigma^2\prod_{l=1}^{p}\int_{C_l}f_l(\mathbf x_l)\,d\mathbf x_l
\end{eqnarray}
and then, as $n\rightarrow \infty$, we have 
\begin{eqnarray}\label{5eqn2}
\bar{{\tau}}_{n,a,b}(\bar{C}, A_{n}(I_p))-\bar{{\sigma}}_{n,a,b}(\bar{C}, A_{n}(I_p))\rightarrow 0,
\end{eqnarray}
which will complete the proof of (\ref{t1var}). Now, consider (\ref{5eqn1}). Applying (\ref{lem1eqn1}) of Lemma \ref{lem1}, with $g(\cdot)=f_l(\cdot)$ and $K(\cdot)=K_l(\cdot+\mathbf t_l)K_l(\cdot+\mathbf s_l)$, $t_l,s_l\in\R^{d_l}$, we get, for each $(\mathbf t_l,\mathbf s_l)$, as $n\rightarrow \infty$, for almost every $\mathbf x_l\in\R^{d_l}$, hence for almost every $\mathbf x_l \in C_l$,
\begin{eqnarray}\label{varcon1}
&&h_{n}^{-d_l}\mathbb E\left[ K_l\left(\frac{\mathbf x_l-\mathbb X^{l}}{h_{n}}+\mathbf t_l\right)K_l\left(\frac{\mathbf x_l-\mathbb X^{l}}{h_{n}}+\mathbf s_l\right
)\right]\nonumber\\&& \hspace{3cm}\rightarrow f_l(\mathbf x_l)\int K_l(\mathbf u_l+\mathbf t_l)K_l(\mathbf u_l+\mathbf s_l)\,d\mathbf u_l.
\end{eqnarray}
Moreover, we get with $g(\cdot)=f_l(\cdot)$ and $K(\cdot)=K_l(\cdot+\mathbf t_l)$, as $n\rightarrow \infty$, for almost every $\mathbf x_l \in C_l$,
\begin{eqnarray}\label{varcon2}
&&h_{n}^{-d_l}\mathbb E\left[ K_l\left(\frac{\mathbf x_l-\mathbb X^{l}}{h_{n}}+\mathbf t_l\right)\right]\rightarrow f_l(\mathbf x_l).
\end{eqnarray}
Now, by (\ref{key1})-(\ref{eqnvareta}) in combination with (\ref{varcon1}) and (\ref{varcon2}), we get, for each $\tilde{\mathbf t}$, as $n\rightarrow \infty$, for almost every $\tilde{\mathbf x}\in \bar{C}$,
\begin{eqnarray}\label{eqn38}
h_n^{d}k_{n,a}( \tilde{\mathbf x}+
h_n\tilde{\mathbf t})\rightarrow \tilde{\bf{K}} \prod_{l=1}^{p}f_l(\mathbf x_l).
\end{eqnarray}
Thus, by (\ref{rho}), (\ref{varcon1}) and (\ref{eqn38}), we have, for each $\tilde{\mathbf t}$ and almost every $\tilde{\mathbf x}\in \bar{C}$, as $n\rightarrow \infty$,
\begin{eqnarray*}
	\rho_{n,a,b}(\tilde{\mathbf x},\tilde{\mathbf x} + h_n\tilde{\mathbf t})\rightarrow \prod_{l=1}^{p}\rho_l(\mathbf t_l)=\rho(\tilde{\mathbf t}),
\end{eqnarray*}
thus we obtain 
\begin{eqnarray*}
	\mathrm{Cov}\left(|Z^{(a,b)}_{n,1}(\tilde{\mathbf x})|,|Z^{(a,b)}_{n,2}(\tilde{\mathbf x} + h_n\tilde{\mathbf t})|\right)\rightarrow
	\mathrm{Cov}\left( \left|\sqrt{1-\rho^{2}(\tilde{\mathbf t}) }Z_{1}+ \rho(\tilde{\mathbf t}) Z_{2}\right|,|Z_{2}|\right).
\end{eqnarray*}
Combining this with Lemma \ref{lem4b} and the continuity of $f_l(\cdot)$ on $C_l$ , $l=1, \ldots, p$, we readily conclude that 
$g_{n,a,b}(\tilde{\mathbf x},\tilde{\mathbf t})$ converges in
measure on $\bar{C}\times\tilde{\B}$ to
\begin{eqnarray*}
	\tilde{\bf{K}}\mathds{1}_{\bar{C}}(\tilde{\mathbf x}) \mathrm{Cov}\Bigg( \left|\sqrt{1-\rho^{2}(\tilde{\mathbf t})}Z_{1}+ \rho(\tilde{\mathbf t}) Z_{2}\right|,|Z_{2}|\Bigg) \prod_{l=1}^{p}f_l(\mathbf x_l).
\end{eqnarray*}
Since, for each $l=1,\ldots,p$, $f_l(\cdot)$ is bounded on $C_l$, the function $g_{n,a,b}(\tilde{\mathbf x},\tilde{\mathbf t})$ is for all $n\geq1$ uniformly bounded on
$\bar{C}\times \tilde{\B}$. Thus, we get by the Lebesgue bounded
convergence theorem, as $n\rightarrow \infty$, $\bar{\tau}_{a,b,n}(\bar{C}, A_{n}(I_p))$ converges to
\begin{eqnarray*}
	\sigma^{2}\prod_{l=1}^{p}\int_{C_l}f_l(\mathbf x_l)d\mathbf x_l.
\end{eqnarray*}
This completes the proof of (\ref{5eqn1}). Now, we consider the relation (\ref{5eqn2}). Set
\begin{eqnarray*} G_{n}(\tilde{\mathbf x},\tilde{\mathbf t})&:=&\tilde{\bf{K}}\mathds{1}_{\bar{C}}(\tilde{\mathbf x})\mathds{1}_{\bar{C}}(\tilde{\mathbf x}+h_n\tilde{\mathbf t})
	\prod_{l=1}^{p} \sqrt{f_l(\mathbf x_l)f_l(\mathbf x_l + h_n\mathbf t_l)}.
\end{eqnarray*}
Notice that, since, for each $l=1,\ldots,p$, $f_l(\cdot)$ is bounded on $C_l$,
\begin{eqnarray}\label{5eqn3}
\int_{\bar{C}}\int_{\tilde{\B}} G_{n}(\tilde{\mathbf x},\tilde{\mathbf t}) \,d\tilde{\mathbf x}\, d\tilde{\mathbf t}
\leq \beta_1,
\end{eqnarray}
where $\beta_1$ is a positive constant. 
We see that
\begin{eqnarray*}\label{5eq4}
	\lefteqn{|\bar{\tau}_{n,a,b}(\bar{C},A_{n}(I_p))-\bar{\sigma}_{n,a,b}(\bar{C},A_{n}(I_p))|}\nonumber
	\\& \leq&
	\int_{\bar{C}}\int_{\tilde{ \B}}\big|\mathbb E|Z^{(a,b)}_{n,1}(\tilde{\mathbf x})| \mathbb E|Z^{(a,b)}_{n,2}(\tilde{\mathbf x}+h_n\tilde{\mathbf t})|\nonumber
	-\mathbb E|T_{\eta,a}(\tilde{\mathbf x})|\mathbb E|T_{\eta,b}(
	\tilde{\mathbf x}+h_n\tilde{\mathbf t})|\big| G_{n}(\tilde{\mathbf x},\tilde{\mathbf t})\,d\tilde{\mathbf x}\, d\tilde{\mathbf t}\nonumber
	\\&&+\int_{\bar{C}}\int_{\tilde{ \B}}\big|\mathbb E|Z^{(a,b)}_{n,1}(\tilde{\mathbf x})Z^{(a,b)}_{n,2}(\tilde{\mathbf x}+h_n\tilde{\mathbf t})|-\mathbb E|T_{\eta,a}(\tilde{\mathbf x})T_{\eta,b}(\tilde{\mathbf x}+h_n\tilde{\mathbf t})|\big| G_{n}(\tilde{\mathbf x},\tilde{\mathbf t}) \,d\tilde{\mathbf x}\, d\tilde{\mathbf t}\nonumber
	\\& :=&\zeta_{n,a,b}(1)+\zeta_{n,a,b}(2).
\end{eqnarray*}
First, using (\ref{5eqn3}) and again the statement (\ref{equa6.18Fact6.1}) with (\ref{l4eqn4}), we obtain
\begin{eqnarray*}
	\zeta_{n,a,b}(1)=O\left(\frac{1}{\sqrt{nh_n^{d}}}\right).
\end{eqnarray*}
Choose any $0<\varepsilon<1$ and set
\begin{eqnarray*}
	E_n(\varepsilon)=\big\{(\tilde{\mathbf x},\tilde{\mathbf t}): 1-\big(\rho_{n,a,b}(\tilde{\mathbf x},\tilde{\mathbf x}+h_n\tilde{\mathbf t})\big)^{2}\geq\varepsilon\big\}.
\end{eqnarray*}
Now, we infer that
\begin{eqnarray*}
	\zeta_{n,a,b}(2)&\leq&\int_{\bar{C}}\int_{\bar{\B}}\big|1-\mathbb E|Z^{(a,b)}_{n,1}(\tilde{\mathbf x})Z^{(a,b)}_{n,2}(\tilde{\mathbf x}+ h_n\tilde{\mathbf t})|\big|
	\mathds{1}_{E_n^{c}(\varepsilon)}(\tilde{\mathbf x},\tilde{\mathbf t})
	G_{n}(\tilde{\mathbf x},\tilde{\mathbf t}) \,d\tilde{\mathbf x}\, d\tilde{\mathbf t}\\&&+\int_{\bar{C}}\int_{\bar{\B}}\big|1-\mathbb E|T_{\eta,a}(\tilde{\mathbf x})T_{\eta,b}(\tilde{\mathbf x}+h_n\tilde{\mathbf t})|\big|
	\mathds{1}_{E_n^{c}(\varepsilon)}(\tilde{\mathbf x},\tilde{\mathbf t})
	G_{n}(\tilde{\mathbf x},\tilde{\mathbf t}) \,d\tilde{\mathbf x}\, d\tilde{\mathbf t}\\&& +\int_{\bar{C}}\int_{\bar{\B}}\big|\mathbb E|Z^{(a,b)}_{n,1}(\tilde{\mathbf x})Z^{(a,b)}_{n,2}(\tilde{\mathbf x}+ h_n\tilde{\mathbf t})|-\mathbb E|T_{\eta,a}(\tilde{\mathbf x})T_{\eta,b}(\tilde{\mathbf x}+h_n\tilde{\mathbf t})|\big|\\&&\times \mathds{1}_{E_n(\varepsilon)}(\tilde{\mathbf x},\tilde{\mathbf t})
	G_{n}(\tilde{\mathbf x},\tilde{\mathbf t}) \,d\tilde{\mathbf x}\, d\tilde{\mathbf t}\\&
	:=&\bar\zeta^{(a,b)}_{n,1}(2,\varepsilon)+\bar\zeta^{(a,b)}_{n,2}(2,\varepsilon)+\zeta^{(a,b)}_{n}(2,\varepsilon)\\&
	:=&\bar\zeta^{(a,b)}_{n}(2,\varepsilon)+\zeta^{(a,b)}_{n}(2,\varepsilon).
\end{eqnarray*}
To bound $\bar\zeta_{n}^{a,b}(2,\varepsilon)$, we use the elementary fact
that if $X$ and $Y$ are mean-zero and variance-one random variables
with $\rho=\mathbb E(XY)$, then $1 -\mathbb E|XY|\leq1-|\rho|\leq 1-\rho^{2}$, in combination with (\ref{5eqn3}), to get 
$$\bar\zeta^{(a,b)}_{n}(2,\varepsilon)\leq 2\varepsilon \beta_1. $$
Next, using the inequality (\ref{5eqn3}) and the statement (\ref{equa6.19Fact6.1}) along with the statement (\ref{l4eqn4}), we obtain
$$
\zeta^{(a,b)}_{n}(2,\varepsilon)=O\left(\frac{1}{\sqrt{nh_{n}^{d}}}\right).
$$
Thus, for any $0< \varepsilon<1$, we have
$$ \limsup_{n\rightarrow\infty}
\left|\bar{\tau}_{n,a,b}(\bar{C},A_{n}(I_p))-\bar{\sigma}_{n,a,b}(\bar{C},A_{n}(I_p))\right|\leq 2\varepsilon \beta_1,$$ which yields the statement (\ref{5eqn2}). This finishes the proof of (\ref{t1var}). Now we turn to (\ref{t2var}). 
Let us show that, for each $I\subsetneq I_p$, as $n\rightarrow \infty$,
\begin{eqnarray}\label{5eqn4}
\sigma_{n,a,b}(\bar{C},A_{n}(I))\rightarrow 0.
\end{eqnarray} 
Set, for each $I\subsetneq I_p$,
\begin{eqnarray*}
	&&\tau_{n,a,b}(\bar{C},A_{n}(I)):=\int_{\bar{C}}\int_{\bar{C}}
	\mathrm{Cov}\big(|Z^{(a,b)}_{n,1}(\tilde{\mathbf x})|, |Z^{(a,b)}_{n,2}(\tilde{\mathbf y})|\big)	\mathds{1}_{A_{n}(I)}(\tilde{\mathbf x},\tilde{\mathbf y})	 
	\\&&\hspace{4cm}\times 	\sqrt{k_{n,a}(\tilde{\mathbf x})k_{n,b}(\tilde{\mathbf y})}
	\,d\tilde{\mathbf x}\,\tilde{\mathbf y}.
\end{eqnarray*}
We will show that, for each $I\subsetneq I_p$, as $n\rightarrow \infty$,
\begin{eqnarray}\label{5eqn12}
\tau_{n,a,b}(\bar{C}, A_{n}(I))\rightarrow 0
\end{eqnarray}
and then, as $n\rightarrow \infty$,
\begin{eqnarray}\label{5eqn22}
\tau_{n,a,b}(\bar{C}, A_{n}(I))-\sigma_{n,a,b}(\bar{C},A_{n}(I))\rightarrow 0,
\end{eqnarray}
which completes the proof of (\ref{5eqn4}).
Notice that, for each $I\subsetneq I_p$ and each $(\tilde{\mathbf x},\tilde{\mathbf y})\in A_n(I)$, 
\begin{eqnarray*}
	\rho_{n,a,b}(\tilde{\mathbf x},\tilde{\mathbf y})&=&\left[p(p-1)ab \prod_{l=1}^{p}\mathbb E f_{n,l}(\mathbf x_l)\mathbb E f_{n,l}(\mathbf y_l)\right.\nonumber\\&&+\left.
	(ab-a-b)\sum_{l\in I}w_{n,l}(\mathbf x_l,\mathbf y_l)\times\prod_{j\neq l}\mathbb Ef_{n,j}(\mathbf x_j)\mathbb Ef_{n,j}(y_j)\right]\bigg{/}\sqrt{
		k_{n,a}(\tilde{\mathbf x})k_{n,b}(\tilde{\mathbf y})}.
\end{eqnarray*}
Therefore, by (\ref{majora1}), (\ref{l4eqn2}) and the fact that, for each $l=1,\ldots,p$, $K_l$ is bounded, we get, for each $I\subsetneq I_p$, for all large enough $n$ uniformly in
$(\tilde{\mathbf x}, \tilde{\mathbf y})\in A_{n}(I) \cap \bar{C}\times \bar{C}$ and for some constant $\beta_2>0$,
\begin{eqnarray}\label{5eqn5}
|\rho_{n,a,b}(\tilde{\mathbf x}, \tilde{\mathbf y})|\leq \beta_2 \prod_{l\in I_p\setminus I}h_n^{d_l}.
\end{eqnarray}
Now, by using \cite{nabeya} formulas, we have 
\begin{eqnarray*}
	\mathrm{Cov}\left(|Z^{(a,b)}_{n,1}(\tilde{\mathbf x})|,|Z^{(a,b)}_{n,2}(\tilde{\mathbf y})|\right)
	&=&\mathrm{Cov}\big( \big|\sqrt{1-(\rho_{n,a,b}(\tilde{\mathbf x},\tilde{\mathbf y}))^{2}}Z_{1}+ \rho_{n,a,b}(\tilde{\mathbf x},\tilde{\mathbf y})Z_{2}\big|,|Z_{2}|\big)\\
	&:= &\varphi(\rho_{n,a,b}(\tilde{\mathbf x},\tilde{\mathbf y})),
\end{eqnarray*}
where
\begin{eqnarray*}
	\varphi(\rho)=\frac{2}{\pi}\left(\rho \arcsin \rho+\sqrt{1-\rho^{2}}-1\right), \quad \rho\in \lbrack-1,1\rbrack,
\end{eqnarray*}
which in turn, for each $\rho\in\lbrack-1,1\rbrack$, is less than or equal to $\frac{(\pi-2)\rho^2}{\pi}$.
Thus, for each $I\subsetneq I_p$, for all large enough $n$ uniformly in
$(\tilde{\mathbf x}, \tilde{\mathbf y})\in A_{n}(I)\cap \bar{C}\times \bar{C}$ and for some constant $\beta_3>0$, by (\ref{5eqn5}) we get 
\begin{eqnarray*}\label{5eqn6b}
	\mathrm{Cov}\big(|Z^{(a,b)}_{n,1}(\tilde{\mathbf x}, \tilde{\mathbf y})|,|Z^{(a,b)}_{n,2}(\tilde{\mathbf x}, \tilde{\mathbf y})|\big)\leq \beta_3\, \prod_{l\in I_p \setminus I}h_n^{2d_l}.
\end{eqnarray*} 
By change of variables $\mathbf y_l=\mathbf x_l+h_n\mathbf t_l$ for each $l\in I$, for some constant $\beta_4>0$, we obtain
\begin{eqnarray*}
	\int_{\bar{C}}\int_{\bar{C}}\mathds{1}_{A_{n}(I)} (\tilde{\mathbf x},\tilde{\mathbf y})\,d\tilde{\mathbf x}\, d\tilde{\mathbf y}
	\leq \beta_4 \prod_{l\in I}h_{n}^{d_l}.
\end{eqnarray*}
Hence, by (\ref{l4eqn2}) we have 
\begin{eqnarray}\label{5eqn32}
\int_{\bar{C}}\int_{\bar{C}}\mathds{1}_{A_{n}(I) }(\tilde{\mathbf x},\tilde{\mathbf y})\sqrt{k_{n,a}(\tilde{\mathbf x})k_{n,b}(\tilde{\mathbf y})}\,d\tilde{\mathbf x}\, d\tilde{\mathbf y}
\leq \beta_5 \prod_{l\in I_p\setminus I}h_{n}^{-d_l}
\end{eqnarray}
for some $\beta_5>0$. This completes the proof of (\ref{5eqn12}). 
We next evaluate the difference in (\ref{5eqn22}).
We can see that
\begin{eqnarray*}\label{5eq4LLLL}
	\lefteqn{|\tau_{n,a,b}(\bar{C}, A_{n}(I))-\sigma_{n,a,b}(\bar{C},A_{n}(I))|}\nonumber
	\\ &\leq&
	\int_{\bar{C}}\int_{\bar{C}}\big|\mathbb E|Z^{(a,b)}_{n,1}(\tilde{\mathbf x})| \mathbb E|Z^{(a,b)}_{n,2}(\tilde{\mathbf y})|-\mathbb E|T_{\eta,a}(\tilde{\mathbf x})|\mathbb E|T_{\eta,b}(
	\tilde{\mathbf y})|\big|\nonumber
	\\&&\times \mathds{1}_{A_{n}(I) }(\tilde{\mathbf x},\tilde{\mathbf y})\sqrt{k_{n,a}(\tilde{\mathbf x})k_{n,b}(\tilde{\mathbf y})}\,d\tilde{\mathbf x}\, d\tilde{\mathbf y}\nonumber
	\\&&+\int_{\bar{C}}\int_{\bar{C}}\big|\mathbb E|Z^{(a,b)}_{n,1}(\tilde{\mathbf x})Z^{(a,b)}_{n,2}(
	\tilde{\mathbf y})|-\mathbb E|T_{\eta,a}(\tilde{\mathbf x})T_{\eta,b}(
	\tilde{\mathbf y})|\big|\nonumber\\&&\times\mathds{1}_{A_{n}(I)}(\tilde{\mathbf x},\tilde{\mathbf y})\sqrt{k_{n,a}(\tilde{\mathbf x})k_{n,b}(\tilde{\mathbf y})}\,d\tilde{\mathbf x}\, d\tilde{\mathbf y}\nonumber
	\\& :=&\xi_{n,a,b}(1)+\xi_{n,a,b}(2).
\end{eqnarray*}
To bound $\xi_{n,a,b}(1)$, we use (\ref{5eqn32}) and (\ref{equa6.18Fact6.1}) in combination with (\ref{l4eqn4}), to get 
\begin{eqnarray*}
	\xi_{n,a,b}(1)=O\left(\frac{1}{\sqrt{nh_n^{3d}}}\right).
\end{eqnarray*}
Next, we use (\ref{5eqn32}) and (\ref{equa6.19Fact6.1}) with (\ref{l4eqn4}) and (\ref{5eqn5}) to get
\begin{eqnarray*}\xi_{n,a,b}(2) =O\left(\frac{1}{\sqrt{nh_n^{3d}}}\right),
\end{eqnarray*}
which gives (\ref{5eqn22}).
This completes the proof of Lemma \ref{lem5}. \hfill $\blacksquare$
\\
~

\noindent In the proof of the next lemmas, we shall apply Lemma \ref{lem5} with $a =b=0$. Note that in this situation, whenever $ (\tilde{\mathbf x},\tilde{\mathbf y})\in A_{n}(I)$ with $I\subsetneq I_p$, the random variables $| \Gamma_{\eta,0}(\tilde{\mathbf x})|$ and $| \Gamma_{\eta,0}(\tilde{\mathbf y})|$ are independent. This follows from the fact that being functions of independent increments of a Poisson process. Therefore, $\sigma_{n,0,0}(\bar{C},A_{n}(I))=0$ for all $n$ and each $I\subsetneq I_p$, which implies that 
\begin{eqnarray}\label{Varf}
\lim_{n\rightarrow\infty} n\mathrm{Var}\big(U_{\eta,0}(\bar{C}) \big)=\sigma^2\prod_{l=1}^{p}\int_{C_l}f_l(\mathbf x_l)\,d\mathbf x_l
\end{eqnarray}
can be obtained under the bandwidth condition $nh_n^{d}\rightarrow\infty$ instead of $nh_n^{3d}\rightarrow\infty$ as $n\rightarrow\infty$. 

\begin{lemma}\label{lem51}
	Suppose that, for each $l=1,\ldots,p$,  $C_l$ satisfies (\ref{lem1eqn2})-(\ref{lem1eqn3})-(\ref{lem1eqn4}) of Lemma \ref{lem1} with $g(\cdot)=f_l(\cdot)$ and $\HH(\cdot)=\HH_l(\cdot)$. If
	$h_n\rightarrow 0$ and $nh_n^{3d}\rightarrow\infty$, then we have 
	\begin{eqnarray*}
		\sqrt{n}\big\lbrack U_{n,1}(\bar{C})-U_{n,0}(\bar{C})-\mathbb E\big(U_{n,1}(\bar{C})-U_{n,0}(\bar{C})\big)\big\rbrack\stackrel{\PP}{\rightarrow}0.
	\end{eqnarray*}
\end{lemma}
\noindent{\bf{Proof of Lemma \ref{lem51}.}} 
For each $l=1,\ldots,p$, we can find a measurable partition $C_{l,1},\ldots, C_{l,k_l}$ of $C_{l}$ so that $$0<\int_{ C_{l,j_l}^{h/2}}f_l(\mathbf u_l)\,d\mathbf u_l<1/(2p),$$ for each $j_l =1,\ldots, k_l$ and all $h_n > 0$ small enough, where $C_{l,j_l}^{h/2}$ is the $h/2$-neighborhood of $C_{l,j_l}$. 
We will apply Lemma $2.1$ of \cite{mason2001} to the semigroup $D$ generated by the point masses, $$D=\left\{ 0, \sum_{i=1}^{n}\delta_{\tilde{\mathbf x}_i}: n\in\N, \tilde{\mathbf x}_i\in \R^{d_1}\times \cdots\times\R^{d_p}\right\},$$ with the $\sigma$-algebra $\D$ generated by the functions 
\begin{eqnarray*}
	f_{n,B_{j_1,\ldots,j_p}}(\tilde{\mathbf x}_1,\ldots, \tilde{\mathbf x}_n)=\sum_{i=1}^{n}
	\mathds{1}_{\left\lbrace \delta_{\tilde{\mathbf x}_i}\in B_{j_1,\ldots,j_p}\right\rbrace } \delta_{\tilde{\mathbf x}_i},  n\in\N, 
\end{eqnarray*} 
where, for each $(j_1,\ldots,j_p)\in \{1,\ldots,k_1\}\times\cdots\times\{1,\ldots,k_p\} =:J_{k_1,\ldots,k_p}$,
\begin{eqnarray*}
	B_{j_1,\ldots,j_p}=\left\{\delta_{\tilde{\mathbf x}}: \tilde{\mathbf x}\in  \R^{d_1}\times\cdots\times\R^{d_{l-1}}\times C_{l,j_l}^{h/2}\times\R^{d_{l+1}}\times\cdots\times \R^{d_p}, \mbox{ for some } l=1,\ldots,p \right\}.
\end{eqnarray*} 
It is easy to see that, for
any measurable function $G: \R^{d_1}\times\cdots \times\R^{d_p}\mapsto\R$, the map $\mu\mapsto \int G\, d\mu$ is $\D$-measurable [just note that $$f^{-1}_{n, B_{j_1,\ldots,j_p}} \left\{ \mu\in D, \int G\, d\mu\leq t \right\}=\left\{(\tilde{\mathbf x}_1,\ldots, \tilde{\mathbf x}_n): \sum_{i=1}^{n}
\mathds{1}_{\left\lbrace \delta_{\tilde{\mathbf x}_i}\in B_{j_1,\ldots,j_p}\right\rbrace } G(\tilde{\mathbf x}_i)\leq t\right\}$$ is measurable set of $(\R^{d_1}\times\cdots \times\R^{d_p})^{n}$ ].
The function $H=H_{n,j_1,\ldots,j_p}$ has the form 
\begin{eqnarray*}
	\lefteqn{H_{n,j_1,\ldots,j_p}\left(\sum_{i=1}^{n}\mathds{1}_{\left\lbrace \delta_{
				\tilde{\mathbf x}_i}\in B_{j_1,\ldots,j_p}\right\rbrace }\delta_{ \tilde{\mathbf x}_i}\right)}\\&=& \left\{\bar{H}_{n,j_1,\ldots,j_p,{\bf{1}}}\Big(\sum_{i=1}^{n}\mathds{1}_{\left\lbrace \delta_{
			\tilde{\mathbf x}_i}\in B_{j_1,\ldots,j_p}\right\rbrace }\delta_{ \tilde{\mathbf x}_i}\Big)-\bar{H}_{n,j_1,\ldots,j_p,{\bf{0}}}\left(\sum_{i=1}^{n}\mathds{1}_{\left\lbrace \delta_{
			\tilde{\mathbf x}_i}\in B_{j_1,\ldots,j_p}\right\rbrace } \delta_{ \tilde{\mathbf x}_i}\right)\right.\\&&
	-nh_{n}^{d}
	\mathbb	E\Big(U_{\eta,1}(\bar{C}_{j_1,\ldots,j_p})
	-
	U_{\eta,0}(\bar{C}_{j_1,\ldots,j_p}
	)\Big)\Bigg\}^{2},
\end{eqnarray*} 
where $\bar{C}_{j_1,\ldots,j_p}=C_{1,j_1}\times\cdots\times C_{p,j_p}$ and,
for each $a\in\{0,1\}$, 
we can therefore write 
\begin{eqnarray*}
	\lefteqn{\bar{H}_{n,j_1,\ldots,j_p,{\bf{a}}}\left(\sum_{i=1}^{n}\mathds{1}_{\left\lbrace \delta_{
				\tilde{\mathbf x}_i}\in B_{j_1,\ldots,j_p}\right\rbrace } \delta_{ \tilde{\mathbf x}_i}\right)
	}\\&=& \int_{\bar{C}_{j_1,\ldots,j_p}}
	\left| \sum_{i=1} ^{n}\left\{\prod_{l=1}^{p}K_{l}\left(\frac{\mathbf x_l-\mathbf x_{i,l}}{h_n}\right) 
	-a\sum_{l=1}^{p} K_{l}\left(\frac{\mathbf x_l-\mathbf x_{i,l}}{h_n}\right)\prod_{j\neq i}\mathbb EK_j\left(\frac{\mathbf x_j-\mathbb X^{j}}{h_n}\right)
	\right\} \right. \\	&&\times \left.\mathds{1}_{\left\lbrace \delta_{
			\tilde{\mathbf x}_i}\in B_{j_1,\ldots,j_p}\right\rbrace } +
	n(pa-1)\prod_{l=1}^{p}\mathbb EK_{l}\left(\frac{\mathbf x_l-\mathbb X^{l}}{h_n}\right)	\right| \,d\tilde{\mathbf x},
\end{eqnarray*}
were we recall $\tilde{\mathbf x}_i=(x_{i,1}^{\top},\ldots,x_{i,p}^{\top})^{\top}\in \R^{d_1}\times\cdots \times\R^{d_p}$. Now, $H_{n,j_1,\ldots,j_p}$ when considered as a function on $(\R^{d_1}\times\cdots \times\R^{d_p})^{n}$
can be shown to be Borel measurable, which
in this setup is equivalent to be $\D$-measurable.
We now get from Lemma $2.1$ of \cite{mason2001}, for each $(j_1,\ldots,j_p)\in J_{k_1,\ldots,k_p}$,
\begin{eqnarray*}
	\mathbb EH_{n,j_1,\ldots,j_p}\Bigg(\sum_{i=1}^{n}\mathds{1}_{\left\lbrace \delta_{
			\tilde{\mathbf X}_i}\in B_{j_1,\ldots,j_p}\right\rbrace }\delta_{\tilde{\mathbf X}_i}\Bigg) \leq 2\mathbb EH_{n,j_1,\ldots,j_p}\Bigg(\sum_{i=1}^{\eta}\mathds{1}_{\left\lbrace \delta_{
			\tilde{\mathbf X}_i}\in B_{j_1,\ldots,j_p}\right\rbrace }\delta_{\tilde{\mathbf x}_i}\Bigg),
\end{eqnarray*}
with $\eta$ is a Poisson random variable with mean $n$, independent
of $\tilde{\mathbf X}, \tilde{\mathbf X}_1,\tilde{\mathbf X}_2\ldots$ Therefore, by assumption 
\textnormal{({\bf A$_1$})}, for each $(j_1,\ldots,j_p)\in J_{k_1,\ldots,k_p}$, we have
\begin{eqnarray*}
	&&n\mathbb E\bigg\lbrack\Big\{ U_{n,1}(\bar{C}_{j_1,\ldots,j_p})-U_{n,0}(\bar{C}_{j_1,\ldots,j_p})-
	\mathbb E\big(U_{\eta,1}(\bar{C}_{j_1,\ldots,j_p})-
	U_{\eta,0}(\bar{C}_{j_1,\ldots,j_p})\big)\Big\}^{2}\bigg
	\rbrack\\&& \hspace{2cm}\leq 2 n\mathrm{Var}\Big( U_{\eta,1}(\bar{C}_{j_1,\ldots,j_p})-U_{\eta,0}(\bar{C}_{j_1,\ldots,j_p})\Big).
\end{eqnarray*}
Using Lemma \ref{lem5} with $\bar{C}=\bar{C}_{j_1,\ldots,j_p}$, as $ n\rightarrow\infty$, we get
\begin{eqnarray*}
	n\mathrm{Var}\Big( U_{\eta,1}(\bar{C}_{j_1,\ldots,j_p})-U_{\eta,0}(\bar{C}_{j_1,\ldots,j_p})\Big)\rightarrow
	0
\end{eqnarray*}
and thus, for each $(j_1,\ldots,j_p)\in J_{k_1,\ldots,k_p}$,
\begin{eqnarray*}
	\sqrt{n}\big\lbrack U_{n,1}(\bar{C}_{j_1,\ldots,j_p})-U_{n,0}(\bar{C}_{j_1,\ldots,j_p})-\mathbb E\big(U_{\eta,1}
	(\bar{C}_{j_1,\ldots,j_p})-U_{\eta,0}(\bar{C}_{j_1,\ldots,j_p})\big)\big\rbrack\stackrel{\PP}{\rightarrow}0
\end{eqnarray*}
which, when coupled with (\ref{keyeqn1}), (\ref{keyeqn2}) and, for each $l=1,\ldots,p$, $C_{l,1},\ldots, C_{l,k_l}$ being a partition of $C_l$, completes the proof of Lemma \ref{lem51}. \hfill $\blacksquare$ 
\\
~

\noindent Suppose that, for each $l=1,\ldots,p$, $C_l$ satisfies Lemma \ref{lem1} with $g(\cdot)=f_l(\cdot)$ and $\HH(\cdot)=\HH_l(\cdot)$. For all $0<\varepsilon<1$ and each $l=1,\ldots,p$, let $M_l, \nu_l$ and $\mu_l$ be the numbers  
from (\ref{lem1eqn2}), namely 
\begin{eqnarray*}
	&&C_l\subset\lbrack -M_l+\nu_l, M_l-\nu_l\rbrack^{d_l}, \qquad\qquad0< \mu_l=\int_{ \R^{d_l}\setminus \lbrack -M_l, M_l\rbrack^{d_l}
	} f_l(\mathbf u_l)\,d\mathbf u_l <1\end{eqnarray*}
and 	\begin{eqnarray*} 
	\int_{C_l} f_l(\mathbf u_l)\,d\mathbf u_l>1-\varepsilon. 
\end{eqnarray*}
Assume that $n$ is large enough such that $h_n\leq \min(\nu_1,\ldots,\nu_p,M_1/2,\ldots,M_p/2)$. Define, for $l=1,\ldots,p$, $m_{n,l}=\lfloor M_l/h_n\rfloor $, where $\lfloor u\rfloor $ denotes the integer part, i.e., $\lfloor u\rfloor \leq u < \lfloor u \rfloor + 1$, and $h_{n,l}=M_l/m_{n,l}$. For $l=1,\ldots,p$, we have $M_l/(2h_n)\leq m_{n,l}\leq M_l/h_n$ and thus $h_n\leq h_{n,l}\leq 2h_n.$
Next, consider the regular grid given by
\begin{eqnarray*}
	A_\textbf{i}&=&[x_{1,i_1},x_{1,i_1+1}]\times\cdots \times [x_{1,i_{d_1}},x_{1,i_{d_1}+1}] \\&&\times\cdots\times
	[x_{p,i_{d-d_{p}+1}},x_{p,i_{d-d_{p}+1}+1}]\times \cdots \times [x_{p,i_{d}},x_{p,i_{d}+1}],
\end{eqnarray*}
where $\textbf{i}=(i_1,\ldots,i_{d})\in \Z^{d}$ and, for all $i\in\Z$, $x_{l,i}=ih_{n,l}$. Define
\begin{eqnarray*}
	R_\textbf{i}=A_\textbf{i}\cap\M,
\end{eqnarray*}
where $\M=[-M_1,M_1]^{d_1}\times \cdots \times [-M_p,M_p]^{d_p}$. With $l_n=\{\textbf{i}:
R_\textbf{i}\neq \emptyset \}$,  we see that $\{R_\textbf{i}: \textbf{i}\in l_n \subset \Z^d\}$ constitutes a partition of $\M$ such that, for each $\textbf{i}\in l_n$,
\begin{eqnarray}\label{lem6eq2b}
\lambda(R_\textbf{i})\leq (2h_n)^d\quad \mbox{and} \quad |l_n|\leq 2^d h_n^{-d}\prod_{l=1}^{p}M_l^{d_l}.
\end{eqnarray}
Set, for any $\textbf{i}\in l_n$,
\begin{eqnarray*}
	\X_{\textbf{i},n}:=\frac{\displaystyle \int_{R_\textbf{i}\cap \bar{C}} \Delta_n(\tilde{\mathbf x})\, d\tilde{\mathbf x}}{\sqrt{n\mathrm{Var}\big(U_{\eta,0}(\bar{C})\big)}},
\end{eqnarray*}
where
\begin{eqnarray*}
	\Delta_n(\tilde{\mathbf x}):=\sqrt{n}\left[| \Gamma_{\eta,0}(\tilde{\mathbf x})| - \mathbb E|\Gamma_{\eta,0}(\tilde{\mathbf x})|\right]. 
\end{eqnarray*}
\begin{lemma}\label{lem8} Whenever
	$h_n\rightarrow 0$, $nh_n^{d}\rightarrow\infty$ 
	and, for each $l=1,\ldots,p$,  $C_l$ satisfies (\ref{lem1eqn2})-(\ref{lem1eqn3})-(\ref{lem1eqn4}) of Lemma \ref{lem1} with $g(\cdot)=f_l(\cdot)$ and $\HH(\cdot)=\HH_l(\cdot)$, there exists a constant $\gamma_0$ such that, uniformly in $\textbf{i}$ and for all $n$ sufficiently large,
	
	\begin{eqnarray}\label{lem8eq1}
	\mathbb E|\X_{\textbf{i},n}|^3\leq \gamma_0 h_{n}^{3d/2}.
	\end{eqnarray}
\end{lemma}
\noindent{\bf{Proof of Lemma \ref{lem8}.}} Observe that
\begin{eqnarray*}
	\Big(n\mathrm{Var}\big(U_{\eta,0}(\bar{C})\big)\Big)^{3/2}\mathbb E|\X_{\textbf{i},n}|^3
	\leq \int_{I_{\textbf{i},n}} \mathbb E \big| \triangle_n(\tilde{\mathbf x})\triangle_n(\tilde{\mathbf x}_1)\triangle_n(\tilde{\mathbf x}_2)\big|\, d\tilde{\mathbf x}\, d\tilde{\mathbf x}_1 \, d\tilde{\mathbf x}_2,
\end{eqnarray*}
where $I_{\textbf{i},n}= (R_\textbf{i}\cap\bar{C})^{3}$. It is obvious to see that
\begin{eqnarray*}
	\mathbb E \big| \triangle_n(\tilde{\mathbf x})\triangle_n(\tilde{\mathbf x}_1)\triangle_n(\tilde{\mathbf x}_2)\big|\leq\mathbb E \left[\big|
	\triangle_n(\tilde{\mathbf x})\big|+\big|\triangle_n(\tilde{\mathbf x}_1)\big|+\big|\triangle_n(\tilde{\mathbf x}_2)\big|\right]^3,
\end{eqnarray*}
which by $c_r$-inequality is, for some constant $\gamma_1 > 0$, less than or equal to
\begin{eqnarray*} \gamma_1n^{3/2} \left[\mathbb E\big|f_\eta(\tilde{\mathbf x})-\mathbb E f_n(\tilde{\mathbf x})\big|^3+\mathbb E\big|f_\eta(\tilde{\mathbf x}_1)-\mathbb E f_n(\tilde{\mathbf x}_1)\big|^3
	+\mathbb E\big|f_\eta(\tilde{\mathbf x}_2)-\mathbb E f_n(\tilde{\mathbf x}_2)\big|^3\right].
\end{eqnarray*}
By Lemma $2.3$ of \cite{mason2001}, we have, for any $\tilde{\mathbf x}\in\bar{C}$, 
\begin{eqnarray*}
	\lefteqn{n^{3/2}\mathbb E\big|f_\eta(\tilde{\mathbf x})-\mathbb E f_n(\tilde{\mathbf x})\big|^3}\\
	&\leq&\Bigg(
	\frac{45}{\log3}\Bigg)^3\max\Bigg\{\Bigg\lbrack
	\frac{1}{h_n^{2d}}\prod_{l=1}^{p}\mathbb EK_{l}^{2}\Bigg(\frac{\mathbf x_l-\mathbb X^{l}}{h_{n}}\Bigg)\Bigg\rbrack^{3/2}, 
	\frac{1}{\sqrt{n}h_n^{3d}}\prod_{l=1}^{p}\mathbb EK_{l}^{3}\Bigg(\frac{\mathbf x_l-\mathbb X^{l}}{h_{n}}\Bigg)
	\Bigg\},
\end{eqnarray*}
which by (\ref{majora1}) is, for $n$ large enough and
uniformly in $\tilde{\mathbf x}\in\bar{C}$, is bounded by 
\begin{eqnarray*}
	\gamma_2\Bigg\lbrack\frac{1}{h_n^{3d/2}}+\frac{1}{\sqrt{n}h_n^{2d}}\Bigg\rbrack,
\end{eqnarray*}
where $\gamma_2$ is some positive constant. Thus, uniformly in $\tilde{\mathbf x}\in\bar{C}$ and for $n$ large enough, we have
\begin{eqnarray*}\mathbb E \big| \triangle_n(\tilde{\mathbf x})\triangle_n(\tilde{\mathbf x}_1)\triangle_n(\tilde{\mathbf x}_2)\big|\leq 3
	\gamma_1\gamma_2\Bigg\lbrack\frac{1}{h_n^{3d/2}}+\frac{1}{\sqrt{n}h_n^{2d}}\Bigg\rbrack.
\end{eqnarray*}
This implies that, for some constant $\gamma_3>0$, uniformly in $\textbf{i}$ and for all large enough $n$,
\begin{eqnarray*}
	\int_{I_{\textbf{i},n}}\mathbb E \big| \triangle_n(\tilde{\mathbf x})\triangle_n(\tilde{\mathbf x}_1)\triangle_n(\tilde{\mathbf x}_2)\big|\, d\tilde{\mathbf x}\, d\tilde{\mathbf x}_1 \, d\tilde{\mathbf x}_2
	&\leq& 3(2^{3d})\gamma_1\gamma_2
	\Bigg\lbrack
	h_{n}^{3d/2}
	+\frac{h_{n}^{d}}{\sqrt{n}}\Bigg\rbrack\\
	&\leq& \gamma_3 h_{n}^{3d/2}.
\end{eqnarray*}
Now, (\ref{lem8eq1}) follows from (\ref{Varf}). \hfill$\blacksquare$
\\
~

\noindent Recall that $\eta$ is a Poisson random variable with mean $n$, independent of
$\tilde{\mathbf X}, \tilde{\mathbf X}_1,\tilde{\mathbf X}_2\ldots$. Define
\begin{eqnarray*}\label{eqnSn}
	S_{n}=\sum_{\textbf{i}\in l_n} \X_{\textbf{i},n},
\end{eqnarray*}
\begin{eqnarray*}
	\U_n=\frac{1}{\sqrt{n}}\left[\sum_{i=1}^{\eta} \mathds{1}_{\left\lbrace \tilde{\mathbf X}_i\in \M\right\rbrace }-n \mathbb P\big(\tilde{\mathbf X}\in \mathcal{M}\big)\right]
\end{eqnarray*}
and
\begin{eqnarray*}
	\V_n=\frac{1}{\sqrt{n}}\left[\sum_{i=1}^{\eta}\mathds{1}_{\left\lbrace \tilde{\mathbf X}_i\notin\M\right\rbrace } -n\mathbb P\big(\tilde{\mathbf X}\notin\M \big)\right].
\end{eqnarray*}
It is obvious to see that the random vector $(S_{n},\U_{n})$ and the
random variable $\V_{n}$ are independent. Observe also that
\begin{eqnarray}\label{eqnC1}
\mathrm{Var}( S_{n})=1 ~~ \mbox{ and } ~~ \mathrm{Var}(\U_n)=\prod_{l=1}^{p}(1-\mathbf u_l).
\end{eqnarray}
\begin{lemma}\label{lemC}
	Whenever
	$h_n\rightarrow 0$, $nh_n^{d}\rightarrow\infty$ 
	and, for each $l=1,\ldots,p$,  $C_l$ satisfies \textnormal{(\ref{lem1eqn3})} and \textnormal{(\ref{lem1eqn4})} of Lemma \ref{lem1} with $g(\cdot)=f_l(\cdot)$ and $\HH(\cdot)=\HH_l(\cdot)$, there exists a constant $\tau$ such that, for all $n$ sufficiently large,
	\begin{eqnarray}\label{eqnC}
	|\mathrm{Cov}(S_{n},
	\U_{n})|\leq \frac{\tau}{\sqrt{n}h_{n}^{d}}.
	\end{eqnarray}
\end{lemma}
\noindent{\bf{Proof of Lemma \ref{lemC}.}}  Observe that
\begin{eqnarray*}
	&&\Big(n\mathrm{Var}\big(U_{\eta,0}(\bar{C})\big)\Big)^{1/2} |\mathrm{Cov}(S_{n},
	\U_{n})|
	=\left| \mathrm{Cov}\Big(
	\sqrt{n}
	\int_{\bar{C}}
	|\Gamma_{0,\eta}(\tilde{\mathbf x})|
	\, d\tilde{\mathbf x}, \U_{n}\Big)\right|.
\end{eqnarray*}
In order to get (\ref{eqnC}), it is sufficient, by using (\ref{Varf}), to prove that there exists a constant $\tau_1$ such that, for all $n$ sufficiently large,
\begin{eqnarray}\label{eqn1lem9}
\left| \mathrm{Cov}\left(
\sqrt{n}
\int_{\bar{C}}
|\Gamma_{\eta,0}(\tilde{\mathbf x})|
\, d\tilde{\mathbf x}, \U_{n}\right)\right| \leq \frac{\tau_1}{\sqrt{n}h_{n}^{d}}.
\end{eqnarray}
Now, for any $\tilde{\mathbf x}\in \bar{C}$, 
\begin{eqnarray*}
	\left(\frac{\sqrt{n}\Gamma_{\eta,0}(\tilde{\mathbf x})}{\sqrt{	k_{n,0}(\tilde{\mathbf x})}},\frac{\U_n}{\sqrt{\mathbb P\big(\tilde{\mathbf x}\in
			\mathcal{M}\big)}}\right)\stackrel{\D}{=}\frac{1}{\sqrt{n}} \sum_{i=1}^{n} \big(\Y_{0,n}^{(i)}(\tilde{\mathbf x}), \U^{(i)}\big),
\end{eqnarray*}
where $\big(\Y_{0,n}^{(i)}(\tilde{\mathbf x}), \U^{(i)}\big)_{1\leq i \leq n}$ are independent and $\big(\Y_{0,n}(\tilde{\mathbf x}),\U\big)$- identically distributed random vectors with
\begin{eqnarray*}
	\Y_{0,n}(\tilde{\mathbf x})=\frac{1}{h_{n}^{d}\sqrt{ k_{n,0}(\tilde{\mathbf x})}}\left[ \sum_{j\leq \eta_1}\prod_{l=1}^{p}K_{l}\bigg(\frac{\mathbf x_l-\mathbb X^{l}_j}{h_n}\bigg) - \prod_{l=1}^{p}\mathbb E K_{l}\bigg(\frac{\mathbf x_l-\mathbb X^{l}}{h_n}\bigg) \right]
\end{eqnarray*}
and
\begin{eqnarray*}
	\U=\left[ \sum_{j\leq \eta_{1}} \mathds{1}_{\left\lbrace \tilde{\mathbf x}_i\in \M\right\rbrace }-\mathbb P\big(\tilde{\mathbf x}\in \M \big)\right] \big /\sqrt{\mathbb P\big(\tilde{\mathbf x}\in \mathcal{M}\big)}.
\end{eqnarray*}
Here, as above, $\eta_1$ denotes a Poisson random variable with mean $1$, independent of $\tilde{\mathbf X}_1,\tilde{\mathbf X}_2\ldots$. Notice that $\mathbb E \Y_{0,n}(\tilde{\mathbf x})=\mathbb E \U=0$,
$\mathrm{Var}(\Y_{0,n}(\tilde{\mathbf x}))=\mathrm{Var} (\U)=1$ and
\begin{eqnarray*}\Big|\mathrm{Cov}\big(\Y_{0,n}(\tilde{\mathbf x}), \U\big)\Big| &=&\left|\frac{\displaystyle \mathbb E\left[\prod_{l=1}^{p} K_{l}\Big(\frac{\mathbf x_l-\mathbb X^{l}}{h_n}\Big) \mathds 1_{\{\tilde{\mathbf x}\in \mathcal{M}\}}\right] }{\displaystyle h_{n}^{d}\sqrt{ k_{n,0}(\tilde{\mathbf x})} \sqrt{\mathbb P\big(\tilde{\mathbf x}\in \mathcal{M}\big)}  }\right|\\
	&\leq& \frac{\displaystyle \prod_{l=1}^{p} \mathbb EK_{l}\Big(\frac{\mathbf x_l-\mathbb X^{l}}{h_n}\Big)
	}{\displaystyle h_{n}^{d}\sqrt{ k_{n,0}(\tilde{\mathbf x})} \sqrt{\mathbb P\big(\tilde{\mathbf x}\in \mathcal{M}\big)} }.
\end{eqnarray*}
Therefore, by (\ref{l4eqn2}) and (\ref{majora1}), this last term is, for all large enough $n$
uniformly in $\tilde{\mathbf x}\in \bar{C}$ for some real number $\tau_2>0$, less than or equal to $ \tau_2 \sqrt{ h_{n}^{d} }$, which, in turn, is less than or equal to $\varepsilon$ for all large
enough $n$ and any $ 0 < \varepsilon <1$. This fact, in combination with (\ref{l4eqn4}), gives, by using (\ref{equa6.20Fact6.1}) and Lemma 2.3 of \cite{mason2001}, that, for some constant $\tau_3>0$, uniformly in $\tilde{\mathbf x}\in \bar{C}$,
\begin{eqnarray*}\label{lem6eq6}
	\lefteqn{\Big|\mathrm{Cov}\big(\sqrt{n} \big|f_\eta(\tilde{\mathbf x})-\mathbb E f_n(\tilde{\mathbf x})\big|,\U_n\big)\Big|}\\\nonumber
	&&\qquad\quad =\left|\mathrm{Cov}\left(\frac{\displaystyle \left|\sum_{i=1}^{n}\Y_{0,n}^{(i)}(\tilde{\mathbf x})\right|}{\displaystyle \sqrt{n}}, \frac{\displaystyle\sum_{i=1}^{n}\U^{(i)}}{\displaystyle \sqrt{n}}\right)\right| \sqrt{ k_{n,0}(\tilde{\mathbf x})\mathbb P\big(\tilde{\mathbf X}\in \mathcal{M}\big)}\\\nonumber
	&&\qquad\quad \leq \frac{\tau_3}{\sqrt{nh_{n}^{d}}} \sqrt{ k_{n,0}(\tilde{\mathbf x})},
\end{eqnarray*}
which, when combined with (\ref{l4eqn2}) and $\lambda(\bar{C})<\infty$, gives (\ref{eqn1lem9}). This finishes the proof of Lemma \ref{lemC}. \hfill$\blacksquare$
\\
~

\begin{lemma}\label{lem6}
	Suppose that, for each $l=1,\ldots,p$,  $C_l$ satisfies \textnormal{(\ref{lem1eqn3})} and \textnormal{(\ref{lem1eqn4})} of Lemma \ref{lem1} with $g(\cdot)=f_l(\cdot)$ and $\HH(\cdot)=\HH_l(\cdot)$. If
	$h_n\rightarrow 0$ and $\sqrt{n}h_n^{d}\rightarrow\infty$, then  
	\begin{eqnarray}\label{lem6eq41}
	( S_n,\U_n) \stackrel{\D}{\rightarrow}\left(Z_1,\prod_{l=1}^{p}\sqrt{(1-\mathbf u_l)}Z_2 \right)
	\end{eqnarray}
	as $n\rightarrow\infty$, where $Z_1$ and $Z_2$ are independent standard normal random variables.
\end{lemma}
\noindent{\bf{Proof of Lemma \ref{lem6}.}}  We will show that, for any $\xi_{1}$ and $\xi_{2}$, as $n\rightarrow \infty$,
\begin{eqnarray*}\label{lem6eq1}
	\xi_{1}S_n+\xi_{2} \U_{n} \stackrel{\D}{\rightarrow} \xi_{1} Z_1+\xi_2\prod_{l=1}^{p}\sqrt{(1-\mathbf u_l)}Z_2. 
\end{eqnarray*}
For each $\textbf{i}\in l_n$, set
\begin{eqnarray*}
	\U_{\textbf{i},n} :=\frac{1}{\sqrt{n}}\left[\sum_{j=1}^{\eta}\mathds{1}_{\left\lbrace \tilde{\mathbf X}_i\in R_{\textbf{i}}\right\rbrace } -n \mathbb P\big(\tilde{\mathbf X}\in R_\textbf{i}\big)\right]
\end{eqnarray*}
and
\begin{eqnarray*}
	\ZZ_{\textbf{i},n} :=\xi_1 \X_{\textbf{i},n} +\xi_2 \U_{\textbf{i},n}.
\end{eqnarray*}
By $c_r$-inequality, for all $r\geq 1$, we have 
\begin{eqnarray}\label{eqn4lem10}
\mathbb E|\ZZ_{i,n}|^r\leq 2^{r-1}\Big(|\xi_1|^r \mathbb E|\X_{i,n}|^r+|\xi_2|^r \mathbb E|\U_{i,n}|^r\Big).
\end{eqnarray}
Now, by H\"{o}lder inequality, for all $2<r<3$ and all $\textbf{i}\in l_n$,
\begin{eqnarray*}
	\mathbb E|\X_{\textbf{i},n}|^r\leq \big(\mathbb E|\X_{\textbf{i},n}|^3\big)^{r/3},
\end{eqnarray*}
which by (\ref{lem6eq2b}) and (\ref{lem8eq1}) implies that, for all $2<r<3$
\begin{eqnarray} \label{eqn5lem10}
\sum_{\textbf{i}\in l_n}\mathbb E|\X_{\textbf{i},n}|^r \leq 2^{d} \gamma_{0}^{r/3} h_n^{d(r-2)/2}\prod_{l=1}^{p}M_l^{d_l}\rightarrow0.
\end{eqnarray}
Set now
\begin{eqnarray*}
	p_{\textbf{i},n}=\mathbb P\big(\tilde{\mathbf X}\in R_{\textbf{i}}\big).
\end{eqnarray*}
By Lemma 2.3 in \cite{mason2001}, we can find a constant $\zeta>0$ such that
\begin{eqnarray}\label{eqn6lem10}
\sum_{\textbf{i}\in l_n}\mathbb E|\U_{\textbf{i},n}|^r &\leq& \zeta n^{-r/2}\sum_{\textbf{i}\in l_n} \left((np_{\textbf{i},n})^{r/2}+
np_{\textbf{i},n}\right)\nonumber\\
&\leq& \zeta\max_{\textbf{i}\in l_n} \left((p_{\textbf{i},n})^{(r-2)/2}+ n^{(2-r)/2}\right)\rightarrow 0. 
\end{eqnarray}
Combining (\ref{eqn4lem10})-(\ref{eqn6lem10}), we obtain, for all $2<r<3$,
\begin{eqnarray*}\label{eqn6lem11}
	\lim_{n\rightarrow \infty} \sum_{\textbf{i}\in l_n}\mathbb E|\ZZ_{\textbf{i},n}|^r =0. 
\end{eqnarray*}
Moreover, note that $\{\ZZ_{\textbf{i},n}:\textbf{i}\in l_n \}$ is a triangular array of mean zero
one-dependent random fields. Consequently, by the statements (\ref{eqnC1}) and (\ref{eqnC}), it follows that
\begin{eqnarray*}\label{lem6eq3}
	\mathrm{Var}\left(\sum_{\textbf{i}\in l_n}
	\ZZ_{\textbf{i},n}\right)= \mathrm{Var}\big (\xi_{1}S_n+\xi_2 \U_n\big)\rightarrow \xi_1^2+\xi_2^2 \prod_{l=1}^{p}(1-\mathbf u_l)\quad \text{as}\, n\rightarrow \infty.
\end{eqnarray*}
Thus, by Theorem $1$ of \cite{shergin}, we can infer that 
\begin{eqnarray*}
	\sum_{\textbf{i}\in l_n} \ZZ_{\textbf{i},n}\stackrel{\D}{\rightarrow} \xi_{1} Z_1+\xi_2\prod_{l=1}^{p}\sqrt{(1-\mathbf u_l)}Z_2 \quad \text{as}\, n\rightarrow \infty.
\end{eqnarray*}
Since 
\begin{eqnarray*}
	\sum_{\textbf{i}\in l_n} \ZZ_{\textbf{i},n}=\xi_{1}S_n+\xi_{2} \U_{n},
\end{eqnarray*} 
the statement (\ref{lem6eq41}) is proved by the Cram\'er-Wold device, for instance, see, \cite{billy}.
\hfill $\blacksquare$\\
~

\noindent
\begin{lemma}\label{lem7}
	Suppose that, for each $l=1,\ldots,p$,  $C_l$ satisfies \textnormal{(\ref{lem1eqn3})} and \textnormal{(\ref{lem1eqn4})} of Lemma \ref{lem1} with $g(\cdot)=f_l(\cdot)$ and $\HH(\cdot)=\HH_l(\cdot)$. If
	$h_n\rightarrow 0$ and $nh_n^{3d}\rightarrow\infty$, then 
	\begin{eqnarray}\label{eqnnor}
	\frac{\sqrt{n}\big(V_{n}(\bar{C})-\mathbb E V_{n}(\bar{C})\big)}{\sqrt{n\mathrm{Var}\big(U_{\eta,0}(\bar{C})\big)}}\stackrel{\D}{\rightarrow} Z\qquad \mbox{ as } n\rightarrow\infty,
	\end{eqnarray}
	where $Z$ is a standard normal random variable.
\end{lemma}
\noindent{\bf{Proof of Lemma \ref{lem7}.}} 
Note that
\begin{eqnarray*}
	S_n=\frac{\displaystyle \sqrt{n}\big(U_{\eta,0}(\bar{C})-\mathbb E U_{\eta,0}(\bar{C})\big)}{\displaystyle \sqrt{n\mathrm{Var}\big(U_{\eta,0}(\bar{C})\big)}}.
\end{eqnarray*}
Indeed, conditioned on $\eta=n$, we have
\begin{eqnarray*}
	S_n=\frac{\displaystyle\sqrt{n}\big(U_{n,0}(\bar{C})-\mathbb E U_{\eta,0}(\bar{C})\big)}{\displaystyle\sqrt{n\mathrm{Var}\big(U_{\eta,0}(\bar{C})\big)}}.
\end{eqnarray*}
By Lemma $2.4$ of \cite{mason2001} refer also to \cite{{mason1995}} and Lemma \ref{lem6}, we conclude that
\begin{eqnarray*}
	\frac{\displaystyle\sqrt{n}\big(U_{n,0}(\bar{C})-\mathbb E U_{\eta,0}(\bar{C})\big)}{\displaystyle\sqrt{n\mathrm{Var}\big(U_{\eta,0}(
			\bar{C})\big)}}\stackrel{\D}{\rightarrow} Z.
\end{eqnarray*}
Therefore, by (\ref{Varf}) and Lemma \ref{lem4}, we have
\begin{eqnarray*}
	\frac{\sqrt{n}\big(U_{n,0}(\bar{C})-\mathbb E U_{n,0}(\bar{C})\big)}{\sqrt{n\mathrm{Var}\big(U_{\eta,0}(\bar{C})\big)}}\stackrel{\D}{\rightarrow} Z.
\end{eqnarray*}
Assertion (\ref{eqnnor}) now follows from Lemmas \ref{lem3} and \ref{lem51} in combination with (\ref{Varf}). Hence the proof is complete. \hfill$\blacksquare$
\\
~

\begin{lemma}\label{Lemma12}
	Assume that \textnormal{({\bf A$_1$})} holds and, for each $ l=1, \ldots,p$, there exists $\delta>0$ such that $f_l(\cdot)$ is bounded on $ \C_{l}^{\delta}$. If $h_{n}\rightarrow0$ and $nh_{n}^{3d}\rightarrow\infty$, then we have 
	\begin{eqnarray}\label{TH2eqn}
	\lim_{n\rightarrow \infty}\mathbb E\int_{\tilde{\C}}\left|\sqrt{\hat{\mathcal{L}}_{n}(\tilde{\mathbf x})
	}-\sqrt{\mathcal{L}_{n,1}(\tilde{\mathbf x})}\right|d\tilde{\mathbf x}=0.
	\end{eqnarray}
\end{lemma}
\noindent{\bf Proof of Lemma \ref{Lemma12}.} Notice that
\begin{eqnarray*}
	\mathcal{L}_{n,1}(\tilde{\mathbf x})=\sum_{\{I\subsetneq I_p,|I_p\setminus I|\geq2\}}\prod_{l\in I }\big(\mathbb Ef_{n,l}(\mathbf x_l)\big)^{2}\prod_{j\in I_p\setminus I}\big\lbrack v_{n,j}(\mathbf x_j)-\big(\mathbb E f_{n,j}(\mathbf x_j)\big)^{2}\big\rbrack
\end{eqnarray*}
and 
\begin{eqnarray*}
	\hat{\mathcal{L}}_{n}(\tilde{\mathbf x})=\sum_{\{I\subsetneq I_p,|I_p\setminus I|\geq2\}}\prod_{l\in I }g_{n,l}(\mathbf x_l)\prod_{j\in I_p\setminus I}\big\lbrack \hat{v}_{n,j}(\mathbf x_j)- g_{n,j}(\mathbf x_j)\big\rbrack.
\end{eqnarray*} 
Set, for each $I\subsetneq I_p$,
\begin{eqnarray*}
	\boldsymbol{\Theta}_{n,1,I}(\tilde{\mathbf x}):=\left|\prod_{j\in I_p\setminus I}\big\lbrack v_{n,j}(\mathbf x_j)-\big(\mathbb E f_{n,j}(\mathbf x_j)\big)^{2}\big\rbrack-
	\prod_{j\in I_p\setminus I}\big\lbrack \hat{v}_{n,j}(\mathbf x_j)- g_{n,j}(\mathbf x_j)\big\rbrack\right|
\end{eqnarray*} 
and
\begin{eqnarray*}
	\boldsymbol{\Theta}_{n,2,I}(\tilde{\mathbf x}):=\left|\prod_{l\in I }g_{n,l}(\mathbf x_l)-
	\prod_{l\in I }   
	\big(\mathbb Ef_{n,l}(\mathbf x_l)\big)^{2}\right|.
\end{eqnarray*} 
It is obvious to see that
\begin{eqnarray}\label{TH2eqn11}
\big|\mathcal{L}_{n,1}(\tilde{\mathbf x})-\hat{\mathcal{L}}_{n}(\tilde{\mathbf x})
\big|&\leq&\sum_{\{I\subsetneq I_p,|I_p\setminus I|\geq2\}} \prod_{l\in I }g_{n,l}(\mathbf x_l) \boldsymbol{\Theta}_{n,1,I}(\tilde{\mathbf x})\nonumber \\&&+\sum_{\{I\subsetneq I_p,|I_p\setminus I|\geq2\}}\prod_{j\in I_p\setminus I} v_{n,j}(\mathbf x_j) \boldsymbol{\Theta}_{n,2,I}(\tilde{\mathbf x})\nonumber\\&=:&\Psi_{n,1}(\tilde{\mathbf x})+\Psi_{n,2}(\tilde{\mathbf x}).
\end{eqnarray}
For each $I\subsetneq I_p$, observe that
\begin{eqnarray}\label{TH2eqn1}
\nonumber\mathbb E\boldsymbol{\Theta}_{n,1,I}(\tilde{\mathbf x})\leq\sum_{A\subsetneq I_p\setminus I} \prod_{l\in A} v_{n,l}(\mathbf x_l)
\prod_{j\in I_p\setminus (I\cup A)}\Big\{\mathbb E| \hat{v}_{n,j}(\mathbf x_j)-v_{n,j}(\mathbf x_j)|+\mathbb E\big| g_{n,j}(\mathbf x_j)-\big(\mathbb E f_{n,j}(\mathbf x_j)\big)^{2}
\big|\Big\}\\
\end{eqnarray} 
and 
\begin{eqnarray}\label{TH2eqn12}
\mathbb E\boldsymbol{\Theta}_{n,2,I}(\tilde{\mathbf x})\leq
\sum_{A\subsetneq I} \prod_{l\in A} \big(\mathbb Ef_{n,l}(\mathbf x_l)\big)^{2}
\prod_{j\in I\setminus A}\mathbb E\left| g_{n,j}(\mathbf x_j)-\big(\mathbb Ef_{n,l}(\mathbf x_j)\big)^{2}\right|.
\end{eqnarray} 
Now, for each $l\in I_p$, we have 
\begin{eqnarray*}
	\mathbb E\big| \hat{v}_{n,l}(\mathbf x_l)-v_{n,l}(\mathbf x_l)\big|&\leq&\Big( \mathbb E\big| \hat{v}_{n,l}(\mathbf x_l)-v_{n,l}(\mathbf x_l)\big|^{2}\Big)^{1/2}\\&\leq&\frac{1}{\sqrt{n}
		h_{n}^{2d_{l}}}\left(\int_{\R^{d_l}} K_{l}^{4}\Bigg(\frac{\mathbf x_l-\mathbf u_l}{h_{n}}\Bigg)f_{l}(\mathbf u_l)\,d\mathbf u_l\right)^{1/2}.
\end{eqnarray*}
Since, for each $l\in I_p$, $K_l(\cdot)$ is bounded on $\R^{d_l}$ and $f_l(\cdot)$ is bounded on $\C_l^{h_n/2}$ for $n$ large enough, we obtain 
\begin{eqnarray}\label{TH2eqn2}
\sup_{x_l\in\C_l}\mathbb E\big| \hat{v}_{n,l}(\mathbf x_l)-v_{n,l}(\mathbf x_l)\big|=O\left(\frac{1}{\sqrt{nh_{n}^{3d_{l}}}}\right).
\end{eqnarray}
Observe that
\begin{eqnarray*}
	\lefteqn{\mathbb E\bigg| g_{n,l}(\mathbf x_l)-\big(\mathbb E f_{n,l}(\mathbf x_l)\big)^{2}\bigg|}\\&\leq&\frac{1}{n(n-1)
		h_{n}^{2d_{l}}}\mathbb E\left|\sum_{i=1}^{n}\sum_{j\neq i}K_{l}\Bigg(\frac{\mathbf x_l-\mathbb X_{i}^{l}}{h_{n}}\Bigg)\Bigg\lbrack
	K_{l}\Bigg(\frac{\mathbf x_l-\mathbb X^{l}_{j}}{h_{n}}\Bigg)-\mathbb EK_{l}\Bigg(\frac{\mathbf x_l-
		\mathbb X^{l}_{j}}{h_{n}}\Bigg)
	\Bigg \rbrack\right.\\&&\left.+\frac{1}{n(n-1)
		h_{n}^{2d_{l}}}\mathbb E\Bigg|\sum_{i=1}^{n}\sum_{j\neq i}\mathbb EK_{l}\Bigg(\frac{\mathbf x_l-\mathbb X_{j}^{l}}{h_{n}}\Bigg)\Bigg\lbrack
	K_{l}\Bigg(\frac{\mathbf x_l-\mathbb X^{l}_{i}}{h_{n}}\Bigg)-\mathbb EK_{l}\Bigg(\frac{\mathbf x_l-
		\mathbb X^{l}_{i}}{h_{n}}\Bigg)
	\Bigg \rbrack\right|.
\end{eqnarray*}
For each $l\in I_p$, since $K_l(\cdot)$ is bounded on $\R^{d_l}$, it follows that, for some constant $\mathfrak H_0>0$, we have 
\begin{eqnarray*}
	\mathbb E\left| g_{n,l}(\mathbf x_l)-\big(\mathbb E f_{n,l}(\mathbf x_l)\big)^{2}\right|&\leq&\frac{\mathfrak H_0}{n
		h_{n}^{2d_{l}}}\mathbb E\left|\sum_{i=1}^{n}\Bigg\lbrack
	K_{l}\Bigg(\frac{\mathbf x_l-\mathbb X^{l}_{i}}{h_{n}}\Bigg)-\mathbb EK_{l}\Bigg(\frac{\mathbf x_l-\mathbb 
		X^{l}_{i}}{h_{n}}\Bigg)
	\Bigg \rbrack\right|\\&\leq&\frac{\mathfrak H_0}{\sqrt{n}
		h_{n}^{2d_{l}}}\left(\int_{\R^{d_l}} K_{l}^{2}\left(\frac{\mathbf x_l-\mathbf u_l}{h_{n}}\right)f_{l}(\mathbf u_l)\,d\mathbf u_l\right)^{1/2}.
\end{eqnarray*}
Since, for each $l\in I_p$, $f_l(\cdot)$ is bounded on $\C_l^{h_n/2}$ for $n$ large enough, we have
\begin{eqnarray}\label{TH2eqn3}
\sup_{\mathbf x_l\in\C_l} \mathbb E\big| g_{n,l}(\mathbf x_l)-\big(\mathbb E f_{n,l}(\mathbf x_l)\big)^{2}\big|=O\left(\frac{1}{\sqrt{nh_{n}^{3d_{l}}}}\right).
\end{eqnarray}
By the same arguments, we get
\begin{eqnarray}\label{TH2eqn4}
\sup_{\mathbf x_l\in\C_l} v_{n,l}(\mathbf x_l)\leq\frac{\mathfrak H_1}
{h_{n}^{d_{l}}},\,\,
\,\, \sup_{\mathbf x_l\in\C_l} \mathbb Eg_{n,l}(\mathbf x_l)\leq\mathfrak H_1 \,\, \mbox{ and }\,\,
\sup_{\mathbf x_l\in\C_l} \mathbb Ef_{n,l}(\mathbf x_l)\leq \mathfrak H_2,
\end{eqnarray}
for some positive constants $\mathfrak H_1, \mathfrak H_2$ and $\mathfrak H_3$. Therefore, by (\ref{TH2eqn1})-(\ref{TH2eqn4}), we infer that 

\begin{eqnarray}\label{TH2eqn51}
\sup_{\tilde{\mathbf x}\in \tilde{\C} }\mathbb E\Psi_{n,1}(\tilde{\mathbf x})=O\left(\frac{1}
{\sqrt{nh_n^{3d}}}\right).
\end{eqnarray}
Similarly, by (\ref{TH2eqn12}) in combination with (\ref{TH2eqn3}) and (\ref{TH2eqn4}), we obtain 
\begin{eqnarray}\label{TH2eqn52}
\sup_{\tilde{\mathbf x}\in \tilde{\C} }\mathbb E\Psi_{n,2}(\tilde{\mathbf x})=O\left(\frac{1}
{\sqrt{nh_n^{3d}}}\right).
\end{eqnarray}
By using (\ref{TH2eqn11}), (\ref{TH2eqn51}) and (\ref{TH2eqn52}) in combination with the fact that $\tilde{\C}$ is bounded, we get 
\begin{eqnarray*}
	\mathbb E\int_{\tilde{\C}}\left|\sqrt{\mathcal{L}_{n,1}(\tilde{\mathbf x})}-\sqrt{\hat{\mathcal{L}}_{n}(\tilde{\mathbf x})
	}\right|d\tilde{\mathbf x}\nonumber&\leq&\int_{\bar{\C}}\sqrt{\mathbb E|\mathcal{L}_{n,1}(\tilde{\mathbf x})-\hat{\mathcal{L}}_{n}(\tilde{\mathbf x})
		|}d\tilde{\mathbf x}\\&=&O\left(\frac{1}
	{\sqrt{nh_n^{3d}}}\right).
\end{eqnarray*}
We so obtain (\ref{TH2eqn}), as sought.\hfill$\blacksquare$

\renewcommand{\thesection}{C}
\section{}
For completeness, 
 we recall the following fact (Fact 6.1 of \cite{mason2001}), which follows from Theorem 1 of \cite{sweeting} and is related to the classical Berry-Esseen theorem.

\begin{theorem} \label{Thesweeting}Let $(\omega, \zeta),\left(\omega_1, \zeta_1\right),\left(\omega_2, \zeta_2\right), \ldots$ be a sequence of i.i.d. random vectors such that each component has variance 1, mean 0, and finite absolute moments of the third order. Further, let $\left(\bar{Z}_1, \bar{Z}_2\right)$ be bivariate normal with mean vector $0,$$ \operatorname{Var}\left(\bar{Z}_1\right)=\operatorname{Var}\left(\bar{Z}_2\right)=1, \mbox{~and ~} \operatorname{Cov}\left(\bar{Z}_1, \bar{Z}_2\right)=\operatorname{Cov}(\omega, \zeta)=\rho.$ Then there exist universal positive constants $A_1, A_2$ and $A_3$ such that
	\begin{equation}
	\left|\mathbb E\left| \frac{\sum_{i=1}^n \zeta_i}{\displaystyle \sqrt{n}}\right|-\mathbb E| \bar{Z}_1|\right| \leq \frac{A_1}{\sqrt{n}} \mathbb E|\zeta|^3, \label{equa6.18Fact6.1}
	\end{equation}
	and, whenever $\rho^2<1$,
	\begin{equation}\left |\mathbb E\left| \frac{ \sum_{i=1}^n \omega_i}{\displaystyle \sqrt{n}} \cdot \frac{\sum_{i=1}^n \zeta_i}{\displaystyle \sqrt{n}}\right|-\mathbb E| \bar{Z}_1 \bar{Z}_2|\right| \leq \frac{A_2}{\left(1-\rho^2\right)^{3 / 2}} \frac{1}{\sqrt{n}}\left(\mathbb E|\omega|^3+\mathbb E|\zeta|^3\right),
	\label{equa6.19Fact6.1}
	\end{equation}
	and
	\begin{equation}
	\left|\mathbb E\left[\frac{ \sum_{i=1}^n \omega_i}{\displaystyle \sqrt{n}} \left|\frac{ \sum_{i=1}^n \zeta_i}{\displaystyle \sqrt{n}}\right|\right]\right| \leq \frac{A_3}{\left(1-\rho^2\right)^{3 / 2}} \frac{1}{\sqrt{n}}\left(\mathbb E|\omega|^3+\mathbb E|\zeta|^3\right).
	\label{equa6.20Fact6.1}
	\end{equation}
\end{theorem}


\begin{thebibliography}{108}

\bibitem[\protect\citeauthoryear{Ahmad and Li}{1997}]{Ahmad1997}
\begin{barticle}[author]
\bauthor{\bsnm{Ahmad},~\bfnm{Ibrahim~A.}\binits{I.~A.}} \AND
  \bauthor{\bsnm{Li},~\bfnm{Qi}\binits{Q.}}
(\byear{1997}).
\btitle{Testing independence by nonparametric kernel method}.
\bjournal{Statist. Probab. Lett.}
\bvolume{34}
\bpages{201--210}.
\bdoi{10.1016/S0167-7152(96)00183-6}
\end{barticle}
\endbibitem

\bibitem[\protect\citeauthoryear{Akaike}{1954}]{Akaike1954}
\begin{barticle}[author]
\bauthor{\bsnm{Akaike},~\bfnm{Hirotugu}\binits{H.}}
(\byear{1954}).
\btitle{An approximation to the density function}.
\bjournal{Ann. Inst. Statist. Math., Tokyo}
\bvolume{6}
\bpages{127--132}.
\bdoi{10.1007/bf02900741}
\end{barticle}
\endbibitem

\bibitem[\protect\citeauthoryear{Araujo and Gin\'{e}}{1980}]{MR576407}
\begin{bbook}[author]
\bauthor{\bsnm{Araujo},~\bfnm{Aloisio}\binits{A.}} \AND
  \bauthor{\bsnm{Gin\'{e}},~\bfnm{Evarist}\binits{E.}}
(\byear{1980}).
\btitle{The {C}entral {L}imit {T}heorem for {R}eal and {B}anach {V}alued
  {R}andom {V}ariables}.
\bseries{Wiley Series in Probability and Mathematical Statistics}.
\bpublisher{John Wiley \& Sons, New York-Chichester-Brisbane}.
\end{bbook}
\endbibitem

\bibitem[\protect\citeauthoryear{Auddy, Deb and Nandy}{2024}]{auddy2023exact}
\begin{barticle}[author]
\bauthor{\bsnm{Auddy},~\bfnm{Arnab}\binits{A.}},
  \bauthor{\bsnm{Deb},~\bfnm{Nabarun}\binits{N.}} \AND
  \bauthor{\bsnm{Nandy},~\bfnm{Sagnik}\binits{S.}}
(\byear{2024}).
\btitle{Exact detection thresholds and minimax optimality of {C}hatterjee's
  correlation coefficient}.
\bjournal{Bernoulli}
\bvolume{30}
\bpages{1640--1668}.
\bdoi{10.3150/23-bej1648}
\end{barticle}
\endbibitem

\bibitem[\protect\citeauthoryear{Bach and Jordan}{2003}]{Bach}
\begin{barticle}[author]
\bauthor{\bsnm{Bach},~\bfnm{Francis~R.}\binits{F.~R.}} \AND
  \bauthor{\bsnm{Jordan},~\bfnm{Michael~I.}\binits{M.~I.}}
(\byear{2003}).
\btitle{Kernel independent component analysis}.
\bjournal{J. Mach. Learn. Res.}
\bvolume{3}
\bpages{1--48}.
\bdoi{10.1162/153244303768966085}
\end{barticle}
\endbibitem

\bibitem[\protect\citeauthoryear{Beirlant and Mason}{1995}]{mason1995}
\begin{barticle}[author]
\bauthor{\bsnm{Beirlant},~\bfnm{J.}\binits{J.}} \AND
  \bauthor{\bsnm{Mason},~\bfnm{D.~M.}\binits{D.~M.}}
(\byear{1995}).
\btitle{On the asymptotic normality of {$L\sb p$}-norms of empirical
  functionals}.
\bjournal{Math. Methods Statist.}
\bvolume{4}
\bpages{1--19}.
\end{barticle}
\endbibitem

\bibitem[\protect\citeauthoryear{Berger}{1991}]{Berg}
\begin{barticle}[author]
\bauthor{\bsnm{Berger},~\bfnm{Erich}\binits{E.}}
(\byear{1991}).
\btitle{Majorization, exponential inequalities and almost sure behavior of
  vector-valued random variables}.
\bjournal{Ann. Probab.}
\bvolume{19}
\bpages{1206--1226}.
\bdoi{10.1214/aop/1176990341}
\end{barticle}
\endbibitem

\bibitem[\protect\citeauthoryear{Berrahou, Bouzebda and
  Douge}{2024}]{L1testBBD}
\begin{barticle}[author]
\bauthor{\bsnm{Berrahou},~\bfnm{Noureddine}\binits{N.}},
  \bauthor{\bsnm{Bouzebda},~\bfnm{Salim}\binits{S.}} \AND
  \bauthor{\bsnm{Douge},~\bfnm{Lahcen}\binits{L.}}
(\byear{2024}).
\btitle{Supplementary materials to ``A nonparametric distribution-free test of
  independence among continuous random vectors based on {$L_1$}-norm''}.
\end{barticle}
\endbibitem

\bibitem[\protect\citeauthoryear{Berrett}{2022}]{MR4472835}
\begin{barticle}[author]
\bauthor{\bsnm{Berrett},~\bfnm{T.~B.}\binits{T.~B.}}
(\byear{2022}).
\btitle{Discussion of `{M}ulti-scale {F}isher's independence test for
  multivariate dependence' [ 4472834]}.
\bjournal{Biometrika}
\bvolume{109}
\bpages{589--592}.
\bdoi{10.1093/biomet/asac023}
\end{barticle}
\endbibitem

\bibitem[\protect\citeauthoryear{Berrett, Grose and Samworth}{2018}]{MINT1}
\begin{barticle}[author]
\bauthor{\bsnm{Berrett},~\bfnm{Thomas~B.}\binits{T.~B.}},
  \bauthor{\bsnm{Grose},~\bfnm{Daniel~J.}\binits{D.~J.}} \AND
  \bauthor{\bsnm{Samworth},~\bfnm{Richard~J.}\binits{R.~J.}}
(\byear{2018}).
\btitle{\textnormal{IndepTest:} {\it{Nonparametric Independence Tests Based on
  Entropy Estimation.}}}
\bjournal{\textnormal{R Package version 0.2.0. Available at
  \url{https://cran.r-project.org/web/ packages/IndepTest/index.html.}}}
\end{barticle}
\endbibitem

\bibitem[\protect\citeauthoryear{Berrett, Kontoyiannis and
  Samworth}{2021}]{MR4338371}
\begin{barticle}[author]
\bauthor{\bsnm{Berrett},~\bfnm{Thomas~B.}\binits{T.~B.}},
  \bauthor{\bsnm{Kontoyiannis},~\bfnm{Ioannis}\binits{I.}} \AND
  \bauthor{\bsnm{Samworth},~\bfnm{Richard~J.}\binits{R.~J.}}
(\byear{2021}).
\btitle{Optimal rates for independence testing via {$U$}-statistic permutation
  tests}.
\bjournal{Ann. Statist.}
\bvolume{49}
\bpages{2457--2490}.
\bdoi{10.1214/20-aos2041}
\end{barticle}
\endbibitem

\bibitem[\protect\citeauthoryear{Berrett and Samworth}{2019}]{berrett}
\begin{barticle}[author]
\bauthor{\bsnm{Berrett},~\bfnm{T.~B.}\binits{T.~B.}} \AND
  \bauthor{\bsnm{Samworth},~\bfnm{R.~J.}\binits{R.~J.}}
(\byear{2019}).
\btitle{Nonparametric independence testing via mutual information}.
\bjournal{Biometrika}
\bvolume{106}
\bpages{547--566}.
\bdoi{10.1093/biomet/asz024}
\end{barticle}
\endbibitem

\bibitem[\protect\citeauthoryear{Billingsley}{1999}]{billy}
\begin{bbook}[author]
\bauthor{\bsnm{Billingsley},~\bfnm{Patrick}\binits{P.}}
(\byear{1999}).
\btitle{Convergence of {P}robability {M}easures},
\bedition{second} ed.
\bseries{Wiley Series in Probability and Statistics: Probability and
  Statistics}.
\bpublisher{John Wiley \& Sons, Inc., New York}
\bnote{A Wiley-Interscience Publication}.
\bdoi{10.1002/9780470316962}
\end{bbook}
\endbibitem

\bibitem[\protect\citeauthoryear{Binet and Vaschide}{1897}]{Binet}
\begin{barticle}[author]
\bauthor{\bsnm{Binet},~\bfnm{Alfred}\binits{A.}} \AND
  \bauthor{\bsnm{Vaschide},~\bfnm{Nicolas}\binits{N.}}
(\byear{1897}).
\btitle{Corrélation des épreuves physiques}.
\bjournal{L'Année psychologique}
\bvolume{4}
\bpages{142--172}.
\bdoi{10.3406/psy.1897.2892}
\end{barticle}
\endbibitem

\bibitem[\protect\citeauthoryear{Blum, Kiefer and Rosenblatt}{1961}]{Blum}
\begin{barticle}[author]
\bauthor{\bsnm{Blum},~\bfnm{J.~R.}\binits{J.~R.}},
  \bauthor{\bsnm{Kiefer},~\bfnm{J.}\binits{J.}} \AND
  \bauthor{\bsnm{Rosenblatt},~\bfnm{M.}\binits{M.}}
(\byear{1961}).
\btitle{Distribution free tests of independence based on the sample
  distribution function}.
\bjournal{Ann. Math. Statist.}
\bvolume{32}
\bpages{485--498}.
\bdoi{10.1214/aoms/1177705055}
\end{barticle}
\endbibitem

\bibitem[\protect\citeauthoryear{B\"{o}ttcher, Keller-Ressel and
  Schilling}{2019}]{bottcher2019}
\begin{barticle}[author]
\bauthor{\bsnm{B\"{o}ttcher},~\bfnm{Bj\"{o}rn}\binits{B.}},
  \bauthor{\bsnm{Keller-Ressel},~\bfnm{Martin}\binits{M.}} \AND
  \bauthor{\bsnm{Schilling},~\bfnm{Ren\'{e}~L.}\binits{R.~L.}}
(\byear{2019}).
\btitle{Distance multivariance: new dependence measures for random vectors}.
\bjournal{Ann. Statist.}
\bvolume{47}
\bpages{2757--2789}.
\bdoi{10.1214/18-AOS1764}
\end{barticle}
\endbibitem

\bibitem[\protect\citeauthoryear{Bouzebda}{2011}]{Bouzebda2011}
\begin{barticle}[author]
\bauthor{\bsnm{Bouzebda},~\bfnm{S.}\binits{S.}}
(\byear{2011}).
\btitle{Some new multivariate tests of independence}.
\bjournal{Math. Methods Statist.}
\bvolume{20}
\bpages{192--205}.
\bdoi{10.3103/S1066530711030021}
\end{barticle}
\endbibitem

\bibitem[\protect\citeauthoryear{Bouzebda}{2014}]{Bouzebda2014}
\begin{barticle}[author]
\bauthor{\bsnm{Bouzebda},~\bfnm{Salim}\binits{S.}}
(\byear{2014}).
\btitle{General tests of independence based on empirical processes indexed by
  functions}.
\bjournal{Stat. Methodol.}
\bvolume{21}
\bpages{59--87}.
\bdoi{10.1016/j.stamet.2014.03.001}
\end{barticle}
\endbibitem

\bibitem[\protect\citeauthoryear{Bouzebda and Nemouchi}{2023}]{MR4562252}
\begin{barticle}[author]
\bauthor{\bsnm{Bouzebda},~\bfnm{Salim}\binits{S.}} \AND
  \bauthor{\bsnm{Nemouchi},~\bfnm{Boutheina}\binits{B.}}
(\byear{2023}).
\btitle{Weak-convergence of empirical conditional processes and conditional
  {$U$}-processes involving functional mixing data}.
\bjournal{Stat. Inference Stoch. Process.}
\bvolume{26}
\bpages{33--88}.
\bdoi{10.1007/s11203-022-09276-6}
\end{barticle}
\endbibitem

\bibitem[\protect\citeauthoryear{Bouzebda and Taachouche}{2023}]{taachoucheB}
\begin{barticle}[author]
\bauthor{\bsnm{Bouzebda},~\bfnm{Salim}\binits{S.}} \AND
  \bauthor{\bsnm{Taachouche},~\bfnm{Nourelhouda}\binits{N.}}
(\byear{2023}).
\btitle{On the variable bandwidth kernel estimation of conditional
  {$U$}-statistics at optimal rates in sup-norm}.
\bjournal{Phys. A}
\bvolume{625}
\bpages{Paper No. 129000, 72}.
\bdoi{10.1016/j.physa.2023.129000}
\end{barticle}
\endbibitem

\bibitem[\protect\citeauthoryear{Chakraborty and Zhang}{2019}]{chakraborty2019}
\begin{barticle}[author]
\bauthor{\bsnm{Chakraborty},~\bfnm{Shubhadeep}\binits{S.}} \AND
  \bauthor{\bsnm{Zhang},~\bfnm{Xianyang}\binits{X.}}
(\byear{2019}).
\btitle{Distance metrics for measuring joint dependence with application to
  causal inference}.
\bjournal{J. Amer. Statist. Assoc.}
\bvolume{114}
\bpages{1638--1650}.
\bdoi{10.1080/01621459.2018.1513364}
\end{barticle}
\endbibitem

\bibitem[\protect\citeauthoryear{Chan, Lee and Peng}{2010}]{chan2010}
\begin{barticle}[author]
\bauthor{\bsnm{Chan},~\bfnm{Ngai-Hang}\binits{N.-H.}},
  \bauthor{\bsnm{Lee},~\bfnm{Thomas C.~M.}\binits{T.~C.~M.}} \AND
  \bauthor{\bsnm{Peng},~\bfnm{Liang}\binits{L.}}
(\byear{2010}).
\btitle{On nonparametric local inference for density estimation}.
\bjournal{Comput. Statist. Data Anal.}
\bvolume{54}
\bpages{509--515}.
\bdoi{10.1016/j.csda.2009.09.021}
\end{barticle}
\endbibitem

\bibitem[\protect\citeauthoryear{Chatterjee}{2021}]{MR4353729}
\begin{barticle}[author]
\bauthor{\bsnm{Chatterjee},~\bfnm{Sourav}\binits{S.}}
(\byear{2021}).
\btitle{A new coefficient of correlation}.
\bjournal{J. Amer. Statist. Assoc.}
\bvolume{116}
\bpages{2009--2022}.
\bdoi{10.1080/01621459.2020.1758115}
\end{barticle}
\endbibitem

\bibitem[\protect\citeauthoryear{Chen and Bickel}{2006}]{chen}
\begin{barticle}[author]
\bauthor{\bsnm{Chen},~\bfnm{Aiyou}\binits{A.}} \AND
  \bauthor{\bsnm{Bickel},~\bfnm{Peter~J.}\binits{P.~J.}}
(\byear{2006}).
\btitle{Efficient independent component analysis}.
\bjournal{Ann. Statist.}
\bvolume{34}
\bpages{2825--2855}.
\bdoi{10.1214/009053606000000939}
\end{barticle}
\endbibitem

\bibitem[\protect\citeauthoryear{Cox and Czanner}{2016}]{MR3513521}
\begin{barticle}[author]
\bauthor{\bsnm{Cox},~\bfnm{Trevor~F.}\binits{T.~F.}} \AND
  \bauthor{\bsnm{Czanner},~\bfnm{Gabriela}\binits{G.}}
(\byear{2016}).
\btitle{A practical divergence measure for survival distributions that can be
  estimated from {K}aplan-{M}eier curves}.
\bjournal{Stat. Med.}
\bvolume{35}
\bpages{2406--2421}.
\bdoi{10.1002/sim.6868}
\end{barticle}
\endbibitem

\bibitem[\protect\citeauthoryear{Cs\"{o}rg\H{o}}{1985}]{csorgHo1985}
\begin{barticle}[author]
\bauthor{\bsnm{Cs\"{o}rg\H{o}},~\bfnm{S\'{a}ndor}\binits{S.}}
(\byear{1985}).
\btitle{Testing for independence by the empirical characteristic function}.
\bjournal{J. Multivariate Anal.}
\bvolume{16}
\bpages{290--299}.
\bdoi{10.1016/0047-259X(85)90022-3}
\end{barticle}
\endbibitem

\bibitem[\protect\citeauthoryear{Cs\"{o}rg\H{o} and
  Horv\'{a}th}{1988}]{csorgHo1988}
\begin{barticle}[author]
\bauthor{\bsnm{Cs\"{o}rg\H{o}},~\bfnm{Mikl\'{o}s}\binits{M.}} \AND
  \bauthor{\bsnm{Horv\'{a}th},~\bfnm{Lajos}\binits{L.}}
(\byear{1988}).
\btitle{Central limit theorems for {$L_p$}-norms of density estimators}.
\bjournal{Probab. Theory Related Fields}
\bvolume{80}
\bpages{269--291}.
\bdoi{10.1007/BF00356106}
\end{barticle}
\endbibitem

\bibitem[\protect\citeauthoryear{Deb, Ghosal and Sen}{2020}]{deb2020measuring}
\begin{barticle}[author]
\bauthor{\bsnm{Deb},~\bfnm{Nabarun}\binits{N.}},
  \bauthor{\bsnm{Ghosal},~\bfnm{Promit}\binits{P.}} \AND
  \bauthor{\bsnm{Sen},~\bfnm{Bodhisattva}\binits{B.}}
(\byear{2020}).
\btitle{Measuring Association on Topological Spaces Using Kernels and Geometric
  Graphs. arXiv:2010.01768}.
\end{barticle}
\endbibitem

\bibitem[\protect\citeauthoryear{Deb and Sen}{2023}]{MR4571116}
\begin{barticle}[author]
\bauthor{\bsnm{Deb},~\bfnm{Nabarun}\binits{N.}} \AND
  \bauthor{\bsnm{Sen},~\bfnm{Bodhisattva}\binits{B.}}
(\byear{2023}).
\btitle{Multivariate rank-based distribution-free nonparametric testing using
  measure transportation}.
\bjournal{J. Amer. Statist. Assoc.}
\bvolume{118}
\bpages{192--207}.
\bdoi{10.1080/01621459.2021.1923508}
\end{barticle}
\endbibitem

\bibitem[\protect\citeauthoryear{Deheuvels}{1977}]{MR539513}
\begin{barticle}[author]
\bauthor{\bsnm{Deheuvels},~\bfnm{Paul}\binits{P.}}
(\byear{1977}).
\btitle{Estimation non param\'{e}trique de la densit\'{e} par histogrammes
  g\'{e}n\'{e}ralis\'{e}s. {II}}.
\bjournal{Publ. Inst. Statist. Univ. Paris}
\bvolume{22}
\bpages{1--23}.
\end{barticle}
\endbibitem

\bibitem[\protect\citeauthoryear{Deheuvels}{1981}]{Deheuvels1981}
\begin{barticle}[author]
\bauthor{\bsnm{Deheuvels},~\bfnm{Paul}\binits{P.}}
(\byear{1981}).
\btitle{A {K}olmogorov-{S}mirnov type test for independence and multivariate
  samples}.
\bjournal{Rev. Roumaine Math. Pures Appl.}
\bvolume{26}
\bpages{213--226}.
\end{barticle}
\endbibitem

\bibitem[\protect\citeauthoryear{Devroye}{1987}]{Devroye12}
\begin{bbook}[author]
\bauthor{\bsnm{Devroye},~\bfnm{Luc}\binits{L.}}
(\byear{1987}).
\btitle{A {C}ourse in {D}ensity {E}stimation}.
\bseries{Progress in Probability and Statistics}
\bvolume{14}.
\bpublisher{Birkh\"{a}user Boston, Inc., Boston, MA}.
\bdoi{10.1002/bimj.4710300618}
\end{bbook}
\endbibitem

\bibitem[\protect\citeauthoryear{Devroye}{1989}]{MR1045250}
\begin{barticle}[author]
\bauthor{\bsnm{Devroye},~\bfnm{Luc}\binits{L.}}
(\byear{1989}).
\btitle{The double kernel method in density estimation}.
\bjournal{Ann. Inst. H. Poincar\'{e} Probab. Statist.}
\bvolume{25}
\bpages{533--580}.
\end{barticle}
\endbibitem

\bibitem[\protect\citeauthoryear{Devroye}{1997}]{devroye1997}
\begin{barticle}[author]
\bauthor{\bsnm{Devroye},~\bfnm{Luc}\binits{L.}}
(\byear{1997}).
\btitle{Universal smoothing factor selection in density estimation: theory and
  practice}.
\bjournal{Test}
\bvolume{6}
\bpages{223--320}.
\bnote{With discussion and a rejoinder by the author}.
\bdoi{10.1007/BF02564701}
\end{barticle}
\endbibitem

\bibitem[\protect\citeauthoryear{Devroye and Gy\"{o}rfi}{1985}]{devroye1}
\begin{bbook}[author]
\bauthor{\bsnm{Devroye},~\bfnm{Luc}\binits{L.}} \AND
  \bauthor{\bsnm{Gy\"{o}rfi},~\bfnm{L\'{a}szl\'{o}}\binits{L.}}
(\byear{1985}).
\btitle{Nonparametric density estimation}.
\bseries{Wiley Series in Probability and Mathematical Statistics: Tracts on
  Probability and Statistics}.
\bpublisher{John Wiley \& Sons, Inc., New York}
\bnote{The $L_1$ view}.
\end{bbook}
\endbibitem

\bibitem[\protect\citeauthoryear{Devroye, Gy\"{o}rfi and
  Lugosi}{2002}]{MR1930019}
\begin{barticle}[author]
\bauthor{\bsnm{Devroye},~\bfnm{Luc}\binits{L.}},
  \bauthor{\bsnm{Gy\"{o}rfi},~\bfnm{L\'{a}szl\'{o}}\binits{L.}} \AND
  \bauthor{\bsnm{Lugosi},~\bfnm{G\'{a}bor}\binits{G.}}
(\byear{2002}).
\btitle{A note on robust hypothesis testing}.
\bjournal{IEEE Trans. Inform. Theory}
\bvolume{48}
\bpages{2111--2114}.
\bdoi{10.1109/TIT.2002.1013154}
\end{barticle}
\endbibitem

\bibitem[\protect\citeauthoryear{Devroye and Lugosi}{2001}]{devroye2001}
\begin{bbook}[author]
\bauthor{\bsnm{Devroye},~\bfnm{Luc}\binits{L.}} \AND
  \bauthor{\bsnm{Lugosi},~\bfnm{G\'{a}bor}\binits{G.}}
(\byear{2001}).
\btitle{Combinatorial {M}ethods in {D}ensity {E}stimation}.
\bseries{Springer Series in Statistics}.
\bpublisher{Springer-Verlag, New York}.
\bdoi{10.1007/978-1-4613-0125-7}
\end{bbook}
\endbibitem

\bibitem[\protect\citeauthoryear{Drton, Han and Shi}{2020}]{MR4185806}
\begin{barticle}[author]
\bauthor{\bsnm{Drton},~\bfnm{Mathias}\binits{M.}},
  \bauthor{\bsnm{Han},~\bfnm{Fang}\binits{F.}} \AND
  \bauthor{\bsnm{Shi},~\bfnm{Hongjian}\binits{H.}}
(\byear{2020}).
\btitle{High-dimensional consistent independence testing with maxima of rank
  correlations}.
\bjournal{Ann. Statist.}
\bvolume{48}
\bpages{3206--3227}.
\bdoi{10.1214/19-AOS1926}
\end{barticle}
\endbibitem

\bibitem[\protect\citeauthoryear{Dudley}{1989}]{dudley}
\begin{bbook}[author]
\bauthor{\bsnm{Dudley},~\bfnm{Richard~M.}\binits{R.~M.}}
(\byear{1989}).
\btitle{Real {A}nalysis and {P}robability}.
\bseries{The Wadsworth \& Brooks/Cole Mathematics Series}.
\bpublisher{Wadsworth \& Brooks/Cole Advanced Books \& Software},
  \baddress{Pacific Grove, CA}.
\end{bbook}
\endbibitem

\bibitem[\protect\citeauthoryear{Dugu\'{e}}{1975}]{dugue1975}
\begin{barticle}[author]
\bauthor{\bsnm{Dugu\'{e}},~\bfnm{Daniel}\binits{D.}}
(\byear{1975}).
\btitle{Sur des tests d'ind\'{e}pendance ``ind\'{e}pendants de la loi''}.
\bjournal{C. R. Acad. Sci. Paris S\'{e}r. A-B}
\bvolume{281}
\bpages{Aii, A1103--A1104}.
\end{barticle}
\endbibitem

\bibitem[\protect\citeauthoryear{Duong and Hazelton}{2005}]{duong2005}
\begin{barticle}[author]
\bauthor{\bsnm{Duong},~\bfnm{Tarn}\binits{T.}} \AND
  \bauthor{\bsnm{Hazelton},~\bfnm{Martin~L.}\binits{M.~L.}}
(\byear{2005}).
\btitle{Cross-validation bandwidth matrices for multivariate kernel density
  estimation}.
\bjournal{Scand. J. Statist.}
\bvolume{32}
\bpages{485--506}.
\bdoi{10.1111/j.1467-9469.2005.00445.x}
\end{barticle}
\endbibitem

\bibitem[\protect\citeauthoryear{Eggermont and LaRiccia}{2001}]{Eggermont2001}
\begin{bbook}[author]
\bauthor{\bsnm{Eggermont},~\bfnm{P.~P.~B.}\binits{P.~P.~B.}} \AND
  \bauthor{\bsnm{LaRiccia},~\bfnm{V.~N.}\binits{V.~N.}}
(\byear{2001}).
\btitle{Maximum {P}enalized {L}ikelihood {E}stimation. {V}ol. {I}}.
\bseries{Springer Series in Statistics}.
\bpublisher{Springer-Verlag}, \baddress{New York}.
\bnote{Density estimation}.
\bdoi{10.1007/978-1-0716-1244-6}
\end{bbook}
\endbibitem

\bibitem[\protect\citeauthoryear{Epsne\v{c}nikov}{1969}]{MR0250422}
\begin{barticle}[author]
\bauthor{\bsnm{Epsne\v{c}nikov},~\bfnm{V.~A.}\binits{V.~A.}}
(\byear{1969}).
\btitle{Nonparametric estimation of a multidimensional probability density}.
\bjournal{Teor. Verojatnost. i Primenen.}
\bvolume{14}
\bpages{156--162}.
\bdoi{10.1137/1114019}
\end{barticle}
\endbibitem

\bibitem[\protect\citeauthoryear{Fan et~al.}{2017}]{fan2017}
\begin{barticle}[author]
\bauthor{\bsnm{Fan},~\bfnm{Yanan}\binits{Y.}}, \bauthor{\bparticle{Lafaye~de}
  \bsnm{Micheaux},~\bfnm{Pierre}\binits{P.}},
  \bauthor{\bsnm{Penev},~\bfnm{Spiridon}\binits{S.}} \AND
  \bauthor{\bsnm{Salopek},~\bfnm{Donna}\binits{D.}}
(\byear{2017}).
\btitle{Multivariate nonparametric test of independence}.
\bjournal{J. Multivariate Anal.}
\bvolume{153}
\bpages{189--210}.
\bdoi{10.1016/j.jmva.2016.09.014}
\end{barticle}
\endbibitem

\bibitem[\protect\citeauthoryear{Felber, Kohler and
  Krzy{z}ak}{2015}]{MR3352508}
\begin{barticle}[author]
\bauthor{\bsnm{Felber},~\bfnm{Tina}\binits{T.}},
  \bauthor{\bsnm{Kohler},~\bfnm{Michael}\binits{M.}} \AND
  \bauthor{\bsnm{Krzy{z}ak},~\bfnm{Adam}\binits{A.}}
(\byear{2015}).
\btitle{Adaptive density estimation from data with small measurement errors}.
\bjournal{IEEE Trans. Inform. Theory}
\bvolume{61}
\bpages{3446--3456}.
\bdoi{10.1109/TIT.2015.2421297}
\end{barticle}
\endbibitem

\bibitem[\protect\citeauthoryear{Feuerverger}{1993}]{MR4205649}
\begin{barticle}[author]
\bauthor{\bsnm{Feuerverger},~\bfnm{Andrey}\binits{A.}}
(\byear{1993}).
\btitle{A consistent test for bivariate dependence}.
\bjournal{Int. Stat. Rev.}
\bvolume{61}
\bpages{419--433}.
\bdoi{10.2307/1403753}
\end{barticle}
\endbibitem

\bibitem[\protect\citeauthoryear{Fukumizu, Bach and Jordan}{2004}]{MR2247974}
\begin{barticle}[author]
\bauthor{\bsnm{Fukumizu},~\bfnm{Kenji}\binits{K.}},
  \bauthor{\bsnm{Bach},~\bfnm{Francis~R.}\binits{F.~R.}} \AND
  \bauthor{\bsnm{Jordan},~\bfnm{Michael~I.}\binits{M.~I.}}
(\byear{2004}).
\btitle{Dimensionality reduction for supervised learning with reproducing
  kernel {H}ilbert spaces}.
\bjournal{J. Mach. Learn. Res.}
\bvolume{5}
\bpages{73--99}.
\bdoi{10.1162/153244303768966111}
\end{barticle}
\endbibitem

\bibitem[\protect\citeauthoryear{Gan, Narisetty and Liang}{2019}]{gan}
\begin{barticle}[author]
\bauthor{\bsnm{Gan},~\bfnm{Lingrui}\binits{L.}},
  \bauthor{\bsnm{Narisetty},~\bfnm{Naveen~N.}\binits{N.~N.}} \AND
  \bauthor{\bsnm{Liang},~\bfnm{Feng}\binits{F.}}
(\byear{2019}).
\btitle{Bayesian regularization for graphical models with unequal shrinkage}.
\bjournal{J. Amer. Statist. Assoc.}
\bvolume{114}
\bpages{1218--1231}.
\bdoi{10.1080/01621459.2018.1482755}
\end{barticle}
\endbibitem

\bibitem[\protect\citeauthoryear{Gao and Gijbels}{2008}]{MR2504206}
\begin{barticle}[author]
\bauthor{\bsnm{Gao},~\bfnm{Jiti}\binits{J.}} \AND
  \bauthor{\bsnm{Gijbels},~\bfnm{Ir\`ene}\binits{I.}}
(\byear{2008}).
\btitle{Bandwidth selection in nonparametric kernel testing}.
\bjournal{J. Amer. Statist. Assoc.}
\bvolume{103}
\bpages{1584--1594}.
\bdoi{10.1198/016214508000000968}
\end{barticle}
\endbibitem

\bibitem[\protect\citeauthoryear{Gin{\'e}, Mason and Zaitsev}{2003}]{mason2001}
\begin{barticle}[author]
\bauthor{\bsnm{Gin{\'e}},~\bfnm{Evarist}\binits{E.}},
  \bauthor{\bsnm{Mason},~\bfnm{David~M.}\binits{D.~M.}} \AND
  \bauthor{\bsnm{Zaitsev},~\bfnm{Andrei~Yu.}\binits{A.~Y.}}
(\byear{2003}).
\btitle{The {$L\sb 1$}-norm density estimator process}.
\bjournal{Ann. Probab.}
\bvolume{31}
\bpages{719--768}.
\bdoi{10.1214/aop/1048516534}
\end{barticle}
\endbibitem

\bibitem[\protect\citeauthoryear{Gorsky and Ma}{2022}]{MR4472834}
\begin{barticle}[author]
\bauthor{\bsnm{Gorsky},~\bfnm{S.}\binits{S.}} \AND
  \bauthor{\bsnm{Ma},~\bfnm{L.}\binits{L.}}
(\byear{2022}).
\btitle{Multi-scale {F}isher's independence test for multivariate dependence}.
\bjournal{Biometrika}
\bvolume{109}
\bpages{569--587}.
\bdoi{10.1093/biomet/asac013}
\end{barticle}
\endbibitem

\bibitem[\protect\citeauthoryear{Gretton and Gy{\"o}rfi}{2010}]{gretton2010}
\begin{barticle}[author]
\bauthor{\bsnm{Gretton},~\bfnm{Arthur}\binits{A.}} \AND
  \bauthor{\bsnm{Gy{\"o}rfi},~\bfnm{L{\'a}szl{\'o}}\binits{L.}}
(\byear{2010}).
\btitle{Consistent nonparametric tests of independence}.
\bjournal{J. Mach. Learn. Res.}
\bvolume{11}
\bpages{1391--1423}.
\end{barticle}
\endbibitem

\bibitem[\protect\citeauthoryear{Gretton et~al.}{2005}]{gretton2005}
\begin{binproceedings}[author]
\bauthor{\bsnm{Gretton},~\bfnm{Arthur}\binits{A.}},
  \bauthor{\bsnm{Bousquet},~\bfnm{Olivier}\binits{O.}},
  \bauthor{\bsnm{Smola},~\bfnm{Alex}\binits{A.}} \AND
  \bauthor{\bsnm{Sch\"{o}lkopf},~\bfnm{Bernhard}\binits{B.}}
(\byear{2005}).
\btitle{Measuring statistical dependence with {H}ilbert-{S}chmidt norms}.
In \bbooktitle{Algorithmic learning theory}.
\bseries{Lecture Notes in Comput. Sci.}
\bvolume{3734}
\bpages{63--77}.
\bpublisher{Springer, Berlin}.
\bdoi{10.1007/11564089\_7}
\end{binproceedings}
\endbibitem

\bibitem[\protect\citeauthoryear{Hahn}{2005}]{HAHN200578}
\begin{barticle}[author]
\bauthor{\bsnm{Hahn},~\bfnm{T.}\binits{T.}}
(\byear{2005}).
\btitle{Cuba---a library for multidimensional numerical integration}.
\bjournal{Comput. Phys. Comm.}
\bvolume{168}
\bpages{78--95}.
\bdoi{10.1016/j.cpc.2005.01.010}
\end{barticle}
\endbibitem

\bibitem[\protect\citeauthoryear{Hall and Wand}{1988}]{Hall1992}
\begin{barticle}[author]
\bauthor{\bsnm{Hall},~\bfnm{Peter}\binits{P.}} \AND
  \bauthor{\bsnm{Wand},~\bfnm{Matthew~P.}\binits{M.~P.}}
(\byear{1988}).
\btitle{Minimizing {$L_1$} distance in nonparametric density estimation}.
\bjournal{J. Multivariate Anal.}
\bvolume{26}
\bpages{59--88}.
\bdoi{10.1016/0047-259X(88)90073-5}
\end{barticle}
\endbibitem

\bibitem[\protect\citeauthoryear{Han, Chen and Liu}{2017}]{han2017}
\begin{barticle}[author]
\bauthor{\bsnm{Han},~\bfnm{Fang}\binits{F.}},
  \bauthor{\bsnm{Chen},~\bfnm{Shizhe}\binits{S.}} \AND
  \bauthor{\bsnm{Liu},~\bfnm{Han}\binits{H.}}
(\byear{2017}).
\btitle{Distribution-free tests of independence in high dimensions}.
\bjournal{Biometrika}
\bvolume{104}
\bpages{813--828}.
\bdoi{10.1093/biomet/asx050}
\end{barticle}
\endbibitem

\bibitem[\protect\citeauthoryear{H\"{a}rdle}{1990}]{hardle1990}
\begin{bbook}[author]
\bauthor{\bsnm{H\"{a}rdle},~\bfnm{Wolfgang}\binits{W.}}
(\byear{1990}).
\btitle{Applied {N}onparametric {R}egression}.
\bseries{Econometric Society Monographs}
\bvolume{19}.
\bpublisher{Cambridge University Press, Cambridge}.
\bdoi{10.1017/CCOL0521382483}
\end{bbook}
\endbibitem

\bibitem[\protect\citeauthoryear{Heller, Heller and Gorfine}{2013}]{heller2012}
\begin{barticle}[author]
\bauthor{\bsnm{Heller},~\bfnm{Ruth}\binits{R.}},
  \bauthor{\bsnm{Heller},~\bfnm{Yair}\binits{Y.}} \AND
  \bauthor{\bsnm{Gorfine},~\bfnm{Malka}\binits{M.}}
(\byear{2013}).
\btitle{A consistent multivariate test of association based on ranks of
  distances}.
\bjournal{Biometrika}
\bvolume{100}
\bpages{503--510}.
\bdoi{10.1093/biomet/ass070}
\end{barticle}
\endbibitem

\bibitem[\protect\citeauthoryear{Hoeffding}{1948}]{Hoeffding}
\begin{barticle}[author]
\bauthor{\bsnm{Hoeffding},~\bfnm{Wassily}\binits{W.}}
(\byear{1948}).
\btitle{A non-parametric test of independence}.
\bjournal{Ann. Math. Statist.}
\bvolume{19}
\bpages{546--557}.
\bdoi{10.1214/aoms/1177730150}
\end{barticle}
\endbibitem

\bibitem[\protect\citeauthoryear{Horowitz and Spokoiny}{2001}]{MR1828537}
\begin{barticle}[author]
\bauthor{\bsnm{Horowitz},~\bfnm{Joel~L.}\binits{J.~L.}} \AND
  \bauthor{\bsnm{Spokoiny},~\bfnm{Vladimir~G.}\binits{V.~G.}}
(\byear{2001}).
\btitle{An adaptive, rate-optimal test of a parametric mean-regression model
  against a nonparametric alternative}.
\bjournal{Econometrica}
\bvolume{69}
\bpages{599--631}.
\bdoi{10.1111/1468-0262.00207}
\end{barticle}
\endbibitem

\bibitem[\protect\citeauthoryear{Horv\'{a}th}{1991}]{horvath1991}
\begin{barticle}[author]
\bauthor{\bsnm{Horv\'{a}th},~\bfnm{Lajos}\binits{L.}}
(\byear{1991}).
\btitle{On {$L_p$}-norms of multivariate density estimators}.
\bjournal{Ann. Statist.}
\bvolume{19}
\bpages{1933--1949}.
\bdoi{10.1214/aos/1176348379}
\end{barticle}
\endbibitem

\bibitem[\protect\citeauthoryear{Jin and Matteson}{2018}]{MR3858367}
\begin{barticle}[author]
\bauthor{\bsnm{Jin},~\bfnm{Ze}\binits{Z.}} \AND
  \bauthor{\bsnm{Matteson},~\bfnm{David~S.}\binits{D.~S.}}
(\byear{2018}).
\btitle{Generalizing distance covariance to measure and test multivariate
  mutual dependence via complete and incomplete {V}-statistics}.
\bjournal{J. Multivariate Anal.}
\bvolume{168}
\bpages{304--322}.
\bdoi{10.1016/j.jmva.2018.08.006}
\end{barticle}
\endbibitem

\bibitem[\protect\citeauthoryear{Kendall}{1938}]{kendall1938}
\begin{barticle}[author]
\bauthor{\bsnm{Kendall},~\bfnm{Maurice~G}\binits{M.~G.}}
(\byear{1938}).
\btitle{A new measure of rank correlation}.
\bjournal{Biometrika}
\bvolume{30}
\bpages{81--93}.
\bdoi{10.1093/biomet/30.1-2.81}
\end{barticle}
\endbibitem

\bibitem[\protect\citeauthoryear{Lauritzen}{1996}]{Lauritzen}
\begin{bbook}[author]
\bauthor{\bsnm{Lauritzen},~\bfnm{Steffen~L.}\binits{S.~L.}}
(\byear{1996}).
\btitle{Graphical {M}odels}.
\bseries{Oxford Statistical Science Series}
\bvolume{17}.
\bpublisher{The Clarendon Press, Oxford University Press, New York}
\bnote{Oxford Science Publications}.
\bdoi{10.1093/oso/9780198522195.001.0001}
\end{bbook}
\endbibitem

\bibitem[\protect\citeauthoryear{Le~Cam}{1970}]{MR0410832}
\begin{bincollection}[author]
\bauthor{\bsnm{Le~Cam},~\bfnm{L.}\binits{L.}}
(\byear{1970}).
\btitle{Remarques sur le th\'{e}or\`eme limite central dans les espaces
  localement convexes}.
In \bbooktitle{Les probabilit\'{e}s sur les structures alg\'{e}briques ({A}ctes
  {C}olloq. {I}nternat. {CNRS}, {N}o. 186, {C}lermont-{F}errand, 1969)}
\bpages{233--249}.
\bpublisher{\'{E}ditions Centre Nat. Recherche Sci., Paris}
\bnote{Avec commentaire en anglais par R. M. Dudley}.
\end{bincollection}
\endbibitem

\bibitem[\protect\citeauthoryear{Lee, Song and Whang}{2013}]{lee2013}
\begin{barticle}[author]
\bauthor{\bsnm{Lee},~\bfnm{Sokbae}\binits{S.}},
  \bauthor{\bsnm{Song},~\bfnm{Kyungchul}\binits{K.}} \AND
  \bauthor{\bsnm{Whang},~\bfnm{Yoon-Jae}\binits{Y.-J.}}
(\byear{2013}).
\btitle{Testing functional inequalities}.
\bjournal{J. Econometrics}
\bvolume{172}
\bpages{14--32}.
\bdoi{10.1016/j.jeconom.2012.08.006}
\end{barticle}
\endbibitem

\bibitem[\protect\citeauthoryear{Leung and Drton}{2018}]{leung2018}
\begin{barticle}[author]
\bauthor{\bsnm{Leung},~\bfnm{Dennis}\binits{D.}} \AND
  \bauthor{\bsnm{Drton},~\bfnm{Mathias}\binits{M.}}
(\byear{2018}).
\btitle{Testing independence in high dimensions with sums of rank
  correlations}.
\bjournal{Ann. Statist.}
\bvolume{46}
\bpages{280--307}.
\bdoi{10.1214/17-AOS1550}
\end{barticle}
\endbibitem

\bibitem[\protect\citeauthoryear{Matteson and Tsay}{2017}]{matt}
\begin{barticle}[author]
\bauthor{\bsnm{Matteson},~\bfnm{David~S.}\binits{D.~S.}} \AND
  \bauthor{\bsnm{Tsay},~\bfnm{Ruey~S.}\binits{R.~S.}}
(\byear{2017}).
\btitle{Independent component analysis via distance covariance}.
\bjournal{J. Amer. Statist. Assoc.}
\bvolume{112}
\bpages{623--637}.
\bdoi{10.1080/01621459.2016.1150851}
\end{barticle}
\endbibitem

\bibitem[\protect\citeauthoryear{Moon}{2016}]{moon2016nonparametric}
\begin{bbook}[author]
\bauthor{\bsnm{Moon},~\bfnm{Kevin~R}\binits{K.~R.}}
(\byear{2016}).
\btitle{Nonparametric Estimation of Distributional Functionals and
  Applications}.
\bnote{Thesis (Ph.D.)--University of Michigan}.
\end{bbook}
\endbibitem

\bibitem[\protect\citeauthoryear{Nabeya}{1951}]{nabeya}
\begin{barticle}[author]
\bauthor{\bsnm{Nabeya},~\bfnm{Seiji}\binits{S.}}
(\byear{1951}).
\btitle{Absolute moments in {$2$}-dimensional normal distribution}.
\bjournal{Ann. Inst. Statist. Math., Tokyo}
\bvolume{3}
\bpages{2--6}.
\bdoi{10.1007/bf02949770}
\end{barticle}
\endbibitem

\bibitem[\protect\citeauthoryear{Nadaraya}{1989}]{nadaraya1989}
\begin{bbook}[author]
\bauthor{\bsnm{Nadaraya},~\bfnm{\`E.~A.}\binits{E.~A.}}
(\byear{1989}).
\btitle{Nonparametric {E}stimation of {P}robability {D}ensities and
  {R}egression {C}urves}.
\bseries{Mathematics and its Applications (Soviet Series)}
\bvolume{20}.
\bpublisher{Kluwer Academic Publishers Group, Dordrecht}.
\bdoi{10.1007/978-94-009-2583-0}
\end{bbook}
\endbibitem

\bibitem[\protect\citeauthoryear{Nguyen and Eisenstein}{2017}]{nguyen}
\begin{barticle}[author]
\bauthor{\bsnm{Nguyen},~\bfnm{Dong}\binits{D.}} \AND
  \bauthor{\bsnm{Eisenstein},~\bfnm{Jacob}\binits{J.}}
(\byear{2017}).
\btitle{A kernel independence test for geographical language variation}.
\bjournal{Computational Linguistics}
\bvolume{43}
\bpages{567--592}.
\bdoi{10.1162/COLI\_a\_00293}
\end{barticle}
\endbibitem

\bibitem[\protect\citeauthoryear{Parzen}{1962}]{Parzen1962}
\begin{barticle}[author]
\bauthor{\bsnm{Parzen},~\bfnm{Emanuel}\binits{E.}}
(\byear{1962}).
\btitle{On estimation of a probability density function and mode}.
\bjournal{Ann. Math. Statist.}
\bvolume{33}
\bpages{1065--1076}.
\bdoi{10.1214/aoms/1177704472}
\end{barticle}
\endbibitem

\bibitem[\protect\citeauthoryear{Pearl}{2009}]{pearl}
\begin{bbook}[author]
\bauthor{\bsnm{Pearl},~\bfnm{Judea}\binits{J.}}
(\byear{2009}).
\btitle{Causality},
\bedition{Second} ed.
\bpublisher{Cambridge University Press, Cambridge}
\bnote{Models, reasoning, and inference}.
\bdoi{10.1017/CBO9780511803161}
\end{bbook}
\endbibitem

\bibitem[\protect\citeauthoryear{Pearson}{1920}]{PEARSONKARL}
\begin{barticle}[author]
\bauthor{\bsnm{Pearson},~\bfnm{Karl}\binits{K.}}
(\byear{1920}).
\btitle{{Notes on the history of correlation}}.
\bjournal{Biometrika}
\bvolume{13}
\bpages{25--45}.
\bdoi{10.2307/2331722}
\end{barticle}
\endbibitem

\bibitem[\protect\citeauthoryear{Peters et~al.}{2014}]{peters}
\begin{barticle}[author]
\bauthor{\bsnm{Peters},~\bfnm{Jonas}\binits{J.}},
  \bauthor{\bsnm{Mooij},~\bfnm{Joris~M.}\binits{J.~M.}},
  \bauthor{\bsnm{Janzing},~\bfnm{Dominik}\binits{D.}} \AND
  \bauthor{\bsnm{Sch\"{o}lkopf},~\bfnm{Bernhard}\binits{B.}}
(\byear{2014}).
\btitle{Causal discovery with continuous additive noise models}.
\bjournal{J. Mach. Learn. Res.}
\bvolume{15}
\bpages{2009--2053}.
\bdoi{10.15496/publikation-1672}
\end{barticle}
\endbibitem

\bibitem[\protect\citeauthoryear{Pfister and Peters}{2017}]{Pfister2017}
\begin{barticle}[author]
\bauthor{\bsnm{Pfister},~\bfnm{Niklas}\binits{N.}} \AND
  \bauthor{\bsnm{Peters},~\bfnm{J.}\binits{J.}}
(\byear{2017}).
\btitle{d\textnormal{HSIC:} {\it{Independence Testing via Hilbert-Schmidt
  Independence Criterion.}}}
\bjournal{\textnormal{R Package version 2.1. Available at
  \url{https://cran.r-project.org/web/packages/dHSIC/index.html.}}}
\end{barticle}
\endbibitem

\bibitem[\protect\citeauthoryear{Pfister et~al.}{2018}]{pfister2018}
\begin{barticle}[author]
\bauthor{\bsnm{Pfister},~\bfnm{Niklas}\binits{N.}},
  \bauthor{\bsnm{B\"{u}hlmann},~\bfnm{Peter}\binits{P.}},
  \bauthor{\bsnm{Sch\"{o}lkopf},~\bfnm{Bernhard}\binits{B.}} \AND
  \bauthor{\bsnm{Peters},~\bfnm{Jonas}\binits{J.}}
(\byear{2018}).
\btitle{Kernel-based tests for joint independence}.
\bjournal{J. R. Stat. Soc. Ser. B. Stat. Methodol.}
\bvolume{80}
\bpages{5--31}.
\bdoi{10.1111/rssb.12235}
\end{barticle}
\endbibitem

\bibitem[\protect\citeauthoryear{Pinelis}{1994}]{pinelis}
\begin{barticle}[author]
\bauthor{\bsnm{Pinelis},~\bfnm{Iosif}\binits{I.}}
(\byear{1994}).
\btitle{On a majorization inequality for sums of independent random vectors}.
\bjournal{Statist. Probab. Lett.}
\bvolume{19}
\bpages{97--99}.
\bdoi{10.1016/0167-7152(94)90139-2}
\end{barticle}
\endbibitem

\bibitem[\protect\citeauthoryear{Qiu, Xu and Zhu}{2023}]{MR4544605}
\begin{barticle}[author]
\bauthor{\bsnm{Qiu},~\bfnm{Tao}\binits{T.}},
  \bauthor{\bsnm{Xu},~\bfnm{Wangli}\binits{W.}} \AND
  \bauthor{\bsnm{Zhu},~\bfnm{Lixing}\binits{L.}}
(\byear{2023}).
\btitle{Independence tests with random subspace of two random vectors in high
  dimension}.
\bjournal{J. Multivariate Anal.}
\bvolume{195}
\bpages{Paper No. 105160}.
\bdoi{10.1016/j.jmva.2023.105160}
\end{barticle}
\endbibitem

\bibitem[\protect\citeauthoryear{Rizzo}{2009}]{MR2751846}
\begin{barticle}[author]
\bauthor{\bsnm{Rizzo},~\bfnm{Maria~L.}\binits{M.~L.}}
(\byear{2009}).
\btitle{New goodness-of-fit tests for {P}areto distributions}.
\bjournal{Astin Bull.}
\bvolume{39}
\bpages{691--715}.
\bdoi{10.2143/AST.39.2.2044654}
\end{barticle}
\endbibitem

\bibitem[\protect\citeauthoryear{Rizzo and Sz{\'e}kely}{2016}]{rizzo2016}
\begin{barticle}[author]
\bauthor{\bsnm{Rizzo},~\bfnm{ML}\binits{M.}} \AND
  \bauthor{\bsnm{Sz{\'e}kely},~\bfnm{GJ}\binits{G.}}
(\byear{2016}).
\btitle{E-Statistics: Multivariate inference via the energy of data. R package
  version 1.7-0}.
\bjournal{\textnormal{Available at
  \url{https://cran.r-project.org/web/packages/energy/index.html}}}.
\end{barticle}
\endbibitem

\bibitem[\protect\citeauthoryear{Rosenblatt}{1956}]{Rosenblatt}
\begin{barticle}[author]
\bauthor{\bsnm{Rosenblatt},~\bfnm{Murray}\binits{M.}}
(\byear{1956}).
\btitle{Remarks on some nonparametric estimates of a density function}.
\bjournal{Ann. Math. Statist.}
\bvolume{27}
\bpages{832--837}.
\bdoi{10.1214/aoms/1177728190}
\end{barticle}
\endbibitem

\bibitem[\protect\citeauthoryear{Rosenblatt}{1975}]{rosenblatt1975}
\begin{barticle}[author]
\bauthor{\bsnm{Rosenblatt},~\bfnm{M.}\binits{M.}}
(\byear{1975}).
\btitle{A quadratic measure of deviation of two-dimensional density estimates
  and a test of independence}.
\bjournal{Ann. Statist.}
\bvolume{3}
\bpages{1--14}.
\bdoi{10.1214/aos/1176342996}
\end{barticle}
\endbibitem

\bibitem[\protect\citeauthoryear{Rosenblatt and Wahlen}{1992}]{MR1190260}
\begin{barticle}[author]
\bauthor{\bsnm{Rosenblatt},~\bfnm{Murray}\binits{M.}} \AND
  \bauthor{\bsnm{Wahlen},~\bfnm{Bruce~E.}\binits{B.~E.}}
(\byear{1992}).
\btitle{A nonparametric measure of independence under a hypothesis of
  independent components}.
\bjournal{Statist. Probab. Lett.}
\bvolume{15}
\bpages{245--252}.
\bdoi{10.1016/0167-7152(92)90197-D}
\end{barticle}
\endbibitem

\bibitem[\protect\citeauthoryear{Sain, Baggerly and Scott}{1994}]{sain1994}
\begin{barticle}[author]
\bauthor{\bsnm{Sain},~\bfnm{Stephan~R.}\binits{S.~R.}},
  \bauthor{\bsnm{Baggerly},~\bfnm{Keith~A.}\binits{K.~A.}} \AND
  \bauthor{\bsnm{Scott},~\bfnm{David~W.}\binits{D.~W.}}
(\byear{1994}).
\btitle{Cross-validation of multivariate densities}.
\bjournal{J. Amer. Statist. Assoc.}
\bvolume{89}
\bpages{807--817}.
\bdoi{10.1080/01621459.1994.10476814}
\end{barticle}
\endbibitem

\bibitem[\protect\citeauthoryear{Samworth and Yuan}{2012}]{Samworth}
\begin{barticle}[author]
\bauthor{\bsnm{Samworth},~\bfnm{Richard~J.}\binits{R.~J.}} \AND
  \bauthor{\bsnm{Yuan},~\bfnm{Ming}\binits{M.}}
(\byear{2012}).
\btitle{Independent component analysis via nonparametric maximum likelihood
  estimation}.
\bjournal{Ann. Statist.}
\bvolume{40}
\bpages{2973--3002}.
\bdoi{10.1214/12-AOS1060}
\end{barticle}
\endbibitem

\bibitem[\protect\citeauthoryear{Scott}{1992}]{Scot1992}
\begin{bbook}[author]
\bauthor{\bsnm{Scott},~\bfnm{David~W.}\binits{D.~W.}}
(\byear{1992}).
\btitle{Multivariate {D}ensity {E}stimation}.
\bseries{Wiley Series in Probability and Mathematical Statistics: Applied
  Probability and Statistics}.
\bpublisher{John Wiley \& Sons Inc.}, \baddress{New York}.
\bnote{Theory, practice, and visualization, A Wiley-Interscience Publication}.
\bdoi{10.1002/9780470316849}
\end{bbook}
\endbibitem

\bibitem[\protect\citeauthoryear{Scott and Wand}{1991}]{MR1118245}
\begin{barticle}[author]
\bauthor{\bsnm{Scott},~\bfnm{David~W.}\binits{D.~W.}} \AND
  \bauthor{\bsnm{Wand},~\bfnm{M.~P.}\binits{M.~P.}}
(\byear{1991}).
\btitle{Feasibility of multivariate density estimates}.
\bjournal{Biometrika}
\bvolume{78}
\bpages{197--205}.
\bdoi{10.1093/biomet/78.1.197}
\end{barticle}
\endbibitem

\bibitem[\protect\citeauthoryear{Sejdinovic et~al.}{2013}]{sejdinovic2013}
\begin{barticle}[author]
\bauthor{\bsnm{Sejdinovic},~\bfnm{Dino}\binits{D.}},
  \bauthor{\bsnm{Sriperumbudur},~\bfnm{Bharath}\binits{B.}},
  \bauthor{\bsnm{Gretton},~\bfnm{Arthur}\binits{A.}} \AND
  \bauthor{\bsnm{Fukumizu},~\bfnm{Kenji}\binits{K.}}
(\byear{2013}).
\btitle{Equivalence of distance-based and {RKHS}-based statistics in hypothesis
  testing}.
\bjournal{Ann. Statist.}
\bvolume{41}
\bpages{2263--2291}.
\bdoi{10.1214/13-AOS1140}
\end{barticle}
\endbibitem

\bibitem[\protect\citeauthoryear{Sheng and Yin}{2016}]{MR3474038}
\begin{barticle}[author]
\bauthor{\bsnm{Sheng},~\bfnm{Wenhui}\binits{W.}} \AND
  \bauthor{\bsnm{Yin},~\bfnm{Xiangrong}\binits{X.}}
(\byear{2016}).
\btitle{Sufficient dimension reduction via distance covariance}.
\bjournal{J. Comput. Graph. Statist.}
\bvolume{25}
\bpages{91--104}.
\bdoi{10.1080/10618600.2015.1026601}
\end{barticle}
\endbibitem

\bibitem[\protect\citeauthoryear{Shergin}{1990}]{shergin}
\begin{bincollection}[author]
\bauthor{\bsnm{Shergin},~\bfnm{V.~V.}\binits{V.~V.}}
(\byear{1990}).
\btitle{The central limit theorem for finitely dependent random variables}.
In \bbooktitle{Probability theory and mathematical statistics, {V}ol. {II}
  ({V}ilnius, 1989)}
\bpages{424--431}.
\bpublisher{``Mokslas'', Vilnius}.
\bdoi{10.1007/BF01261286}
\end{bincollection}
\endbibitem

\bibitem[\protect\citeauthoryear{Shi, Drton and Han}{2022a}]{shi}
\begin{barticle}[author]
\bauthor{\bsnm{Shi},~\bfnm{Hongjian}\binits{H.}},
  \bauthor{\bsnm{Drton},~\bfnm{Mathias}\binits{M.}} \AND
  \bauthor{\bsnm{Han},~\bfnm{Fang}\binits{F.}}
(\byear{2022}a).
\btitle{Distribution-free consistent independence tests via center-outward
  ranks and signs}.
\bjournal{J. Amer. Statist. Assoc.}
\bvolume{117}
\bpages{395--410}.
\bdoi{10.1080/01621459.2020.1782223}
\end{barticle}
\endbibitem

\bibitem[\protect\citeauthoryear{Shi, Drton and Han}{2022b}]{shi2022power}
\begin{barticle}[author]
\bauthor{\bsnm{Shi},~\bfnm{H.}\binits{H.}},
  \bauthor{\bsnm{Drton},~\bfnm{M.}\binits{M.}} \AND
  \bauthor{\bsnm{Han},~\bfnm{F.}\binits{F.}}
(\byear{2022}b).
\btitle{On the power of {C}hatterjee's rank correlation}.
\bjournal{Biometrika}
\bvolume{109}
\bpages{317--333}.
\bdoi{10.1093/biomet/asab028}
\end{barticle}
\endbibitem

\bibitem[\protect\citeauthoryear{Silverman}{1986}]{silverman1986}
\begin{bbook}[author]
\bauthor{\bsnm{Silverman},~\bfnm{B.~W.}\binits{B.~W.}}
(\byear{1986}).
\btitle{Density {E}stimation for {S}tatistics and {D}ata {A}nalysis}.
\bseries{Monographs on Statistics and Applied Probability}.
\bpublisher{Chapman \& Hall}, \baddress{London}.
\bdoi{10.1007/978-1-4899-3324-9}
\end{bbook}
\endbibitem

\bibitem[\protect\citeauthoryear{Spearman}{1904}]{spearman1904}
\begin{barticle}[author]
\bauthor{\bsnm{Spearman},~\bfnm{Charles}\binits{C.}}
(\byear{1904}).
\btitle{The proof and measurement of association between two things}.
\bjournal{American journal of Psychology}
\bvolume{15}
\bpages{72--101}.
\bdoi{10.2307/1422689}
\end{barticle}
\endbibitem

\bibitem[\protect\citeauthoryear{Sriperumbudur et~al.}{2012}]{MR2988458}
\begin{barticle}[author]
\bauthor{\bsnm{Sriperumbudur},~\bfnm{Bharath~K.}\binits{B.~K.}},
  \bauthor{\bsnm{Fukumizu},~\bfnm{Kenji}\binits{K.}},
  \bauthor{\bsnm{Gretton},~\bfnm{Arthur}\binits{A.}},
  \bauthor{\bsnm{Sch\"{o}lkopf},~\bfnm{Bernhard}\binits{B.}} \AND
  \bauthor{\bsnm{Lanckriet},~\bfnm{Gert R.~G.}\binits{G.~R.~G.}}
(\byear{2012}).
\btitle{On the empirical estimation of integral probability metrics}.
\bjournal{Electron. J. Stat.}
\bvolume{6}
\bpages{1550--1599}.
\bdoi{10.1214/12-EJS722}
\end{barticle}
\endbibitem

\bibitem[\protect\citeauthoryear{Stein}{1970}]{stein}
\begin{bbook}[author]
\bauthor{\bsnm{Stein},~\bfnm{Elias~M.}\binits{E.~M.}}
(\byear{1970}).
\btitle{Singular {I}ntegrals and {D}ifferentiability {P}roperties of
  {F}unctions}.
\bseries{Princeton Mathematical Series, No. 30}.
\bpublisher{Princeton University Press, Princeton, N.J.}
\bdoi{10.1515/9781400883882}
\end{bbook}
\endbibitem

\bibitem[\protect\citeauthoryear{Sweeting}{1977}]{sweeting}
\begin{barticle}[author]
\bauthor{\bsnm{Sweeting},~\bfnm{T.~J.}\binits{T.~J.}}
(\byear{1977}).
\btitle{Speeds of convergence for the multidimensional central limit theorem}.
\bjournal{Ann. Probability}
\bvolume{5}
\bpages{28--41}.
\bdoi{10.1214/aop/1176995888}
\end{barticle}
\endbibitem

\bibitem[\protect\citeauthoryear{Sz\'{e}kely and Rizzo}{2005}]{MR2231170}
\begin{barticle}[author]
\bauthor{\bsnm{Sz\'{e}kely},~\bfnm{G\'{a}bor~J.}\binits{G.~J.}} \AND
  \bauthor{\bsnm{Rizzo},~\bfnm{Maria~L.}\binits{M.~L.}}
(\byear{2005}).
\btitle{Hierarchical clustering via joint between-within distances: extending
  {W}ard's minimum variance method}.
\bjournal{J. Classification}
\bvolume{22}
\bpages{151--183}.
\bdoi{10.1007/s00357-005-0012-9}
\end{barticle}
\endbibitem

\bibitem[\protect\citeauthoryear{Sz\'{e}kely, Rizzo and
  Bakirov}{2007}]{szekely2007}
\begin{barticle}[author]
\bauthor{\bsnm{Sz\'{e}kely},~\bfnm{G\'{a}bor~J.}\binits{G.~J.}},
  \bauthor{\bsnm{Rizzo},~\bfnm{Maria~L.}\binits{M.~L.}} \AND
  \bauthor{\bsnm{Bakirov},~\bfnm{Nail~K.}\binits{N.~K.}}
(\byear{2007}).
\btitle{Measuring and testing dependence by correlation of distances}.
\bjournal{Ann. Statist.}
\bvolume{35}
\bpages{2769--2794}.
\bdoi{10.1214/009053607000000505}
\end{barticle}
\endbibitem

\bibitem[\protect\citeauthoryear{Sz\'{e}kely and Rizzo}{2009}]{szekely2009}
\begin{barticle}[author]
\bauthor{\bsnm{Sz\'{e}kely},~\bfnm{G\'{a}bor~J.}\binits{G.~J.}} \AND
  \bauthor{\bsnm{Rizzo},~\bfnm{Maria~L.}\binits{M.~L.}}
(\byear{2009}).
\btitle{Brownian distance covariance}.
\bjournal{Ann. Appl. Stat.}
\bvolume{3}
\bpages{1236--1265}.
\bdoi{10.1214/09-AOAS312}
\end{barticle}
\endbibitem

\bibitem[\protect\citeauthoryear{Wand and Jones}{1995}]{wand1995}
\begin{bbook}[author]
\bauthor{\bsnm{Wand},~\bfnm{M.~P.}\binits{M.~P.}} \AND
  \bauthor{\bsnm{Jones},~\bfnm{M.~C.}\binits{M.~C.}}
(\byear{1995}).
\btitle{Kernel {S}moothing}.
\bseries{Monographs on Statistics and Applied Probability}
\bvolume{60}.
\bpublisher{Chapman and Hall Ltd.}, \baddress{London}.
\bdoi{10.1007/978-1-4899-4493-1}
\end{bbook}
\endbibitem

\bibitem[\protect\citeauthoryear{Watson and Leadbetter}{1963}]{MR148149}
\begin{barticle}[author]
\bauthor{\bsnm{Watson},~\bfnm{G.~S.}\binits{G.~S.}} \AND
  \bauthor{\bsnm{Leadbetter},~\bfnm{M.~R.}\binits{M.~R.}}
(\byear{1963}).
\btitle{On the estimation of the probability density. {I}}.
\bjournal{Ann. Math. Statist.}
\bvolume{34}
\bpages{480--491}.
\bdoi{10.1214/aoms/1177704159}
\end{barticle}
\endbibitem

\bibitem[\protect\citeauthoryear{Yang}{2012}]{MR3093953}
\begin{bbook}[author]
\bauthor{\bsnm{Yang},~\bfnm{Guangyuan}\binits{G.}}
(\byear{2012}).
\btitle{The energy goodness-of-fit test for univariate stable distributions}.
\bpublisher{ProQuest LLC, Ann Arbor, MI}
\bnote{Thesis (Ph.D.)--Bowling Green State University}.
\end{bbook}
\endbibitem

\bibitem[\protect\citeauthoryear{Yatracos}{1985}]{MR790571}
\begin{barticle}[author]
\bauthor{\bsnm{Yatracos},~\bfnm{Yannis~G.}\binits{Y.~G.}}
(\byear{1985}).
\btitle{Rates of convergence of minimum distance estimators and {K}olmogorov's
  entropy}.
\bjournal{Ann. Statist.}
\bvolume{13}
\bpages{768--774}.
\bdoi{10.1214/aos/1176349553}
\end{barticle}
\endbibitem

\bibitem[\protect\citeauthoryear{Zhang, Gao and Ng}{2023}]{MR4527082}
\begin{barticle}[author]
\bauthor{\bsnm{Zhang},~\bfnm{Wei}\binits{W.}},
  \bauthor{\bsnm{Gao},~\bfnm{Wei}\binits{W.}} \AND
  \bauthor{\bsnm{Ng},~\bfnm{Hon Keung~Tony}\binits{H.~K.~T.}}
(\byear{2023}).
\btitle{Multivariate tests of independence based on a new class of measures of
  independence in {R}eproducing {K}ernel {H}ilbert {S}pace}.
\bjournal{J. Multivariate Anal.}
\bvolume{195}
\bpages{Paper No. 105144}.
\bdoi{10.1016/j.jmva.2022.105144}
\end{barticle}
\endbibitem

\bibitem[\protect\citeauthoryear{Zhang et~al.}{2018}]{MR3741641}
\begin{barticle}[author]
\bauthor{\bsnm{Zhang},~\bfnm{Qinyi}\binits{Q.}},
  \bauthor{\bsnm{Filippi},~\bfnm{Sarah}\binits{S.}},
  \bauthor{\bsnm{Gretton},~\bfnm{Arthur}\binits{A.}} \AND
  \bauthor{\bsnm{Sejdinovic},~\bfnm{Dino}\binits{D.}}
(\byear{2018}).
\btitle{Large-scale kernel methods for independence testing}.
\bjournal{Stat. Comput.}
\bvolume{28}
\bpages{113--130}.
\bdoi{10.1007/s11222-016-9721-7}
\end{barticle}
\endbibitem

\end{thebibliography}

\end{document}